\documentclass[12pt]{amsart}

\usepackage{latexsym,amsthm,verbatim,ifthen,graphicx,color,mathrsfs,amsmath,mathtools,bm,bbm}
\usepackage[all]{xy}
\usepackage{color}
\usepackage[symbol]{footmisc}
\makeatletter

\def\subsection{\@startsection{subsection}{2}%
  \z@{.5\linespacing}{.5\linespacing}%
  {\normalfont\bfseries}}
\makeatother

\makeatletter
\renewcommand{\l@subsection}{\@tocline{2}{0pt}{2.5em}{3.5em}{}}
\makeatother

\makeatletter
\newlength{\abstractwidth}
\renewenvironment{abstract}{%
  \ifx\maketitle\relax
    \ClassWarningNoLine{amsart}{Abstract should precede \protect\maketitle}%
  \fi
  \global\setbox\abstractbox=\vbox\bgroup
  \normalfont
  \vspace*{1.5\baselineskip}%
  \begin{center}\textsc{Abstract}\end{center}%
  \vspace{0.5\baselineskip}%
  \begin{center}\begin{minipage}{\abstractwidth}
  \normalfont\normalsize
  \noindent
}{%
  \par\end{minipage}\end{center}%
  \egroup
  \ifx\@setabstract\relax \@setabstracta \fi
}
\makeatother

\makeatletter
\renewcommand{\subsubsection}{\@startsection{subsubsection}{3}%
  \z@{1.0\linespacing \@plus .7\linespacing}{.5\linespacing}%
  {\normalfont\itshape}}
\makeatother

\usepackage{color}

\usepackage{lscape}
 \usepackage[ansinew]{inputenc}
 \usepackage{multicol}
 \usepackage[all]{xy}\xyoption{poly}
 \usepackage{mathrsfs}
\usepackage[psamsfonts]{eucal}
\usepackage{amssymb, amsmath}
\usepackage{geometry}
\usepackage{enumerate}
\usepackage{float}

\newtheorem{theorem}{Theorem}[section]

\newtheorem{corollary}[theorem]{Corollary}
 \newtheorem{lemma}[theorem]{Lemma}
 \newtheorem{proposition}[theorem]{Proposition}
 \theoremstyle{definition}
 \newtheorem{definition}[theorem]{Definition}

 \newtheorem{rem}[theorem]{Remark}

\newtheorem*{theorem*}{Theorem}

\def\sus{\Sigma}
\def\suse{\mathbf{\Sigma}}
\def\suseb{\mathbf{\Sigma}_B}
\def\susel{\mathbf{\Sigma}_L}
\def\susb{\Sigma_B}
\def\lup{\Omega}
\def\lupb{\Omega_B}
\def\lupl{\Omega_L}

\def\susl{\Sigma_L}
\newcommand\lupe{\mathbf{\Omega}}
\newcommand\lupeb{\mathbf{\Omega}_B}
\newcommand\lupel{\mathbf{\Omega}_L}
\newcommand\End{\rm End}

\def\quic{{\mathscr C}}
\def\quid{{\mathscr D}}

\def\esp#1{{\quad\text{#1}\quad}}

\newcommand{\bz}{\mathbb Z}
\newcommand{\compul}{\widehat{UL}}

\def\bq{\mathbb{Q}}

\def\bs{\mathbb{S}}

\newcommand{\adj}{{\scriptscriptstyle\vee}}

\newcommand{\cgm}{\operatorname{{\bf cgm
}}}
\newcommand{\cdgml}{\operatorname{{\bf cdgm
}}_L}
\newcommand{\cdgmul}{\operatorname{{\bf cdgm
}}_{\widehat{UL}}}
\newcommand{\cdgmulpri}{\operatorname{{\bf cdgm
}}_{\widehat{UL'}}}
\newcommand{\cdgma}{\operatorname{{\bf cdgm
}}_A}
\newcommand{\cdgm}{\operatorname{\bf cdgm
}}

    \newcommand{\lasu}{{\mathfrak{L}}}
    \newcommand{\lasub}{\boldsymbol{\mathfrak{L}}}

     \newcommand{\spec}{{{\bf Sp}}}
        \newcommand{\specb}{{{\bf Sp}_B}}
         \newcommand{\specl}{{{\bf Sp}_L}}
        \newcommand{\specmodu}{{{\bf Sp}_{\compul}^0}}  
        \newcommand{\specmodul}{{{\bf Sp}_{\compul}}}

       \newcommand{\ext}{{\rm Ext}}

 \newcommand{\lib }{\mathbb{L}}

 \newcommand{\cyl}{\operatorname{\text{\rm Cyl}\,}}

  \newcommand{\Hom}{\operatorname{\text{\rm Hom}}}
    \newcommand{\Endog}{\operatorname{\mathcal{E}\text{\rm nd}}}
        \newcommand{\Endo}{\operatorname{\text{\rm End}}}

\newcommand{\Span}{\operatorname{{Span}}}

\newcommand{\catcdgl}{\operatorname{{\bf cdgl}}}
\newcommand{\cdgll}{\operatorname{{\bf cdgl}}_{\sslash L}}

\newcommand{\catss}{\operatorname{{\bf sset}}}

\newcommand{\catvect}{\operatorname{{\bf vect}}}

\newcommand{\Ho}{\operatorname{{\rm Ho}}}

 \newcommand{\map}{\operatorname{{\rm map}}}

    \newcommand{\id}{\operatorname{{\rm id}}}
        \newcommand{\ad}{\operatorname{{\rm ad}}}
        \newcommand{\retb}{\operatorname{\bf sset}_{\sslash B}}
         \newcommand{\reta}{\operatorname{\bf sset}_{\sslash A}}

  \newcommand{\pie}{ \pi^{{\rm st}}}

 \newcommand{\MC}{\operatorname{{\rm MC}}}

     \newcommand{\eva}{\operatorname{{\rm ev}}}
\newcommand{\mc}{{\MC}}

  \newcommand{\otimesc}{\,\widehat{\otimes}\,}

   \newcommand{\libc}{{\widehat\lib}}
      \newcommand{\libcb}{\boldsymbol{\widehat\lib}}

\def\ret#1{\operatorname{\bf sset}_{\sslash #1}}

\newcommand{\quicr}{\quic_{\sslash c}}

\newcommand{\quidr}{\quid_{\sslash d}}

\usepackage{hyperref}
\usepackage{physics}
\hypersetup{
    colorlinks=true,
    linkcolor=blue,
    filecolor=magenta,
    urlcolor=cyan,
    citecolor=blue,
}

\newcommand{\acento}{{\scriptscriptstyle \wedge}}

\newcommand{\calker}{\operatorname{\mathfrak{K}}}

\newtheorem{example}[theorem]{Example}				
\theoremstyle{remark}

\usepackage{stmaryrd}

\usepackage{ textcomp }
\newcommand{\cdgl}{\operatorname{\bf cdgl}}

\theoremstyle{boldbody}
\newtheorem{introthm}{Theorem}[section]

\newtheorem{introprop}[introthm]{Proposition}

\begin{document}

\title{A new approach to rational stable parametrized homotopy theory}

\author{Yves F\'elix, Aniceto Murillo and Alejandro Saiz}

\thanks{This work was partially supported by the Spanish Government under grant PID2023-149804NB-I00 and by the Junta de Andaluc{\'\i}a under grant ProyExcel-00827-PY21.}

\begin{abstract}
This work develops a comprehensive algebraic model for rational stable parametrized homotopy theory over arbitrary  base spaces. Building on the simplicial analogue of the foundational framework of May--Sigurdsson for parametrized spectra and on the homotopy theory of complete differential graded Lie algebras, we construct an explicit sequence of Quillen equivalences that identify the homotopy theory of rational spectra of retractive simplicial sets  with the purely algebraic framework of complete differential graded modules over the completed universal enveloping algebra $\compul$ of a Lie model $L$ of the base simplicial set $B$. More precisely, there is a sequence of Quillen adjunctions
$$
\xymatrix{ 
\spec_B 
  \ar@<0.75ex>[r] 
& 
\specl 
  \ar@<0.75ex>[l]
    \ar@<-0.80ex>[r]
& 
\specmodu 
  \ar@<-0.60ex>[l]
  \ar@<0.50ex>[r]
& 
\cdgmul, 
  \ar@<0.90ex>[l] 
}
$$
which induces a natural, strong monoidal equivalence of categories
$$
 \Ho\spec_B^\bq\cong\Ho\cdgmul.
  $$
This equivalence is highly effective in practice as it provides direct computational access to invariants of simplicial spectra by translating them into homotopy invariants of $\compul$-modules. Here $\specb$ is the stable model category of spectra of retractive simplicial sets over $B$, $\specl$ denotes the stable model category of spectra of retractive complete differential graded Lie algebras over $L$, $\specmodu$ is the stable model category of connected $\compul$-module
 spectra, and $\cdgmul$ denotes the category of complete differential graded $\compul$-modules.
\end{abstract}

\maketitle

\tableofcontents

\section*{Introduction}

The development of a stable homotopy theory for parametrized spaces has undergone substantial progress over the past two decades. Building on foundational work by May and Sigurdsson \cite{maysi}, parametrized spectra have emerged as the correct homotopical generalization of classical bundles of spectra, allowing for a coherent treatment of fiberwise constructions and base change. Such developments have led, for instance, to a rich theory of twisted generalized (co)homology, enriched categories of spectra, and structured module categories over ring spectra.

Our focus is on the rational aspects of parametrized spectra and on their computability by means of faithful  models of  their algebraic counterparts. We prove:

\begin{introthm}\label{A} Let $B$ be a  simplicial set with Lie model $L$. Then there is a natural and strong symmetric monoidal equivalence between the homotopy category $\Ho\spec_B^\bq$ of rational spectra of retractive simplicial sets over $B$ and the homotopy category $\Ho\cdgmul$ of complete modules over the completed universal enveloping algebra of $L$. 
\end{introthm}

This theorem extends the main results in \cite{brau2} to the non-simply connected setting, and can also be viewed as  a  rational parametrized analogue of \cite[Thm.~1.1]{ship} or \cite[Thm.~5.1.6]{schship}. More significantly, it may further be regarded as
 providing the stable counterpart of the unstable extension of Quillen's approach to rational homotopy theory given in \cite{bufemutan0}, which until now  has remained confined to the non-parametrized unstable context. 
 
In fact, the theorem above is the visible tip of a deeper structural theory.
Its proof unfolds through a sequence of explicit Quillen equivalences linking parametrized  spectra first to spectra of complete Lie algebras, then to spectra of $\compul$-modules and finally to the underlying module category itself. Each step isolates a different homotopical feature, making the passage from parametrized spectra to modules  fully transparent.
These intermediate equivalences are of independent interest, and the core of our work lies precisely in the construction and analysis of this pathway, which we now describe in detail.

We begin by recalling  that the aforementioned extension of Quillen's rational homotopy theory beyond the simply connected setting relies on the existence of a Quillen adjunction, given by the model and realization functors, 
$$
\xymatrix{ \catss& \catcdgl, \ar@<0.75ex>[l]^(.50){\langle\,\cdot\,\rangle}
\ar@<0.75ex>[l];[]^(.50){\lasu}\\}
$$
between the categories of simplicial sets and complete differential graded Lie algebras, the latter being equipped with a model structure that reflects the geometric behavior encoded by this correspondence.   This adjunction extends directly  to a Quillen pair
\begin{equation}\label{uno}
\xymatrix{ \retb & \cdgll \ar@<0.75ex>[l]^(.45){{\langle\,\cdot\,\rangle}_{\sslash B}}
\ar@<0.75ex>[l];[]^(.50){\lasu_{\sslash B}}\\}
\end{equation}
between the categories of retractive simplicial sets over $B$ (ex-spaces in the topological setting of \cite{maysi}) and  retractive cdgl's over $L=\lasu_B$.
We then gather the main homotopical aspects of the category $\cdgll$  which  provides the algebraic setting in which parametrization is encoded directly in the algebraic model of the given retractive cdgl.

Building on this, we adopt  Hovey's general stabilization template \cite{ho1}, initiated by Schwede \cite{schwe}, to develop the stable model category $\specl$ of spectra of retractive cdgl's ($L$-spectra).  This requires establishing that $\cdgll$ is proper and combinatorial (Proposition \ref{combi}), and addressing the fact that, 
although $\cdgll$ carries a form of simplicial enrichment \cite[\S12.4.4]{bufemutan0}, it is not a simplicial model category in the usual sense.
Consequently, the adjoint endofunctors required for stabilization, retractive loops and retractive suspension,  are not available a priori and must be constructed from scratch.
In our approach, the retractive loop functor is  explicitly obtained via a path object construction in the category of retractive cdgl's, involving completed tensor products with the Sullivan interval, see Definition \ref{loopsl}, Proposition \ref{lupquasi} and Corollary \ref{lupquasicoro}. The retractive suspension is also given explicitly up to weak equivalence in Theorem \ref{suspensionex}.

 Then, see Theorems \ref{primero}, \ref{primeroq} and \ref{primeroqq}, we prove:

\begin{introthm}\label{B} The retractive model and realization functors in (\ref{uno})  prolong to a Quillen pair
$$
\xymatrix{ \specb& \specl \ar@<0.75ex>[l]^(.50){\bm{\langle\,\cdot\,\rangle}}
\ar@<0.75ex>[l];[]^(.50){\lasub}\\}
$$
between the stable model categories of spectra of  retractive simplicial sets over  $B$ and spectra of retractive cdgl's over $L$. 
Furthermore, this adjunction induces a Quillen equivalence
$$
\xymatrix{ \spec_B^\bq& \specl \ar@<0.75ex>[l]^(.50){\bm{\langle\,\cdot\,\rangle}}
\ar@<0.75ex>[l];[]^(.50){\lasub}\\}
$$
after performing a Bousfield localization of the  stable model category of retractive simplicial sets at rational equivalences of $B$-spectra.
\end{introthm}
We emphasize that no connectivity or finiteness assumptions on the base are required for this equivalence. In fact, it allows for a natural encoding of the (possibly non-connected) base space and of the parametrization of a given $B$-spectrum in terms of Maurer-Cartan data in the associated $L$-spectrum. Furthermore,  a model of such a $B$-spectrum can be assembled from Lie models at each level (see Theorem \ref{freespec}) and, in line with the general simplifications that emerge upon passing to the stable setting, this model is stably equivalent to both its linear and its indecomposable reductions (Proposition \ref{estahomolin}).

Also, the Quillen equivalence in the above result provides a homological characterization of stable equivalences in terms of derived fiberwise homology (see Theorem \ref{homotohomolo} and Corollary \ref{sorpresa}). More precisely, the fiberwise homotopy groups of a $B$-spectrum correspond, under this equivalence, to the stable homology groups of the associated $L$-spectrum. These are, see Definition \ref{sorpresa2}, the derived  (stable)  homology  obtained from perturbations  at Maurer-Cartan elements.

\medskip

As one might expect, the essential homotopy-theoretic content of a retractive cdgl over $L$ can be extracted from its $L$-module part. This comes from decomposing a retractive cdgl as a twisted product of the base $L$ and its fiber (the kernel of the projection)  which naturally inherits a complete $L$-module equipped with the adjoint action. In other words, the category $\cdgll$ is closely related to the category $\cdgm_L$ of complete $L$-modules which is in turn equivalent, see Theorem \ref{equivalencia}, to the category $\cdgmul$ of complete differential graded modules over the completion of the universal enveloping algebra of $L$. 

This relationship is made precise in Definition \ref{calker}, by an adjunction
$$
\xymatrix{ \cdgmul& \cdgll \ar@<0.75ex>[l]^(.43){\calker}
\ar@<0.75ex>[l];[]^(.55){\libc_L}\\}
$$
in which the ``fiber'' functor $\calker$ sends each retractive cdgl to the fiber of its projection over $L$ while the left adjoint $\libc_L$ assigns to a complete $\compul$-module $R$,  the retractive cdgl whose fiber over $L$ is the (complete) free Lie algebra generated by $R$. 

On $\cdgmul$, suspension and desuspension define mutually inverse endofunctors, forming a suitable framework to develop its stabilization in the sense of Hovey. However, to connect the stable model category of $\compul$-spectra with that of $L$-spectra we need to restrict at this stage to the connected setting: $L$ itself is assumed connected and we consider only connected $\compul$-modules. This restriction, harmless in the stable category, ensures the detection of weak equivalences by both functors involved, see Theorem \ref{mainmod}, and the natural characterization of stable equivalences of $\compul$-modules, given in Proposition \ref{homomodes}, in terms of their stable homology. All of that is needed to prove (see Theorem \ref{segundo}):

\begin{introthm}\label{C} The  Quillen pair  
$$
\xymatrix{ \cdgmul^0& \cdgll \ar@<0.75ex>[l]^(.43){\calker^0}
\ar@<0.75ex>[l];[]^(.55){\libc_L}\\}
$$
 prolongs to a Quillen equivalence
$$
\xymatrix{ \specmodu & \specl. \ar@<0.75ex>[l]^(.43){\calker^0}
\ar@<0.75ex>[l];[]^(.55){\libcb_L}\\}
$$
\end{introthm}

We complete the algebraization process  by proving that spectra of connected $\compul$-modules
are Quillen equivalent to $\compul$-modules themselves.  This is a simple adaptation of the classical recognition principle for $\Omega$-spectra of modules in the ungraded, non-complete setting, see for instance \cite[Prop.~4.7]{ship}: in the connected context, the stable homotopy type of a fibrant module spectrum is entirely determined by the $0$th level and its module structure so that stabilization adds no new information. In Proposition \ref{casifin} we verify:

\begin{introprop} \label{D}
The  adjunction
$$
\xymatrix{ \specmodu & \cdgmul, \ar@<0.75ex>[l]^(.54){\mathcal{C}}
\ar@<0.75ex>[l];[]^(.44){\mathcal{D}}\\}
$$
defined for $R\in\specmodu$ and $T\in\cdgmul$ by 
$$
\mathcal{D}R=\varinjlim_n s^{-n+1}R^n,\quad 
(\mathcal{C}T)^n=(s^{n-1}T)^{(0)},
$$
 is a Quillen equivalence.
\end{introprop}

Combining  Proposition \ref{D} with Theorems \ref{B} and \ref{C}, we conclude (see Theorem \ref{mainmain}):

\begin{introthm}\label{E} Let $B$ be a connected simplicial set and let $L=\lasu^*_B$. Then, the functor
  $$
 \Psi\colon  \Ho\spec_B \longrightarrow\Ho\cdgmul
  $$
induced by the composition of the  sequence of Quillen adjunctions, 
$$
\xymatrix{ 
\spec_B 
  \ar@<0.75ex>[r]^(.50){\lasub^*} 
& 
\specl 
  \ar@<0.75ex>[l]^(.45){\bm{\langle\,\cdot\,\rangle}}
    \ar@<-0.80ex>[r]_(.45){\calker^0} 
& 
\specmodu 
  \ar@<-0.60ex>[l]_(.50){\libc_L} 
  \ar@<0.50ex>[r]^(.45){\mathcal{D}} 
& 
\cdgmul, 
  \ar@<0.90ex>[l]^(.55){\mathcal{C}} 
}
$$
induces an equivalence of categories
$$
 \Psi_\bq\colon  \Ho\spec_B^\bq\stackrel{\cong}{\longrightarrow}\Ho\cdgmul.
  $$
\end{introthm} 

The naturality of $\Psi$ is established, mainly through  direct arguments, in Proposition \ref{cambiobase} (see also Corollary \ref{coronatural}). Given a map of simplicial sets
 $f\colon B\to B'$ let $\varphi\colon L\to  L'$ denote the induced morphism $\lasu_f\colon\lasu_B\to\lasu_{B'}$. This gives rise to a change of base functor  (see Proposition \ref{changespec}) and its algebraic counterpart, the derived extension of scalars: 
$$
 \specb\longrightarrow\spec_{B'} \quad\text{and}\quad \Ho\cdgmul\to \Ho\cdgmulpri
$$
induced by $f$ and $\varphi$ respectively.
 Then:

  \begin{introprop}\label{G}
  The following diagram commutes:
  $$
\xymatrix{
\Ho\specb\ar[d]\ar[r]^(.40){\Psi}&\Ho\cdgmul\ar[d]\\
\Ho\spec_{B'}\ar[r]_(.40){\Psi}&\Ho\cdgmulpri.
}
$$
\end{introprop}

Finally, the strong monoidal character of our construction is proven in Theorem \ref{monoidal} which relies on a detailed analysis of Lie models of retractive smash products:

\begin{introthm}\label{F} The functor $\Psi$ is strong monoidal.
\end{introthm}

Combining all the results above we obtain an explicit formulation of Theorem A that makes the underlying strategy fully transparent (see Theorem \ref{mainmain}).
Furthermore, from a computational perspective and through explicit models, our approach gives direct access to stable invariants, bypassing the intermediate stages of spectra of retractive cdgl's and of $\compul$-modules, see Theorems \ref{teorema1} and \ref{teocomputable}. This leads to concrete and immediately applicable correspondences.
For example, see Proposition \ref{homotospec} the fiberwise stable homotopy groups of a $B$-spectrum correspond exactly to the suspension  of the  homology  of the underlying $\compul$-module.
In particular, the suspension isomorphisms in stable homotopy translate directly into the usual degree-shift isomorphisms in homology, so that stable suspension on the simplicial side becomes literal suspension (degree shift) in the algebraic category. In the same direction, see Corollary  \ref{fffinal}, the stable homotopy groups of the internal smash product correspond to the homology of the complete tensor product while, dually, the bifunctor of fiberwise stable homotopy classes of maps corresponds to  the differential $\ext_{\compul}$, see Proposition \ref{extfun}.

\medskip

A brief historical note is now in order.
The first systematic connection between rational fiberwise stable constructions and module categories over the rational model of the base appears in \cite{femutan}.
There, the authors develop a model for fiberwise rational stable homotopy theory in terms of module spectra over the commutative differential graded algebra modeling the base space, providing algebraic descriptions of fiberwise suspension spectra and their homotopy invariants.

The broader perspective underlying the present work was later raised in \cite{sch} where U. Schreiber  suggested the existence of a general equivalence between the stable category of rational parametrized spectra and the homotopy category of modules over algebraic models of the base. In that discussion, however, and based on a counterexample arising outside the rational context he was advised to restrict attention to the simply connected case.  This line of thought set the stage for Braunack-Mayer's subsequent contribution.

Indeed, as previously mentioned, Theorems 1.1, 1.4 and 1.5 of  \cite{brau2} established a natural strong monoidal equivalence between the homotopy category of rational parametrized spectra over a simply connected base $B$ and that of $UL$-modules, where $L$ is the   Quillen model of $B$.
This is a deep and conceptually powerful contribution representing a significant advance. 
Nevertheless, the simply-connected hypothesis inevitably limits its range of applicability, and  the  algebraic machinery involved is not primarily designed for explicit computations.

\medskip

Although the table of contents reflects the structure of the paper, we briefly outline its organization here.

Section 1 contains the first stage of our program. We begin with a short review of the homotopy theory of retractive simplicial sets and their spectra, together with the essential features of the homotopy theory of complete differential graded Lie algebras.
On this basis, we then construct the model category $\cdgll$, along with the endofunctors of retractive loops and suspension, leading to the stable model category of $L$-spectra, for which we introduce explicit models and the stable homology. Finally, we construct the bridge between $\specb$ and $\specl$, and prove Theorem \ref{B}.

In Section 2, we briefly recall  the foundations of complete modules over a general complete differential graded Lie algebra, and show that for a cdgl $L$, complete $L$-modules coincide with complete $\compul$-modules.
We then construct spectra on this category, relate them to $\specl$, and prove Theorem \ref{C}.

Section 3 is brief and contains the proof of Proposition \ref{D}.

Section 4 first gathers the preceding results to obtain Theorem \ref{E}, and then illustrates the computational power of our framework through explicit examples of the correspondence between homotopy invariants of parametrized $B$-spectra and those of $\compul$-modules. This section also addresses  functoriality under change of base and the proof of Proposition \ref{G}. 

Section 5 is devoted to the proof of Theorem \ref{F} which requires a careful construction of Lie models of  both, the  internal and external parametrized smash products.

Finally, the Appendix  of Section 6  recalls Hovey's machinery for constructing spectra in general model categories, together with some general results on retractive model categories and their spectra.
We also include a discussion of transferred model structures through an adjunction, extended to Bousfield localizations, retractive categories, and spectral categories.

\medskip

As  general conventions throughout the paper, we do not distinguish between a category and the class of its objects.
The symbol $\sim$ denotes a weak equivalence in the relevant category, while $\simeq$ denotes a quasi-isomorphism, which, depending on the context, may not coincide with a weak equivalence. Finally, whenever an adjoint pair of functors is displayed, the upper arrow will always denote the left adjoint.

\section{From spectra of retractive spaces to spectra of retractive Lie algebras}

\subsection{Spectra of retractive spaces}\label{prelimi}

None of the material in this section is new. We merely collect the concepts and results from retractive homotopy theory that are required to formally develop their spectra. Classical references for this are \cite{crabbja} and \cite{maysi}. However, we find the treatment in   \cite{brau1} and \cite{mal} to be particularly well-suited to our purposes. As usual, we denote by $\catss$ and $\catss^*$ the categories of  simplicial and pointed simplicial sets, respectively.

Although we are aware of the important and subtle differences between working with simplicial sets and spaces in the parametrized setting, we will, unless strictly necessary, treat simplicial sets and their geometric realizations interchangeably. Indeed, geometric realization does not interact as transparently with fiberwise constructions as one might expect: pullbacks and fiberwise smash products generally require homotopical corrections and suitable fibrancy assumptions in order to preserve the expected homotopy type. Nevertheless, since all constructions and results in this paper are formulated in the simplicial framework, and since fibrant replacements are introduced whenever required to ensure the correct homotopy type of the constructions involved, these phenomena do not play any essential role here. We therefore use topological terminology only for intuition and convenience, and freely identify a simplicial set with its geometric realization whenever no ambiguity can arise. In particular, in keeping with this convention, points of topological spaces will also be identified with the corresponding $0$-simplices of simplicial sets.

Let $B$ be a simplicial set. The category $\retb$ of {\em retractive simplicial sets over $B$}, referred to as {\em ex-spaces} in \cite{maysi}, is the retractive category of  $\catss$, see \S\ref{apret}.  In other words, an object $X$ of $\retb$ is a map (retraction or projection) of simplicial sets $ X\to B$ equipped with a section. A morphism of $\retb$, simply denoted by $X\to Y$, is given by a commutative diagram
$$
\xymatrix@R=10pt@C=10pt{
&B\ar[dl]\ar[dr]&\\
X\ar[dr]\ar[rr]&&Y\ar[dl]\\
&B&
}
$$
where the diagonal arrows define $X$ and $Y$ as retractive spaces.

By general arguments, see \S\ref{apret}, the usual combinatorial and proper model structure\footnote{The cofibrant generation of $\catss$ is well known. For its properness and local presentability see for instance \cite[Thm.~13.1.13]{hirsch0} and \cite[5.2.2b]{bor} respectively.} on $\catss$ induces a combinatorial model structure in $\retb$. Moreover, $\retb$ is both a simplicial and $\catss^*$-model category\footnote{As usual, the monoidal structures of $\catss$ and $\catss^*$ are given by the product and the smash product respectively.}. Although the enrichment over certain monoidal categories (such as $\catss$ or $\catss^*$) of a given category naturally extends the corresponding retractive categories, we prefer to briefly make this explicit for retractive spaces, see \cite[\S1]{brau1} or \cite[\S3]{mal} for further details.

By general categorical arguments, see \S\ref{apret}, given a map $f\colon A\to B$ of simplicial sets, the {\em change of base} adjunction  (\ref{ap4}),
$$
\xymatrix{ \catss_{\sslash A} & \retb\ar@<0.75ex>[l]^(.50){f^*}
\ar@<0.75ex>[l];[]^(.50){f_!}\\}
$$
 is always a Quillen pair, and becomes a Quillen equivalence whenever $f$ is a weak equivalence.

Given $X\in\ret{A}$ and $Y\in\ret{B}$,  their {\em external smash product} is defined as the space $X\,\bar\wedge\, Y\in\ret{A\times B}$ given by the following pushout:
\begin{equation}\label{pushout1}
\xymatrix{
(X\times B)\cup_{A\times B}(A\times Y)\ar[d]\ar[r]&X\times Y\ar[d]\\
A\times B\ar[r]&X\,\bar\wedge\, Y.
}
\end{equation}

The  {\em smash product} $X\wedge_B Y\in\retb$ of two objects $X,Y\in\retb$ is then defined as the pullback
\begin{equation}\label{pullback1}
\xymatrix{
X\wedge_B Y\ar[d]\ar[r]&X\,\bar\wedge\, Y\ar[d]\\
B\ar[r]^(.42){\Delta}&B\times B.}
\end{equation}

In particular, viewing $\catss^*$ as $\ret{*}$, for any pointed  $K\in\catss^*$ and any $X\in \retb$ we denote
$$
K\owedge X=K\,\bar\wedge \,X\in\retb
$$
which defines the tensor functor of the $\catss^*$-structure in $\retb$.
For the unpointed tensoring, given $K\in\catss$ and $X\in \retb$ define
$$
K\otimes X\coloneq K_{+}\owedge X\in\retb
$$
which is easily seen to coincide with the pushout
$$
\xymatrix{
K\times B\ar[d]\ar[r]&B\ar[d]\\
K\times X\ar[r]&K\otimes X.
}
$$
Here $K_+$ stands for $K\amalg *$, which is $K_{+*}$ with the notation in \S\ref{apret}.

The external smash product has a right adjoint in each variable: given $Z\in \ret{A\times B}$ and $Y \in \retb$ with sections and retractions $i$, $j$, and $p$, $q$ respectively, define the {\em external mapping space}
$\underline\map(Y,Z)\in\ret{A}$ as the pullback
 $$
\xymatrix{
\underline\map(Y,Z)\ar[d]\ar[r]&\map(Y,Z) \ar[d],\\\
A\times\{*\}\times \{*\}\ar[r]&\map(Y,A)\times \map(Y,B)\times\map(B,Z)
}
$$
where $\map$ stands for the usual simplicial unpointed mapping space, the right vertical map is $\bigl(\map(Y,\pi_A\circ p),\map(Y,\pi_B\circ p)\times \map(j,Z)\bigr)$ and the bottom horizontal map sends $a\in A$ to the triple consisting in the simplicial maps induced by the constant map on $a$, $q$ and $j$. One then checks that
$$
\Hom_{\,\ret{A\times B}}(X\,\bar\wedge\,Y,Z)\cong\Hom_{\,\ret{A}}\bigl(X,\underline\map(Y,Z)\bigr).
$$
By choosing $K\in \catss^*=\ret{*}$ and $X\in\retb$ we get
$$
X^K_*\coloneq\underline\map(K,X)\in \retb
$$
which exhibits the cotensor functor of the $\catss^*$-structure in $\retb$. Again, for the free cotensoring, given $K\in\catss$ and $X\in \retb$ define
$$
X^{K}\coloneq X^{K_{+}}_*\in\retb.
$$

As for any simplicial model category, the sequences 
$$\partial\Delta^1\otimes X\to \Delta^1\otimes X \to X\quad\text{and}\quad X\to X^{\Delta^1}\to X^{\partial\Delta^1}
$$ 
define, respectively, a natural  good cylinder and a good loop object of  $X$. In particular, the  {\em retractive suspension} and {\em retractive loops} of $X$ are functorially defined, respectively,  by the following pushout and pullback
$$
\xymatrix{ \partial\Delta^1\otimes X\ar[d]\ar[r]&\Delta^1\otimes X \ar[d] & \quad & \Omega_B X\ar[d]\ar[r]& X^{\Delta^1}\ar[d]\\ B\ar[r]&\Sigma_BX & \quad & B\ar[r]&X^{\partial\Delta^1}. }
$$

Again, as in any simplicial category, the retractive suspension and loops constitute a Quillen pair in $\retb$. It is convenient to note that 
choosing the pointed circle as $S^1=\Delta^1_+/\partial\Delta^1_+$ we get
$$
\Sigma_BX=S^1\owedge X.
$$
On the other hand, see \cite[Cor.~1.8]{brau1}, the change of base (as in Definition \ref{change}) commutes with retractive suspension: for any map $f\colon A\to B$ there is a natural isomorphism,
\begin{equation}\label{sus!}
f_!\,\sus_A\cong \susb\, f_!.
\end{equation}

We summarize all of the above as follows:

\begin{theorem} The retractive category $\retb$ is a proper combinatorial, simplicial and $\catss^*$-model category endowed with  a Quillen pair of endofunctors provided by the retractive suspension $\Sigma_B$ and retractive loops $\Omega_B$.\hfill$\square$
\end{theorem}

As for the stabilization of $\retb$, since it is proper and combinatorial, we simply adapt to this category the general setting of \S\ref{specmod}, taking  $\susb\dashv\lupb$ as the Quillen pair of endofunctors.

\begin{definition} A {\em spectrum} in $\retb$ or {\em $B$-spectrum} is an object $X$ of the category $\spec_B(\catss)$, see Definition \ref{specret}, which we denote $\spec_B$ henceforth. That is, $X$ is a sequence $\{X_n\}_{n\ge0}$ of retractive spaces over $B$ equipped with a family of {\em structure maps} in $\retb$
$$
\sigma\colon \susb X_n\longrightarrow X_{n+1},\quad\text{or equivalently,}\quad \eta\colon X_n\longrightarrow \lupb X_{n+1},\quad n\ge 0.
$$
A map  $f\colon X\to Y$ of $B$-spectra is a collection $\{f_n\colon X_n\to Y_n\}_{n\ge0}$ of maps in $\retb$ compatible with the structure maps in the obvious sense. From now on we simply write $\spec$ to denote $\spec_*(\catss)$, the category of spectra of pointed simplicial sets.
\end{definition}

Particularizing Sections \S\ref{specmod} and \S\ref{apret} to $\specb$, we list the features we require from this category:

\begin{proposition}\label{resub}
{\em (i)} The functors $\susb$ and $\lupb$ prolong to adjoint endofunctors
$$
\xymatrix{ \spec_B & \spec_B.\ar@<0.75ex>[l]^(.50){\lupeb}
\ar@<0.75ex>[l];[]^(.50){\suseb}\\}
$$

{\em (ii)} For each $k\ge 0$ there are adjoint functors defined as in (\ref{eva}) and (\ref{susk}),
$$
\xymatrix{ \retb & \spec_B\ar@<0.75ex>[l]^(.40){\eva_k}
\ar@<0.75ex>[l];[]^(.50){\susb^{\infty-k}}\\}
$$
\hfill$\square$
\end{proposition}

 Applied to this setting, Definition \ref{proestable}  reads:
 
 \begin{definition}\label{proyes} The {\em projective model structure}  on $\specb$ is defined so that the  fibrations and weak equivalences are the levelwise fibrations and levelwise weak equivalences, respectively. A map $f\colon X\to Y$ is a cofibration (or trivial cofibration)  if $f_0$ and the induced maps $ X_{n}\amalg_{\sus X_{n-1}}\susb Y_{n-1}\to Y_n$, $n\ge 1$, are cofibrations (or trivial cofibrations).

\smallskip 
 
Consider the usual sets
\begin{equation}\label{usugen}
  \mathcal I=\{\partial\Delta^n\hookrightarrow \Delta^n\}_{n\ge0}\esp{and}
  \mathcal J=\{\Lambda^n_i \hookrightarrow \Delta^n\}_{n\ge0, \,i=0,\dots,n},
\end{equation}
of generating cofibrations and trivial cofibrations of $\catss$ and write 
$$
\mathcal I_B=\{B\amalg i,\,i\in\mathcal I\}\esp{and}\mathcal J_B=\{B\amalg j,\,j\in\mathcal J\},
$$
which are sets of generating cofibrations and trivial cofibrations of $\retb$.
Then, see (\ref{generaes}), 
$$
\mathcal I_{\susb}=\cup_{k\ge 0}\susb^{\infty-k}(\mathcal I_B)\esp{and}\mathcal J_{\susb}=\cup_{k\ge 0}\susb^{\infty-k}(\mathcal J_B)
$$
are generating cofibrations and trivial cofibrations of $\specb$.
 \end{definition}

On the other hand, Definition  \ref{eses} particularizes to the following:

\begin{definition}\label{es}
Consider the family $S$ from (\ref{equibous}) (no cofibrant replacement is needed here since all objects in  $\retb$ are cofibrant), and define the {\em stable model structure} on $\specb$ as the Bousfield localization of $\specb$ with respect to $S$. 
\end{definition}

\begin{rem}\label{inutil} By Corollary \ref{changeret}, the pair in Proposition \ref{resub}(i) becomes a Quillen equivalence with respect to the stable structure in which  the fibrant objects are the $\lupb$-spectra. That is, those  $B$-spectra $X$ for which each $X_n$ is fibrant and the structure maps $X_n\to \lupb X_{n+1}$ are weak equivalences of $\retb$ for all $n\ge 0$.  As a consequence, a map  between $\lupb$-spectra is a stable equivalence if and only if it is a projective equivalence.
\end{rem}
 \begin{rem}\label{monoidal}
Note that $\specb$ is a basic example of a category of diagram spectra \cite{manmayschshi}. As such, and since we are not considering symmetric spectra \cite{hoshis}, only the homotopy category of $\specb$ carries a symmetric monoidal structure with respect to the smash product of spectra. For our purposes we use  the retractive handcrafted version  given by the Day convolution: for $X,Y\in\specb$ define $X\wedge_B Y$ levelwise as
$$
(X \wedge_B Y)_n = \mathop{\bigvee\!\raisebox{-0.8ex}{\scriptsize $B$}}\limits_{p+q=n} (X_p \wedge_B Y_q).
$$
The $n$th structure map is defined on the suspension of each term as
$$
\sigma_p\wedge_B\id +\id\wedge_B \sigma_q\colon \Sigma_B(X_p\wedge_B Y_q)\to (X_{p+1}\wedge_B Y_q)\vee_B (X_p\wedge_B Y_{q+1}).
$$
As in the non-retractive setting, this is the map
 $\bigl((\sigma_p\wedge_B\id) \vee_B(\id\wedge_B \sigma_q)\bigr)\circ\nu$  where $\sigma_p$ and $\sigma_q$ are the $p$th and $q$th structure maps of $X$ and $Y$ respectively, and
 $$
 \nu\colon \Sigma_B(X_p\wedge_B Y_q)\to (\Sigma_BX_p\wedge_B Y_q)\vee_B (X_p\wedge_B\Sigma_BY_q)
 $$
 is the retractive pinching followed by the appropriate weak equivalences.
 
The unit object for this monoidal  structure is  the {\em retractive sphere spectrum}  
$$
\bs_B=B\times\bs,
$$
i.e.,  the trivial bundle of spectra over $B$ with  the sphere spectrum $\bs$ as fiber. 
  \end{rem}
  
 \begin{definition}\label{lupinfinity}
 The {\em infinity loop space} $\lupb^\infty X\in\retb$ of a given spectrum $X\in\spec_B$ is the $0$-level space of a fibrant replacement of $X$. Its weak equivalence class in $\retb$ is independent of the choice of  fibrant replacement  since $\lupb^\infty X$ represents the right derived functor of $\eva_0$ applied to $X$.   
 \end{definition}
 In the simplicial context,  Corollary \ref{changeret} together with (\ref{sus!}) yields:
 
 \begin{proposition}\label{changespec}
For any map $f\colon A\to B$, the change of base
$$
\xymatrix{ \catss_{\sslash A} & \retb, \ar@<0.75ex>[l]^(.50){f^*}
\ar@<0.75ex>[l];[]^(.50){f_!}\\}
$$
induces a Quillen pair between the corresponding stable model categories
$$
\xymatrix{ \spec_A & \spec_B \ar@<0.75ex>[l]^(.50){\widetilde{f^*}}
\ar@<0.75ex>[l];[]^(.50){\widetilde{f_!}}\\}
$$
in which $\widetilde{f_!}$ and $\widetilde{f^*}$ denote the prolongations of $f_!$ and $f^*$, respectively. Moreover, this is a Quillen equivalence whenever $f$ is a weak equivalence.\hfill$\square$
\end{proposition}

\begin{rem}\label{conexo2} From a conceptual point of view there is no loss of generality restricting our attention to the case where $B$ is reduced. Indeed on the one hand, there is an equivalence of categories
$$
\spec_B\cong {\textstyle \prod_{B_i\in\pi_0(B)}}\, \spec_{B_i},
$$
and on the other hand, every connected simplicial set is weakly equivalent to a reduced one.
\end{rem}
 
 \begin{definition}\label{stableho}
(i) Given a $B$-spectrum $X$ and  $b\in B$, we denote by $X^b\in\spec$ the {\em fiber spectrum of $X$ at $b$} where
$X^b_n$ is the fiber of the retraction $p_n\colon Y_n\to B$ at $b$, along with the induced structure maps.  Observe that if $X$ is fibrant, then $X^b$ is also fibrant for any $b\in B$. Moreover, for any $X\in\spec_B$ and any $b\in B$, $(\lupb^\infty X)^b$ and $\Omega^\infty(X^b)$ are weakly equivalent. Note also that  each $X^b_n$ is pointed by $s_n(b)$, where $s_n\colon B\to Y_n$ denotes the section of $Y_n$.

\smallskip

(ii) The {\em (fiberwise) stable homotopy groups} $\pie_*(X)$ of a spectrum
$X\in\specb$ are defined as the collection
$$
\pie_*(X)=\{\pie_*\bigl((RX)^b\bigr)\}_{b\in B},
$$
where $R$ denotes a fibrant replacement functor in the projective model structure and
$(RX)^b\in\spec$ is the fiber spectrum at $b$. Equivalently, for each $n\ge0$, the $n$th level map $(RX)_n\to B$ is a fibration. This definition is independent of the chosen replacement and one may simply select  an arbitrary point in each path component of $B$. Sometimes it is useful to regard the fiber spectrum  at a point $b\in B$, as in \cite[Def.~2.14]{brau1}, as the functor $b^*\colon \specb\to \spec$ where $b\colon *\to B$. 

\smallskip

(iii) Given $X,Y\in\spec_B$, the graded abelian group $\{X,Y\}^B$ of {\em (fiberwise) stable homotopy classes of  maps} is defined by
$$
\{X,Y\}^B_k=\Hom_{\,\Ho\spec_B}(X,\susel^k Y),\quad k\in\bz.
$$
One can verify that for any $B$-spectrum $X$ there is a natural isomorphism 
\begin{equation}\label{stamap}
\pie(X)\cong \{\bs_B,X\}^B.
\end{equation}
\end{definition}

Then,
see \cite[Thm.~12.3.14]{maysi}, or \cite[Lemma 2.16]{brau1} in the retractive simplicial setting:

\begin{theorem}\label{eseq}
A morphism of $\specb$ is a stable equivalence if and only if it is a $\pie_*$-isomorphism, that is,   it induces isomorphisms in fiberwise stable homotopy groups.

\hfill$\square$
\end{theorem}

\begin{rem}
One could start the stabilization process using the functor $S^1\owedge -$, with the induced structure maps built upon the permutation isomorphism of $S^1\wedge S^1$. In the non-parametrized case, Corollary 3.5 of \cite{ho1} shows that the stable model structure in $\spec$ coincides with the classical one. In the parametrized version, the remarkable Lemma 2.17 of \cite{brau1} shows that the functor induced by $S^1\owedge -$ in the stable structure $\specb$ is a Quillen equivalence.
\end{rem}

We finish with an essential notion.

\begin{definition}\label{homoq} A map $f$ of $\specb$ is a rational equivalence if $\pie_*(f)\otimes\bq$ is an isomorphism. The {\em rationalization} of $\specb$ is the Quillen functor
$$
(\,\cdot\,)^\bq\colon\specb\longrightarrow \spec_B^\bq
$$
obtained by Bousfield localization with respect to rational equivalences.

Therefore, a $B$-spectrum $X\in\spec_B$ is $\bq$-local (or rational) if and only if the fiberwise stable homotopy groups $\pi^{\rm st}(X^b)$ are $\bq$-vector spaces. Equivalently, for any $b\in B$, the spectrum $X^b\in\spec$ is rational, that is,  stably equivalent to $X^b\wedge H\bq$.
\end{definition}

We will not need the following observation which is included only for completeness.

\begin{rem}\label{localracional}
As in the classical non-parametrized setting, see \cite{bous} for  the original reference, the rationalization of $\specb$ may be seen to be Quillen equivalent to a certain parametrized symmetric smashing \cite[\S2.3]{brau2}. As a consequence, rational parametrized spectra are, up to stable equivalence, symmetric parametrized spectra over $B$ whose fibers are rational Eilenberg-MacLane spectra with structure maps compatible with the parametrization \cite[Prop.~ 2.16]{brau2}.
\end{rem}

\subsection{Retractive cdgl's}

\subsubsection{A few insights into homotopy theory of complete differential graded Lie algebras}\label{homotocdgl} 
Most of the material presented here is covered thoroughly  in the comprehensive reference \cite{bufemutan0} or in  \cite{bufemutan1, bufemutan2}. What follows is a brief summary of the key concepts.
Unless otherwise specified, all algebraic objects are assumed to be $\bz$-graded, and defined over  $\bq$.

A {\em complete differential graded Lie algebra} (cdgl) is a differential graded Lie algebra (dgl)  $L$, or $(L,d)$  to emphasize the differential, equipped with a decreasing  dgl filtration  $\{F^n\}_{n\ge 1}$ such that $F^1=L$, $[F^n,F^m]\subset F^{n+m}$, for all $n,m\ge1$, and the  natural map
$$
L\stackrel{\cong}{\longrightarrow}\varprojlim_n L/F^n
$$
is a dgl isomorphism. Morphisms of cdgl's are required to preserve the filtration and we denote by $\catcdgl$ the corresponding category which is bicomplete \cite[\S3.1]{bufemutan0}. 

The {\em completion} of a filtered dgl $L$ is
$
\widehat L=\varprojlim_nL/F^{n}
$
which is always complete with respect to the filtration
$
\widehat{F}^n=\ker ( \widehat L \to L/F^n)
$ as $\widehat L/\widehat F^n\cong L/F^n$.

 Let $F=\{F^n\}_{n\ge 1}$ and $G=\{G^n\}_{n\ge 1}$ be filtrations of a dgl $L$ with $F^n\subset G^n $ for all $n$. Observe that if $L$ is complete with respect to $G$, then it is also complete with respect to $F$. In particular any cdgl $L$ is automatically complete with respect to the {\em adic} filtration of $L$, that is, the lower central series filtration $\{L^n\}_{n\ge 1}$ where $L^1=L$ and $L^n=[L^{n-1},L]$. 
 
Given $\lib(V)$, the free Lie algebra generated by the graded vector space $V$, we denote by
$$
\libc(V)=\varprojlim_n\lib(V)/\lib(V)^n$$
its completion with respect to the  lower central series filtration and call it the {\em free complete Lie algebra} generated by $V$ \cite[\S3.2]{bufemutan0}. Abusing notation, we refer to a cdgl of the form $(\libc(V),d)$ as a {\em free cdgl}.  

 The {\em Maurer-Cartan set} $\mc(L)$ of a given dgl $L$ consists of elements $a\in L_{-1}$ that satisfy the Maurer-Cartan equation $da=-\frac{1}{2}[a,a]$. For any cdgl $L$ and any  $a\in \mc(L)$, we denote by $d_a= d+\ad_a$ the {\em perturbed differential} where $\ad$ is the usual adjoint operator. The {\em component} of $L$ at $a$ is the connected sub dgl $L^{(a)}$ of $(L,d_a)$ given by
\begin{equation}\label{trunca}
L^{(a)}_p=\begin{cases} \ker d_a&\text{if $p=0$},\\ \,\,\,L_p&\text{if $p>0$}.\end{cases}
\end{equation} 

The homotopy theory in $\catcdgl$ is structured around a pair of adjoint functors,   {\em (global) model} and {\em realization} \cite[\S7]{bufemutan0},
\begin{equation}\label{pair}
\xymatrix{ \catss& \catcdgl \ar@<0.75ex>[l]^(.50){\langle\,\cdot\,\rangle}
\ar@<0.75ex>[l];[]^(.50){\lasu}\\},
\end{equation}
developed from the cosimplicial cdgl $\lasu_\bullet=\{\lasu_n\}_{n\ge 0}$, and defined as follows:
\medskip

For any $n\ge 0$, denote by $s^{-1}\Delta^{n}$  the desuspension of the non-degenerate simplicial chains  on the simplicial set $\underline\Delta^n$. Then, 
$$
\lasu_n=\bigl(\libc(s^{-1}\Delta^{n}),d)
$$
in which  $d$ is the only differential (up to cdgl isomorphism) for which: the generators of $s^{-1}\Delta^{n}$, corresponding to vertices, are MC elements; 
 the linear part of $d$ is induced by the boundary operator of $s^{-1}\Delta^{n}$; and the cofaces and codegeneracies are simply induced in the corresponding free cdgl by the cosimplicial vector space  $s^{-1}\Delta^\bullet$ \cite[\S6]{bufemutan0}. 
 
 For instance $\lasu_1=\libc(s^{-1}\Delta^1),d)$ is the {\em Lawrence-Sullivan interval} \cite{lawsu}. This is the free cdgl $(\libc(a,b,c),d)$ where $a$ and $b$ are Maurer-Cartan elements and
 $$
 dc=\ad_cb+\sum_{n\ge0}\frac{B_n}{n!}ad_c^n(b-a),
 $$
 being $B_n$ the $n$th Bernoulli number. A {\em path} on a cdgl $L$ is a cdgl morphism  $\lasu_1 \to L$ which is therefore characterized by elements $x,y\in \mc(L)$ and $z\in L_0$ whose differential is as above. Then, the {\em gauge} action of $L_0$  (endowed with the BCH product) on $\mc(L)$ is given by $z{\mathcal G}y=x$ \cite[\S4]{bufemutan0}. 
     We denote by $\widetilde\mc(L)$ the orbit set $L_0/{\mathcal G}$.

On the other hand, the realization of a given cdgl $L$ is the simplicial set
     $$
     \langle L\rangle=\Hom_{\catcdgl}(\lasu_\bullet,L),
     $$
  whereas, for any simplicial set $X$,
     $$
     \lasu_X=\varinjlim_{\sigma\in X}\lasu_{|\sigma|}.
     $$
     In other words,   $
\lasu_X=(\libc(s^{-1}X),d)$
 where  $s^{-1}X$ denotes the desuspension of the chain complex of non-degenerate simplicial chains on $X$; the  $0$-simplices of $X$ are Maurer-Cartan elements; 
and the linear part of $d$ is the  boundary operator of $s^{-1}X$.

Among other properties of the model and realization functors \cite[\S7]{bufemutan0} we emphasize that given a cdgl $L$,
\begin{equation}\label{muchas}
\pi_0\langle L\rangle=\widetilde\mc(L) \quad\text{and}\quad \langle L\rangle\simeq \amalg_{a\in\widetilde\mc(L)}\langle L^a\rangle.
  \end{equation}
Furthermore, if $L$ is connected, then $\langle L\rangle$ is  reduced and there are isomorphisms
 \begin{equation}\label{homoto}
 \pi_n\langle L\rangle\cong  H_{n-1}(L),\quad\text{for any $n\ge 1$,}
  \end{equation}
where the group structure in $H_0(L)$ is considered with the Baker-Campbell-Hausdorff (BCH) product.

The category $\catcdgl$ inherits a {\em cofibrantly generated model structure} from $\catss$ via right transfer \cite[Chapter 8]{bufemutan0}, making the functors in (\ref{pair}) into a Quillen pair \cite[\S8]{bufemutan0}:

\begin{itemize}

\item[] Fibrations are cdgl morphisms which are surjective in non negative degrees.
 
\item[]   Weak equivalences are morphisms $f\colon L\to M$ such that the following holds:

(i) $\widetilde\mc(f)\colon \widetilde\mc(L)\stackrel{\simeq}{\to}\widetilde\mc(M)$ is a bijection.

(ii) For any $a\in \mc(L)$, the induced morphism, $f^{(a)}\colon L^{(a)}\stackrel{\simeq}{\to} M^{(f(a))}$ is a quasi-isomorphism.

\item[]  Lastly, cofibrations and acyclic cofibrations are generated, respectively, by $\lasu({\mathcal I})$ and $\lasu({\mathcal J})$ with $\mathcal I$ and $\mathcal J$ as in (\ref{usugen}).
 \end{itemize}

Thus, in general, quasi-isomorphisms are not  weak equivalences and fibrations need not  be surjective. These characterizations hold only in the connected setting.  Nevertheless the {\em Goldman-Millson Theorem} \cite[Thm.~4.33]{bufemutan0} provides certain cases where a quasi-isomorphism is indeed a weak equivalence. Specifically, let  $f\colon L\to L'$ be a morphism of cdgl's filtered by $\{F^n\}_{n\ge 0}$ and $\{G^n\}_{n\ge 0}$, respectively, such that the induced map $F^n/F^{n+1}\stackrel{\simeq}{\to}G^n/G^{n+1}$ is a quasi-isomorphism for all $n$. Then, $f$ is a weak equivalence.

Any  quasi-isomorphism of connected cdgl's of the form
$$
(\libc(V),d)\stackrel{\simeq}{\longrightarrow} L
$$
exhibits  $(\libc(V),d)$ as a cofibrant replacement of $L$ and we say that this is a {\em  Lie model of $L$} \cite[\S8.4]{bufemutan0}. If $d$ has no linear term we say that $(\libc(V),d)$ is the {\em  minimal model of} $L$, and is unique up to cdgl isomorphism.

The {\em minimal model} of a reduced simplicial set $X$ is, by definition, the  minimal model $(\libc(V),d)$ of  $\lasu_X^{(a)}$, where $a$ is the Maurer-Cartan element corresponding to the only $0$-simplex of $X$. This coincides with the minimal model of $\lasu_X/(a)$, where $(a)$ denotes the ideal generated by $a$, since the composition $\lasu_X^{(a)}\hookrightarrow \lasu_X{\to}\lasu_X/(a)$ is a quasi-isomorphism.  

As shown in \cite{fefuenmu0},  the Quillen pair (\ref{pair}) extends, up to homotopy, the classical functors of Quillen \cite{qui}. More generally \cite[Thm.~0.1]{fefuenmu1}, for any reduced simplicial set $X$ and for any $0$-simplex $a\in X$, the unity of (\ref{pair}),
$$
X\longrightarrow \langle \lasu_X^{(a)}\rangle
$$
 is weakly equivalent to the {\em Bousfield-Kan $\bq$-completion} of $X$
$$
X\longrightarrow X^\acento_\bq.
$$
In particular, by (\ref{muchas}), for any simplicial set $X$ with path components $\{X_i\}$, we have
\begin{equation}\label{bousfikan}
\langle \lasu_X\rangle\simeq \amalg_{i}\, {(X_i)}^\acento_\bq\amalg\{*\}
\end{equation}
where the extra point reflects that $\langle \lasu_X^{(0)}\rangle\simeq\{*\}$ for the ubiquitous Maurer-Cartan element $0$.

It is important to note that perturbing a cdgl $L$ by a Maurer-Cartan element does not change the homotopy type of its realization:
\begin{equation}\label{nodif}
\langle L\rangle\cong \langle L,d_a\rangle\quad\text{for  any $a\in \mc(L)$}.
\end{equation}
However, when $L=\lasu_X$ and $a$ is a $0$-simplex of $X$, the realization of the component of  $(L,d_a)$ at $0$ recovers the $\bq$-completion of the component $X_a$ of $X$ at $a$, while the realization of the component of $(L,d_a)$ at $a$ is contractible:
$$
\langle (\lasu_X,d_a)^{(0)} \rangle\simeq (X_a)^\acento_\bq\quad \text{and}\quad \langle (\lasu_X,d_a)^{(a)}\rangle\simeq \{*\}.
$$
We finish by remarking that the model functor  can be slightly altered while still yielding a Quillen pair
\begin{equation}\label{pair2}
\xymatrix{ \catss^*& \catcdgl \ar@<0.75ex>[l]^(.50){\langle\,\cdot\,\rangle}
\ar@<0.75ex>[l];[]^(.50){\lasu^*}\\}
\end{equation}
where now, for each pointed simplicial set $(X,a)$, the {\em pointed model functor } is defined as $\lasu^*_X=\lasu_X/(a)$.

\subsubsection{Retractive cdgl's}

Let $L$ be a cdgl.

\begin{definition}\label{recdgl} The {\em category  of retractive cdgl's over $L$} is the category $\cdgll$, see \S\ref{apret}. As such, an object of this category consists of a pair of cdgl morphisms (called the {\em section}  and the {\em retraction} or {\em projection})
$$
L\stackrel{s}{\longrightarrow}M\stackrel{p}{\longrightarrow} L
$$
whose composition is the identity. Morphisms are defined accordingly. As usual, $\cdgll$ inherits from $\catcdgl$ a model structure. 

Observe that any retractive cdgl $M$ admits a splitting   of the form
\begin{equation}\label{splitting}
M\cong L\oplus K_M,
\end{equation}
where $K_M$, or $K$ when no ambiguity arises, is the ideal $\ker p$, which is a complete dgl. Furthermore, if $a\in \mc(L)$, write $s(a)=a$ and note that the perturbed differential $d_a$ restricts to $K$ so that
\begin{equation}\label{splittingmc}
(M,d_a)\cong(L,d_a)\oplus (K,d_a).
\end{equation}
It also follows that any morphism $M\to N$ of retractive cdgl's induces a morphism of ideals $K_M\to K_N$. 
\end{definition}

\begin{definition}\label{libre} Abusing terminology, we say that a retractive cdgl is  {\em free} if it is of the form $(L\amalg\,\libc(W),d)$ where the coproduct is taken in the category of complete Lie algebras.
For such a retractive cdgl, the decomposition  (\ref{splitting}) takes the form
\begin{equation}\label{esencialk}
(L\amalg\libc(W),d)\cong L\oplus (\libc(T),d)
\end{equation}
where  $T$ is the subspace of $ L\amalg\libc( W)$ spanned by brackets containing exactly one element of $W$. Indeed, 
 $ \libc(T)\subset K$ by definition of $K$. Conversely, consider an element of $K$ consisting of a single bracket  containing exactly $n$ elements of $W$ with $n\ge 1$. By induction and the Jacobi identity, this element can be written  as a linear combination of brackets each of which lies in $ \lib^n(T)$.  This implies that $K\subset \libc(T)$ and equality holds. 
 Observe that  $[L,T]\subset T$.
 
 A free retractive cdgl $(L\amalg\libc(W),d)$ is said to be {\em (retractive)  linear} if, for each $w\in W$, $dw$ is a linear combination of brackets, each of which contains exactly one element from $W$.

Note that the differential of any free retractive cdgl $(L\amalg\libc(W),d)$ can be decomposed as
$$
d=\sum_{i\ge 1}d_i,\quad \text{with}\quad d_iW\subset\Span\{\text{brackets containing exactly $i$ elements from $W$}\}.
$$
 In particular  $d_1$ is a differential in $L\amalg\libc(W)$ and we call  it the {\em retractive linear part} of $d$. Note that
$$
(L\amalg\libc(W),d_1)\cong L\oplus (\libc(T),d_1).
$$
where, on the right-hand side $d_1$ denotes the usual, non-retractive, linear part of the differential in $\libc(T)$.
\end{definition}

Any object in $\cdgll$ admits a cofibrant replacement given by a specific free retractive cdgl:  

\begin{proposition} \label{cofire} For any $M\in
\cdgll$ there is a weak equivalence 
$$
(L \amalg  \libc(W),d)\stackrel{\sim}{\longrightarrow} M
$$
where:

 $W = W_{\geq -1} $ and $W_{-1}$ is generated by Maurer-Cartan elements;

$W_0= W_0' \oplus W_0''$, where  $dW_0'=0$ and $W_0''$ is generated by paths between MC elements in  $L \amalg \libc(W_{-1})$;

for $x\in W_n$, with $n\geq 1$, there is a Maurer-Cartan element $a$ such that $d_ax\in L \amalg {\widehat{\mathbb L}}(W_{<n})$.

\end{proposition}
\begin{proof} By \cite[Prop.~8.14]{bufemutan0}, we have a commutative diagram in  $\catcdgl$
$$
\xymatrix{L \ar[r]^s\ar[rd]&  M\ar[r]^p&L\\&(L\amalg\libc(W),d) \ar[u]_\sim\ar[ru]&}$$
where $(L\amalg\libc(W),d)$ satisfies the required conditions and the vertical arrow is a weak equivalence. 
\end{proof}

\begin{definition}\label{modeloretractivo} We refer to the weak equivalence $(L\amalg \libc(W),d)\stackrel{\sim}{\to} M$, or simply the object $(L\amalg \libc(W),d)$, as  a {\em retractive model of $M$}. By \cite[Thm.~8.12]{bufemutan0} the morphism  $L\to (L\amalg\,\libc(W),d)$ is a cofibration in $\catcdgl$ and consequently such a model is a cofibrant replacement of $M$  in $\cdgll$.
\end{definition}

In the connected setting, the previous result reads:

\begin{corollary}\label{cofirecon} Let $L$ be a connected cdgl. Then, any connected $M\in\cdgll$ admits a model of the form
$$
(L  \amalg  \libc(W),d)
$$
where $W=W_{\ge 0}$ and  $dx\in L \amalg {\widehat{\mathbb L}}(W_{<n})$ for each $x\in W_n$.
\hfill$\square$
\end{corollary}

We conclude with a key result for developing spectra on $\cdgll$.

\begin{theorem}\label{combi}
 $\cdgll$ is a proper combinatorial model category.\hfill$\square$
 \end{theorem}

\begin{proof} In view of Corollary \ref{combret} this result is equivalent to its non-retractive analogue and we are thus reduced to prove that $\catcdgl$ is both proper  and combinatorial. The combinatorial nature is straightforward to verify: as noted in  S\ref{homotocdgl}, $\catcdgl$ admits all colimits and is cofibrantly generated. Moreover, by  standard arguments establishing local presentability for many algebraic categories (see for instance \cite[\S5]{bor}), one  readily checks that the free cdgls form a set of compact objects generating  $\catcdgl$ under colimits. Also, the right properness of $\catcdgl$ follows from the fact that every object is fibrant, see \cite[Thm.~ 5.1.1]{hirsch0}.

To show left properness\footnote{The authors believe that the proof  of the general statements \cite[Thm.~4.17]{whi} or \cite[Thm.~0.1]{baber} could be adapted to show left properness in $\catcdgl$. For completeness and to be precise with details that could be crucial, we have chosen to present a direct proof.} we reduce the problem as follows. Since directed limits preserve weak equivalences in $\catcdgl$, and any cofibration is a (retract of) transfinite composition of pushouts along generating cofibrations, it is enough to show the following: the pushout of a morphism, which is itself the pushout of a map in $ \lasu({\mathcal I})$, along another morphism, remains a weak equivalence when taken along any weak equivalence.

Explicitly, a general morphism of $\lasu({\mathcal I})$ is of the form $\lasu_{\partial\Delta^n}\hookrightarrow \lasu_{\Delta^n}$. Pushing out this map along any morphism   $\lasu_{\partial\Delta^n}\to L$ results in the following diagram, see for instance \cite[Prop.~6.8]{bufemutan0}:
$$
\xymatrix{
\lasu_{\partial\Delta^n}\ar[d]\ar[r]&\lasu_{\Delta^n}\ar[d]\\
L\ar[r]&(L\amalg \libc(x),d)
}
$$
where the coproduct is taken in the category of complete Lie algebras, and $d$ restricts to the given differential in $L$. Additionally,   $x$ corresponds to the (desuspension) of the top simplex of $\Delta^n$. Hence, it has degree $n-1$, it is a Maurer-Cartan element if $n=0$, and is either a cycle or a path between two Maurer-Cartan elements of $L$ if $n=1$.

The pushout of the bottom map along any morphism $\varphi\colon L\to M$ produces
$$
\xymatrix{
L\ar[d]_\varphi\ar[r]&(L\amalg \libc(x),d)\ar[d]^{\overline\varphi}\\
M\ar[r]&(M\amalg \libc(x),d).
}
$$
We are thus compelled to prove that $\overline\varphi$ is a weak equivalence if $\varphi$ is.  

Begin by choosing  cofibrant replacements of $\libc(V)\stackrel{\sim}{\to} L$ and $\libc(W)\stackrel{\sim}{\to}M$ where $V$ and $W$ are as in Proposition \ref{cofire}, see also Definition \ref{modeloretractivo}. We then construct in the natural way  weak equivalences
$$
\libc(V)\amalg\libc(\overline L\oplus d\overline L)\stackrel{\sim}{\longrightarrow} L\quad\text{and}\quad\libc(W)\amalg\libc(\overline M\oplus d\overline M)\stackrel{\sim}{\longrightarrow} M,
$$
where $\overline L=L_{\ge0}$ and $\overline M=M_{\ge0}$. These also serve as  cofibrant replacements of $L$ and $M$ respectively, according to Proposition 8.10 of \cite{bufemutan0}, and they are fibrations since they are  surjective in non-negative degrees.  This yields the following commutative diagram
$$
\xymatrix{
\libc({\mathcal V})\ar[d]_\sim\ar[r]^\sim&L\ar[d]^{\varphi}_\sim\\
\libc({\mathcal W})\ar[r]_\sim&M
}
$$
where, for convenience, we denote ${\mathcal V}=V\oplus\overline L\oplus d\overline L$ and ${\mathcal W}=W\oplus\overline M\oplus d\overline M$. 

Next we check that there is a cdgl extension $\libc({\mathcal V})\to \libc({\mathcal V})\amalg\libc(x)$ of $\libc({\mathcal V})$ (the same holds for $\libc({\mathcal W})$) for which the following diagram is a pushout:
$$
\xymatrix{
\libc({\mathcal V})\ar[d]_\sim\ar[r]&\libc({\mathcal V})\amalg\libc(x)\ar[d]\\
L\ar[r]&L\amalg\libc(x).
}
$$
Indeed this is trivial if $x$ is an MC element  or if $x$ is a cycle of degree $0$. If $x$ is a path between two Maurer-Cartan elements of $L$, as $\libc({\mathcal V})\to L$ is surjective in degree $0$, we may apply \cite[Proposition 5.18]{bufemutan0} to find a path in $\libc({\mathcal V})$ which is sent to the path in $L$. Define $dx\in \libc({\mathcal V})$ accordingly. Finally if $|x|>0$, as $\libc({\mathcal V})\to L$ is surjective in degree zero and an isomorphism in homology for non-negative degrees, we may also find an element $\gamma\in \libc({\mathcal V})$ which is sent to $dx\in L$. In $\libc({\mathcal V})\amalg\libc(x)$, define $dx=\gamma$.

Now, recall that left properness holds whenever the domain and codomain of the cofibration  are cofibrant \cite[Prop.~13.1.2(1)]{hirsch0}. Hence, in the following commutative diagram,
$$
\xymatrix@R=8pt@C=8pt{
    \libc({\mathcal V})\ar[rr] \ar[dd]_\sim \ar[dr]^\sim & & \libc({\mathcal V})\amalg\libc(x) \ar@{-}[d]_(.65)\sim \ar[dr]^\sim \\
    & L \ar[rr] \ar[dd]_(.35)\sim^(.35)\varphi & \ar[d]& L \amalg\libc(x) \ar[dd]^{\overline\varphi} \\
     \libc({\mathcal W})\ar@{-}[r] \ar[dr]_\sim & \ar[r]& \libc({\mathcal W})\amalg\libc(x)  \ar[dr]^\sim \\
    & M \ar[rr] & &  M\amalg\libc(x)
}
$$
the three maps in the right square, and thus  $\overline\varphi$ as well,  are weak equivalences.
 \end{proof}

\subsubsection{Quillen endofunctors in $\cdgll$}
To construct the  endofunctors needed to develop spectra in $\cdgll$, we begin by extending  the path object functor of $\catcdgl$, as  defined in \cite[\S8.3]{bufemutan0}, to  $\cdgll$.

Given a cdgl  $L$, with associated filtration $\{F^n\}_{n\ge 1}$, consider
$$
L^I=L\otimesc \land(t,dt)=\varprojlim_n \bigl(L/F^n\otimes \land(t,dt)\bigr)
$$
where $t$ has degree $0$. Observe that an element of $L^I$ can be written as a series 
\begin{equation}\label{series}
\sum_{j\ge 0} x_j t^j+\sum_{k\ge 0} y_k t^kdt,\quad x_j,y_k\in L,
\end{equation}
where $\sum_{j\ge0} x_j$ and $\sum_{k\ge0}y_k$ are well defined elements in $L$. In particular, the maps
$$
\varepsilon_i\colon L^I\longrightarrow L,\quad t\mapsto i,\quad i=0,1,
$$
are well defined cdgl morphisms.

\begin{definition}\label{pathobjectd}
Let $M\in\cdgll$. Define $M^I_L$ as the pullback
$$
\xymatrix{
M^I_L\ar[d]\ar[r]&L\ar[d]\\
M^I\ar[r]^{p^I}&L^I,
}
$$
in which the right vertical morphism is the natural inclusion $a\mapsto a\otimes 1$ and  $p$ denotes the retraction of $M$. Observe that $M^I_L$ is the sub cdgl of $M^I$ formed by the series
$$
x_0+\sum_{j\ge 1} x_j t^j+\sum_{k\ge 0} y_k t^kdt,\quad x_0\in L,\quad x_j,y_k\in K_M.
$$
In other words, the splitting (\ref{splitting}) results in
\begin{equation}\label{splitl}
M^I_L\cong L\oplus \bigl(K\otimesc \land (t,dt)\bigr).
\end{equation}
Note also that $M^I_L$ is an object of  $\cdgll$ with section and retraction given by the top sequence of the following diagram in which both squares are pullbacks and $s$ denotes the section of $M$,
$$
\xymatrix{L\ar[r]\ar[d]&
M^I_L\ar[d]\ar[r]&L\ar[d]\\
L^I\ar[r]^{s^I}&M^I\ar[r]^{p^I}&L^I.
}
$$
\end{definition}
Consider the following maps in which the left one is again the natural inclusion
\begin{equation}\label{pathobject}
M\longrightarrow M^I_L\stackrel{\scriptscriptstyle(\varepsilon_0,\varepsilon_1)}{\longrightarrow} M{\times_L}M
\end{equation}
Then:

\begin{proposition}\label{pathobjectp}
This sequence constitutes a functorial good path object in $\cdgll$.
\end{proposition}
\begin{proof}
We first check that $M\stackrel{\sim}{\to} M^I_L$, $x\mapsto x\otimes 1$, is a weak equivalence. To see this, recall (see for instance \cite[Prop.~3.4]{bufemutan0})  that if $M=\varprojlim_n M/F^n$, then $M^I$ is complete with respect to the filtration
$$
G^n=\ker\bigl(M^I\longrightarrow M/F^n\otimes \land (t,dt) \bigr),\quad n\ge 0.
$$
As a result $M^I_L$ is complete for the filtration $\{G^n\cap M^I_L\}_{n\ge 1}$. In view of (\ref{splitl}), one verifies that, as vector spaces,
$$
(G^n\cap M^I_L)/(G^{n+1}\cap M^I_L)\cong F^n/F^{n+1}\oplus \bigl((K\cap F^n)/(K\cap F^{n+1})\otimes\land(t,dt)\bigr).
$$
Under this identification, the morphism induced on the associated graded  spaces by  the map $M\to M^I_L$,
$$
F^n/F^{n+1}\longrightarrow F^n/F^{n+1}\oplus \bigl((K\cap F^n)/(K\cap F^{n+1})\otimes\land(t,dt)\bigr),
$$
is simply the inclusion into the first summand. This is a quasi-isomorphism since the second summand is acyclic. Then, by the {\em Goldman-Millson Theorem} (see \S\ref{homotocdgl}), the map $M\stackrel{\sim}{\to} M^I_L$ is a weak equivalence in $\catcdgl$, and consequently,  it is also a weak equivalence in $\cdgll$.

On the other hand, the composition of the maps in (\ref{pathobject}) is  the diagonal in $\cdgll$. Finally, both $\varepsilon_i$, $i=0,1$, are surjective maps and therefore fibrations. We also remark that, by the ``2 out of 3'' property, these maps are also weak equivalences
\end{proof}

\begin{definition}\label{loopsl} The {\em retractive loop} on a retractive cdgl $M$, denoted $\lupl M$, is  the retractive cdgl defined by the pullback
$$
\xymatrix{
\lupl M\ar[d]\ar[r]&M^I_L\ar[d]^{(\varepsilon_0,\varepsilon_1)}\\
L\ar[r]_(.40){\Delta{\circ}s} &M{\times_L}M.
}
$$
When $L=0$ we simply write $\lup M$.
\end{definition}
In view of Proposition \ref{pathobjectp} we deduce:

\begin{corollary}\label{loopobject}
$\lupl$ is a functorial loop object in $\cdgll$.\hfill$\square$
\end{corollary}

The following explicit description of the functor $\lupl$ will be used extensively in the sequel.
\begin{proposition}\label{lupquasi}
Let $M\cong L\oplus K$ a retractive cdgl. Then, there is a natural isomorphism of retractive cdgl's,
$$
\lupl M\cong L\oplus \bigl(K\otimesc (C\oplus dC\oplus \bq dt)\bigr)
$$
\end{proposition}
\begin{proof}
By Definitions \ref{pathobjectd} and \ref{loopsl}, an arbitrary element of $\lupl M$ can be written as a series 
$$
x_0+\sum_{j\ge 1} x_j t^j+\sum_{k\ge 0} y_k t^kdt,\quad x_0\in L,\quad x_j,y_k\in K_M,
$$
where the series $\sum_{j\ge 1}x_j$ has $0$ and $1$ as roots. In other words, 
$$
\lupl M\cong L\oplus \Bigl(K\otimesc \bigl((t^2-t)\bq[t]\oplus dt\bq[t]\bigr)\Big).
$$
But
$$
dt\,\bq [t]\cong (2t-1)\bq[t]dt\oplus \bq dt
$$
so that
$$
\lupl M\cong L\oplus \bigl(K\otimesc (C\oplus dC\oplus \bq dt)\bigr)
$$
with $U=(t^2-t)\bq[t]$. 
\end{proof}

\begin{corollary}\label{lupquasicoro}
There is a map of retractive cdgl's
$$
\lupl M\stackrel{\simeq}{\longrightarrow} L\oplus s^{-1}K
$$
that is both a quasi-isomorphism and a weak equivalence. Here, $s^{-1}K$ is an abelian sub cdgl, with differential $s^{-1}d$ given by the desuspension of the differential in $K$ and 
$$
[x,s^{-1}y]=(-1)^{|x|}s^{-1}[x,y],\quad x\in L,\quad y\in K.
$$
Moreover, this weak equivalence  admits a section.
\end{corollary}

\begin{proof}
Using the notation from the previous result,  the projection $C\oplus dC\oplus \bq dt \stackrel{\simeq}{\to} \bq dt$ induces a quasi-isomorphism of retractive cdgl's
$$
\lupl M\cong L\oplus \bigl(K\otimesc (C\oplus dC\oplus \bq dt)\bigr)\stackrel{\simeq}{\longrightarrow} L\oplus (K\otimesc \bq dt).
$$
But 
$$
L\oplus (K\otimesc \bq dt)= L\oplus (K\otimes \bq dt)\cong L\oplus s^{-1}K.
$$ 
To show that this map  is also a weak equivalence, consider  filtrations $\{F_n\}_{n\ge0}$ and $\{G_n\}_{n\ge 0}$ with respect to  which $L$ and $K$ are complete, and filter the domain and codomain of the preceding quasi-isomorphism  accordingly. The  induced map on the associated graded spaces is then
$$
F^n/F^{n+1} \oplus \bigl(G^n/G^{n+1}\otimesc(C\oplus dC\oplus \bq dt)\bigr)\stackrel{\simeq}{\longrightarrow} F^n/F^{n+1} \oplus (G^n/G^{n+1}\otimesc \bq dt),
$$
which is also a quasi-isomorphism. The assertion then follows from the {\em Goldman-Millson Theorem} (see \S\ref{homotocdgl}). 

Finally, the section is induced by the inclusion $\bq dt\stackrel{\simeq}{\to} C\oplus dC\oplus \bq dt $.
\end{proof}

\begin{proposition}\label{lup*}For any cdgl morphism $f\colon L\to L'$ there is a natural isomorphism
$$
f^*\Omega_{L'}\cong\lupl f^*.
$$
\end{proposition}
\begin{proof}
We first see that for each $M\in\cdgll$ there is a natural isomorphism
\begin{equation}\label{naturi}
 (f^*M)^I_L\cong f^*(M^I_{L'}).
\end{equation}
For it recall that  $f^*M$ is given by the pullback
$$
\xymatrix{
f^*M\ar[d]\ar[r]&L\ar[d]^f\\
M\ar[r]^{p}&L'.
}
$$
As noted in the proof of Proposition \ref{limite}, $(\,\cdot\,)^I$ preserves limits, so that 
$$
\xymatrix{
(f^*M)^I\ar[d]\ar[r]&L^I\ar[d]^(.45){f^I}\\
M^I\ar[r]^{p^I}&{{L'}^I}
}
$$
is also a pullback. Consider then the commutative cube
$$
\xymatrix@R=8pt@C=8pt{
    (f^*M)^I_L\ar[rr] \ar[dd] \ar[dr] & & L \ar@{-}[d]\ar[dr]^f \\
    & M^I_{L'} \ar[rr] \ar[dd] & \ar[d]& L'  \ar[dd] \\
    (f^*M)^I\ar@{-}[r] \ar[dr] & \ar[r]& L^I  \ar[dr]^(.35){f^I} \\
    & M^I \ar[rr] & &  {L^\prime}^I
}
$$
in which the back, bottom and front faces are pullbacks. Hence, the top face is also a pullback and thus  $(f^*M)^I_L\cong f^*(M^I_{L'})$ as required.
 
Next, observe that for any map $N\to N'$ in $\cdgll$ there is a pullback
$$
\xymatrix{
f^*N\ar[d]\ar[r]&^N\ar[d]\\
f^*N'\ar[r]&N'.
}
$$
Thus, in the commutative cube
$$
\xymatrix@R=6pt@C=6pt{
    f^*(\Omega_{L'}M)\ar[rr] \ar[dd] \ar[dr] & & f^*(M_{L'}^I) \ar@{-}[d]\ar[dr] \\
    & \Omega_{L'}M \ar[rr] \ar[dd] & \ar[d]& M_{L'}^I  \ar[dd] \\
    L\ar@{-}[r] \ar[dr] & \ar[r]& f^*(M{\times_{L'}} M)\ar[dr] \\
    & L' \ar[rr] & &  M{\times_{L'}} M,
}
$$
the front and the two lateral faces are pullbacks. Hence, the back face is also a pullback. However, by (\ref{naturi}) and the obvious natural isomorphism 
$$
f^*(M{\times_{L'}} M)\cong f^*M{\times_L}f^*M,
$$
the back face has the form
$$
\xymatrix{
f^*(\Omega_{L'}M)\ar[d]\ar[r]&(f^*M)^I_L\ar[d]\\
L\ar[r]&f^*M{\times_L}f^*M
}
$$
and we deduce that
$$
f^*(\Omega_{L'}M)\cong \lupl(f^*M).
$$
\end{proof}

\begin{proposition}\label{limite} the functor $\lupl$ preserves (small) limits.
\end{proposition}
\begin{proof}
We first see that, in $\catcdgl$, the functor $(\,\cdot\,)^I$ is continuous by checking that it commutes with (small) products and binary equalizers. Let $\{L_i\}_{i\in I}$ be a family of cdgl's each of which is filtered by $\{F_i^n\}_{n\ge 1}$. On the one hand, 
$$
\textstyle\prod_{i\in I}(L_i^I)=\textstyle\prod_{i\in I}\Bigl(\varprojlim_n \bigl(L_i/F_i^n\otimes \land(t,dt)\bigr)\Bigr)\cong \varprojlim_n\bigl(\textstyle\prod_{i\in I}\bigl( L_i/F_i^n\otimes \land(t,dt)\bigr)\Bigr).
$$
On the other hand, as $\textstyle\prod_{i\in I}L_i$ is complete respect to the filtration $\{\textstyle\prod_{i\in I}F_i^n\}_{n\ge 1}$,
$$
(\textstyle\prod_{i\in I} L_i)^I=\varprojlim_n\bigl((\textstyle\prod_{i\in I}L_i/F_i^n)\otimes\land(t,dt)\bigr).
$$
Now, for each $n\ge 1$, the natural map induced by the projections $\prod_{i\in I}(L_i/F_i^n)\to L_i/F_i^n$,
$$
(\textstyle\prod_{i\in I}L_i/F_i^n)\otimes\land(t,dt)\to \textstyle\prod_{i\in I}\bigl( L_i/F_i^n\otimes \land(t,dt)\bigr)
$$
fails to be surjective in general. This is fixed by taking inverse limits, and we have a bijection
$$
\varprojlim_n\bigl((\textstyle\prod_{i\in I}L_i/F_i^n)\otimes\land(t,dt)\bigr)\stackrel{\cong}{\longrightarrow} \varprojlim_n\bigl(\textstyle\prod_{i\in I}\bigl( L_i/F_i^n\otimes \land(t,dt)\bigr)\Bigr).
$$
Indeed, any element in the codomain can be written as a series like in (\ref{series}), but with coefficients now belonging to $\textstyle\prod_{i\in I}L_i$.

Moreover, the representation (\ref{series}) of a generic element in the complete tensor product implies that $(\,\cdot\,)^I$ preserves equalizers. Consequently, this functor preserves (small) limits. 

Thus, by Definition \ref{pathobjectd}, the same conclusion applies to the functor $(\,\cdot\,)^I_L$ in $\cdgll$. Once again, by Definition \ref{loopsl}, and noting that the product in $\cdgll$ commutes with limits, the functor $\lupl$ preserves limits as well.
\end{proof}

By the previous proposition and noting, for instance, that $\cdgll$ is bicomplete and locally presentable, the {\em special adjoint functor theorem}  justifies the following:

\begin{definition}\label{tsuspension}
The retractive loop functor $\lupl$ admits a left adjoint denoted by $\susl$ and referred to as the {\em retractive suspension functor}. When $L=0$ we write it simply as $\Sigma$.
\end{definition}

 We can explicitly describe the suspension functor up to weak equivalence via its left derived functor in the homotopy category. To this end let $(L\amalg\libc(W),d)$ be a retractive model of a given  $M\in \cdgll$ which, in view of  equality (\ref{esencialk}) of Definition \ref{libre}, can be written as 
 $$
 (L\amalg\libc(W),d)\cong L\oplus (\libc(T),d).
 $$
 
\begin{theorem}\label{suspensionex} The retractive suspension $\Sigma_LM$ is weakly equivalent to 
$$
L\oplus (\libc(sT),d_1)
$$
where $d_1$ denotes the suspension of the linear part of $d$:
$$
d_1sx=-sd_1x,\quad x\in T.
$$
Equivalently, this retractive cdgl is isomorphic to
$$
(L\amalg \libc(sW),\partial)
$$
where $\partial$ is retractive linear: for each $w\in W$, each summand of $\partial sw$ is a bracket containing a single element of $sW$.
\end{theorem}

\begin{proof}
Let $N=L\oplus K\in\cdgll$. By Corollary \ref{lupquasicoro}, $\lupl N$ is weakly equivalent to  $L\oplus s^{-1}K$ and, since $s^{-1}K$ is an abelian cdgl, any morphism $\varphi\colon (\libc(T),d)\to  s^{-1}K$ is uniquely determined by the  morphism 
\begin{equation}\label{formula1}
\psi\colon (\libc(sT),d_1)\longrightarrow  K,\quad \psi(sx)=y,\quad \text{with}\quad s^{-1}x=\varphi(x).
\end{equation}
That is, there is a natural bijection
\begin{equation}\label{formula2}
\Hom_{\cdgll}(L\oplus(\libc(sT),d_1),L\oplus K)\cong \Hom_{\cdgll}(L\oplus(\libc(T),d),L\oplus s^{-1}K)
\end{equation} which induces 
$$
\Hom_{\Ho\cdgll}(\susl M,N)\cong \Hom_{\Ho\cdgll}(M\lupl N)
$$
and the first claim follows.

On the other hand, consider the isomorphism of  graded Lie algebras provided by (\ref{esencialk})
$$
L\oplus \libc(sT)\stackrel{\cong}{\longrightarrow} L\amalg\libc(sW)
$$
which is the identity on $L$ and sends a generator $s\bigl[a_1,[a_2,[\dots[a_n,w]\bigr]\dots\bigr]$ of $sT$ to $\pm \bigl[a_1,[a_2,[\dots[a_n,sw]\bigr]\dots\bigr]$. If we equip $L\amalg\libc(W)$ with the differential $\partial$ making it isomorphic to $L\oplus (\libc(sT),d_1)$ as retractive cdgl's, then the second assertion follows.
\end{proof}
 
With the notation of the previous result:

\begin{corollary}\label{nosesi}
For any $n\ge1$, $\susl^nM$ is weakly equivalent to 
$$
L\oplus (\libc(s^nT),d_1)\cong (L\amalg\libc(s^nW),\partial)
$$
where 
$$
d_1(s^nx)=(-1)^{n}s^{n}d_1(x)\quad\text{and}\quad \partial(s^nw)=(-1)^{n-1}s^{n-1}\partial(sw).
$$
 In particular if $M=L\oplus K$
$$
\susl^n(M)\quad\text{is weakly equivalent to}\quad L\oplus \Sigma^nK.
$$
\hfill$\square$
\end{corollary}

\begin{rem}\label{ya vere}
Although implicit in the proof of Theorem \ref{suspensionex}, it is worth noting that for retractive cdgl's $M$ and $N=L\oplus K$,  any map $\susl M\to N$ is weakly equivalent  to a map $\psi\colon (L\oplus (\libc(sT),d_1)\to N$ whose adjoint is weakly equivalent to  $\varphi\colon L\oplus (\libc(T),d)\to L\oplus s^{-1}K$, defined by $\varphi(x)=s^{-1}y$ where $\psi(sx)=y$, for $x\in T$.
\end{rem}

We conclude with:

\begin{proposition}\label{endol}  $\susl{\dashv}\,\, \lupl$ form a Quillen pair of endofunctors.
\end{proposition}

\begin{proof}
Observe that $\lupl$ preserves fibrations. Moreover, since $(\,\cdot\,)^I_L$ provides a functorial path object, $\lupl$  also preserves weak equivalences between fibrant objects. Since every object in $\cdgll$ is fibrant this establishes the proposition.  Note that this reasoning remains valid for any general model category in which every object is fibrant and equipped with a limit-preserving functorial path object functor. Nevertheless, for completeness, we proceed with the argument in our specific case:

 If $f\colon M\stackrel{\sim}{\to} N$ is a weak equivalence in $\cdgll$ so is $f^I_L$ in view of the commutative squares, for  $i=0,1$, 
$$
\xymatrix{
M^I_L\ar[d]_{\varepsilon_i}^\sim\ar[r]^{f^I}&N^I_L\ar[d]^{\varepsilon_i}_\sim\\
M\ar[r]_{f}^\sim &N.
}
$$

On the other hand, for any $M\in\cdgll$, write $\lupl M$ as the result of the successive pullbacks
$$
\xymatrix{ P_M\ar[d]\ar[r]& M^I_L \ar[d]^{\varepsilon_0} & \quad & \lupl M\ar[d]\ar[r]&P_M \ar[d]^{\varepsilon_1} \\ L\ar[r]_s&M, & \quad & L\ar[r]_s&M. } 
$$
Then, in the commutative cubes,
$$
\xymatrix@R=10pt@C=10pt{
    P_M\ar[rr] \ar[dd] \ar[dr]^\sim & & M^I_L \ar@{-}[d]_(.70){\varepsilon_0}\ar[dr]^{f^I}_\sim \\
    & P_N \ar[rr] \ar[dd]& \ar[d]& N^I_L \ar[dd]^(.42){\varepsilon_0}\\
    L\ar@{-}[r] \ar@{=}[dr] & \ar[r]& M  \ar[dr]^(.40)f_\sim \\
    & L \ar[rr] & &  N
}
\qquad \xymatrix@R=10pt@C=10pt{
    \lupl M\ar[rr] \ar[dd] \ar[dr]_\sim^(.60){\lupl f} & & P_M \ar@{-}[d]_(.70){\varepsilon_1}\ar[dr]^\sim \\
    & \lupl N \ar[rr] \ar[dd]& \ar[d]& P_N \ar[dd]^(.42){\varepsilon_1} \\
    L\ar@{-}[r] \ar@{=}[dr] & \ar[r]& M  \ar[dr]_(.40)\sim^(.40)f \\
    & L \ar[rr] & &  N,
}
$$
the map $P_M\stackrel{\sim}{\to} P_N$ is a weak equivalence and thus so is $\lupl f$. For the last assertion see for instance  the dual of Corollary to Theorem B in \cite{ree} which, to the best of our knowledge, is the original source.
\end{proof}

\subsection{Rational parametrized spectra are parametrized cdgl spectra}

\subsubsection{Spectra of retractive cdgl's}

To stabilize the proper and combinatorial model category $\cdgll$ we follow the general framework in \S\ref{specmod} using the Quillen pair of endofunctors  ${\susl}\dashv\,\lupl$.

\begin{definition}\label{falta} Let $L$ be a cdgl. A {\em spectrum} in $\cdgll$ or {\em $L$-spectrum} is an object  of the category $\spec_L(\catcdgl)$, see Definition \ref{specret}, which we denote by  $\spec_L$ henceforth. In other words, an $L$-spectrum is a sequence\footnote{We have chosen to write superscripts to avoid confusion with the homological degree of a given cdgl.} $M=\{M^n\}_{n\ge0}$  in $\cdgll$,
$$
L\stackrel{s^n}{\longrightarrow} M^n \stackrel{p^n}{\longrightarrow}L,
$$
 endowed with  {\em structure morphisms}, also in $\cdgll$,
$$
\sigma\colon \susl M^n\longrightarrow M^{n+1}\quad\text{or equivalently}\quad \eta\colon M^n\longrightarrow \lupl M^{n+1},\quad n\ge 0.
$$
A morphism  $f\colon M\to N$ of $\spec_L$ is a family $\{f^n\colon M^n\to N^n\}_{n\ge0}$ of maps in $\cdgll$ compatible with the structure morphisms.

An $L$-spectrum  is said to be {\em free} if it is levelwise free. That is, it has the form $\{(L\amalg\libc(W^n),d)\}_{n\ge 0}$. We  denote such spectrum simply by $(L\amalg\libc(W),d)$. 
We say that a free $L$-spectrum $(L\amalg\libc(W),d)$ is {\em linear} if its differential is levelwise linear in the sense of Definition \ref{libre}.
\end{definition}
 
 The following definitions are of special relevance in what follows.
 
 \begin{definition}\label{indelineal} {\em (The linear and indecomposable reduction of a spectrum)} Let $M=(L\amalg\libc(W),d)=\{(L\amalg\libc(W^n),d)\}_{n\ge 0}$ be a free $L$-spectrum and, as usual, decompose each level as
 $$
( L\amalg\libc(W^n),d)\cong L\oplus(\libc(T^n),d),\quad n\ge 0.
 $$

 \smallskip
 
 (i) By taking  the retractive linear part of the differential on each level, see Definition \ref{libre}, we obtain a family of retractive cdgl's 
$$
(L\amalg\libc(W^n),d_1)\cong L\oplus ( \libc(T^n),d_1),\quad n\ge0.
$$
The adjoint of the $n$th structure map of $M$,
$$
\eta_n\colon M^n= L\oplus (\libc(T^n),d)\longrightarrow L\oplus (\libc(T^{n+1}),d)\otimesc (C\oplus dC\oplus \bq dt)=\lupl M^{n+1},
$$
restricts to a morphism between the corresponding ideals
\begin{equation}\label{lasmu}
 (\libc(T^n),d)\longrightarrow (\libc(T^{n+1}),d)\otimesc (C\oplus dC\oplus \bq dt)
\end{equation}
to which we may apply the following general fact: 

Any cdgl morphism of the form 
$$
\varphi\colon (\libc(V),d)\longrightarrow (\libc(W),d)\otimesc A,
$$
with $A$ is a commutative differential graded algebra, 
induces  a morphism 
$$
\varphi^1\colon (\libc(V),d_1)\longrightarrow (\libc(W),d_1)\otimesc A\quad\text{such that}\quad
\varphi^1(V)\subset W\otimes A.
$$
In our case this yields a morphism
$$
 (\libc(T^n),d_1)\longrightarrow (\libc(T^{n+1}),d_1)\otimesc (C\oplus dC\oplus \bq dt)
$$
which extends to a map
$$
\eta_n^1\colon L\oplus \libc(T^n,d_1)\to \lupl \bigl(L\oplus\libc(T^{n+1},d_1)\bigr).
$$
The {\em linear reduction of $M$} is the $L$-spectrum,
$$
M_{\rm lin}=(L\amalg \libc(W),d_1)=\{(L\amalg \libc(W^n),d_1)\}_{n\ge 0}
$$
equipped with the adjoint structure maps $\{\eta_n^1\}_{n\ge0}$.

 \smallskip
 
 (ii) On the other hand, 
by projecting every $(L\amalg \libc(W^n),d)$ onto the indecomposables of $\libc(W^n)$ we obtain a family of retractive cdgl's $\{(L\amalg W^n,d_1)\}_{n\ge 0}$ in which each $W^n$ is an abelian cdgl. We then have,
$$
(L\amalg W^n,d_1)\cong L\oplus (T^n,d_1),\quad n\ge 0.
$$
As before, each  morphism in  (\ref{lasmu}) induces a map
$$
\bar\mu_n\colon (T^n,d_1)\to (T^{n+1},d_1)\otimesc (C\oplus dC\oplus \bq dt)
$$
which extends to a morphism 
$$
\bar\eta_n\colon L\oplus (T^n,d_1)\to \lupl \bigl(L\oplus(T^{n+1},d_1)\bigr).
$$
The {\em indecomposable reduction of $M$} is defined as the $L$-spectrum,
$$
M_{\rm ind}=(L\amalg W,d_1)=\{(L\amalg W^n,d_1)\}_{n\ge 0}
$$
with adjoint structure maps $\{\bar\eta_n\}_{n\ge0}$.

 \smallskip
 
Note that, in both cases,  projecting over the indecomposables yields maps of $L$-spectra
\begin{equation}\label{reducciones}
M\longrightarrow M_{\rm ind}\longleftarrow M_{\rm lin}
\end{equation}
\end{definition}

Sections \S\ref{specmod} and \S\ref{apret}, when specialized to $\cdgll$, provide the following facts:

\begin{proposition}\label{generalcdgl} 
There are adjoint endofunctors 
\begin{equation}\label{boul}
\xymatrix{ \specl & \specl\ar@<0.75ex>[l]^(.50){\lupel}
\ar@<0.75ex>[l];[]^(.50){\susel}\\}
\end{equation}
obtained by prolongation of  $\susl$ and $\lupl$. Moreover, for any $k\ge 0$, there are adjoint functors
$$
\xymatrix{ \cdgll & \specl\ar@<0.75ex>[l]^(.38){\eva_k}
\ar@<0.75ex>[l];[]^(.50){\susl^{\infty-k}}\\}
$$
provided by {\em (\ref{eva})} and {\em (\ref{susk})}.\hfill$\square$
\end{proposition}

\begin{definition}\label{projectivecdgl} 
(i) The fibrations and weak equivalence of the {\em projective model structure}  on $\specl$ are, as given in Definition \ref{proestable}, levelwise fibrations and weak equivalences in $\cdgll$. Cofibrations are defined accordingly. In particular, see (\ref{generaes}), if we consider the sets
$$
\mathcal I_L=\{L\amalg i,\,i\in \lasu(I)\}\quad\text{and}\quad \mathcal J_L=\{L\amalg j,\,j\in \lasu(J)\},
$$
with $\lasu(\mathcal I)$ and $\lasu(\mathcal J)$ as in \S\ref{homotocdgl}, of generating cofibration and trivial cofibrations of $\cdgll$, then the families

$$
\mathcal I_{\susl}=\cup_{k\ge 0}\susl^{\infty-k}(\mathcal I_L)\esp{and}\mathcal J_{\susl}=\cup_{k\ge 0}\susl^{\infty-k}(\mathcal J_L)
$$
form generating cofibrations and trivial cofibrations of $\specl$, respectively.

 (ii) On the other hand, the {\em stable model structure} in $\specl$ is then the Bousfield localization of the projective model structure with respect to the family $S$ in (\ref{equibous}). 
 \end{definition}

\begin{rem}\label{cosascdgl} By Corollary \ref{changeret} the pair (\ref{boul}) becomes a Quillen equivalence in the stable structure of $\specl$, in which the fibrant objects are the $\lupl$-spectra, i.e., those $M=\{M^n\}_{n\ge 0}$ for which the map $M^n\to\lupl M^{n+1}$ is a weak equivalence. 
Additionally, a morphism 
 between $\lupl$-spectra  is a stable equivalence if and only if it is a levelwise weak equivalence.
\end{rem}

On the other hand, given any cdgl morphism $f\colon L\to L'$, Proposition \ref{lup*}  and Corollary \ref{changeret}  imply:

\begin{proposition}\label{quillenchange} 
There is a Quillen pair with respect to the stable structure,
\begin{equation}\label{changespecl}
\xymatrix{ \specl & \spec_{L'} \ar@<0.75ex>[l]^(.50){\widetilde{f^*}}
\ar@<0.75ex>[l];[]^(.50){\widetilde{f_!}}\\}
\end{equation}
where $\widetilde{f_!}$ and $\widetilde{f^*}$ are the prolongations of $f_!$ and $f^*$ respectively. Additionally, this is a Quillen equivalence if $f$ is a weak equivalence.\hfill$\square$
\end{proposition}

\begin{definition}\label{modelo} A {\em model} of an $L$-spectrum $N$ is a free $L$-spectrum $M$ together with a stable equivalence
$$
M\stackrel{\sim}{\longrightarrow}N.
$$
We also refer to the free spectrum $M$ itself as a model of $N$.

The linear reduction $M_{\rm lin}$ and the indecomposable  reduction $M_{\rm ind}$ of $M$ are called,   respectively, a {\em linear model} and an {\em indecomposable model}  of $N$.
\end{definition}

 \begin{theorem}\label{freespec} Every $L$-spectrum $N$ admits a model $M$ which can be chosen to be projectively equivalent to $N$. Such a model is unique up to stable equivalence.
\end{theorem}
\begin{proof}
Choose a retractive model for each level of $N$,
$$
\varphi_n\colon (L\amalg \libc(W^n),d)\stackrel{\sim}{\longrightarrow} N^n,\quad n\ge0.
$$ 
 By the standard argument (see the proof of Theorem \ref{combi}) this map can be chosen to be surjective in non-negative degrees, that is, a fibration in $\cdgll$. As $\susl$ is  left Quillen it preserves cofibrant objects and acyclic cofibrations. We thus obtain a map $\gamma_n$ that completes the following commutative square, in which $\sigma_n$ is the structure map and the vertical arrows are weak equivalences:
$$
\xymatrix{
\susl(L\amalg \libc(W^n),d)\ar[d]_{\susl\varphi_n}^\sim\ar[r]^{\gamma_n}&(L\amalg \libc(W^{n+1}),d) \ar[d]^{\varphi_{n+1}}_\sim\\
\susl N^{n}\ar[r]_{\sigma_n} &N^{n+1}.
}
$$
In other words, given the $L$-spectrum 
$
M=\{(L\amalg \libc(W^n),d)\}_{n\ge 0}$, with structure maps $\{\gamma_n\}_{n\ge0}$,
 we have a map $
\varphi\colon M{\longrightarrow} N
$ in $\specl$
which is a projective, and hence stable, equivalence.
\end{proof}

The following will be shown to define a fundamental stable invariant of $\specl$, see Theorem \ref{homotohomolo} and Corollary \ref{sorpresa}:

\begin{definition}\label{sorpresa2} ({\em Stable homology}) Let $M=\{M^n\}_{n\ge0}\in\specl$ with $M^n=L\oplus K^n$. Composing the  adjoint structure map 
$$
\eta_n\colon M^n\longrightarrow\lupl M^{n+1}
$$
with the weak equivalence of Corollary \ref{lupquasicoro}, and restricting to the corresponding ideals, yields a cdgl morphism
$$
\mu_n\colon K^n\longrightarrow s^{-1}K^{n+1}
$$
with abelian codomain. This defines a cdgl sequence
\begin{equation}\label{secuencia}
K^0\stackrel{\mu_0}{\longrightarrow} s^{-1}K^1\to\dots\to s^{-n}K^n\stackrel{s^{-n}\mu_n}{\longrightarrow}s^{-n-1}K^{n+1} \dots
\end{equation}
and the {\em stable homology of $M$} (at 0) is defined by 
$$
H^{\rm st}(M)=H(\varinjlim_n s^{-n}K^n).
$$
Any map $f\colon M\to N$ in $\specl$ induces a morphism between the corresponding sequences, and hence a map
$$
H^{\rm st}(f)\colon H^{\rm st}(M)\longrightarrow H^{\rm st}(N).
$$
Since $
H_k(s^{-n}K^n)\cong H_{n+k}(K^n)
$ and homology commutes with directed colimits,
 we have
 $$
H^{\rm st}(M)=\oplus_{k\in\bz}H_k^{\rm st}(M)
\quad\text{where}\quad
H^{\rm st}_k(M)=\varinjlim_{n} H_{n+k}(K^n).
$$
Now if $a\in\mc(L)$ we may consider the perturbed differential $d_a$ on each $K^n$ and define the {\em stable homology of $M$ at $a$} as
$$
H^{\rm st}(M,a)=H\bigl(\varinjlim_n (s^{-n}K^n,s^{-n}d_a)\bigr).
$$
Again,
$$
H^{\rm st}(M,a)=\oplus_{k\in\bz}H_k^{\rm st}(M,a)\quad\text{where}\quad
H^{\rm st}_k(M,a)=\varinjlim_{n} H_{n+k}(K^n,d_a).
$$
\end{definition}

\begin{rem}\label{remarkestable} Since for each  $n\ge k+1$ we have $H_{n+k}(K^n,d_a)=H_{n+k}(K^n,d_a)^{(0)}$, and $(K^n,d_a)^{(0)}=(K^n)^{(a)}$ for all $n\ge 0$, we  are free to use any of these connected cdgl's to compute the stable homology of $M$ at $a$. 
\end{rem}

Unlike the unstable setting, yet consistent with the simplifications that arise in the stable context, we have:

\begin{proposition}\label{estahomolin} Every  $L$-spectrum is stably equivalent to  its linear and indecomposable  models. 
\end{proposition}
 
\begin{proof} Without loss of generality we may assume that $M$ is a free spectrum.
We rely on Corollary \ref{sorpresa} below and show that the levelwise projections  in (\ref{reducciones}) of Definition \ref{indelineal}
$$
 M\stackrel{\sim}{\longrightarrow} M_{\rm ind}\stackrel{\sim}{\longleftarrow}M_{\rm lin}
$$ 
 induce  stable homology isomorphism at every Maurer-Cartan element of $L$ and thus, they are stable equivalences. We focus on  $q\colon M\to M_{\rm ind}$, the other case being completely analogous. As usual, write $M^n= L\oplus (\libc(T^n),d)$ for any $n\ge 0$, and 
consider the map of sequences as in (\ref{secuencia}) 
$$
\xymatrix{
(\libc(T^0),d)\ar[r]^(.40){\mu_0}\ar[d]& (s^{-1}\libc(T^1),s^{-1}d)\ar[r]\ar[d]&\dots\ar[r]&(s^{-n}\libc(T^n),s^{-n}d)
\ar[r]^(.69){s^{-n}\mu_n}\ar[d]&
\dots&\\
(T^0,d_1)\ar[r]^(.40){\bar\mu_0}& (s^{-1} T^1,s^{-1}d_1)\ar[r]&\dots\ar[r]&(s^{-n}T^n,s^{-n}d_1)\ar[r]^(.69){s^{-n}\bar\mu_n}&
\dots&\\}
$$
where each $\mu^\prime_n$ is the map induced by $\mu_n$ in the indecomposables and the vertical arrows are the projections induced by $q$. As the codomain of  $\mu_n\colon (\libc(T^n),d)\to (s^{-1}\libc(T^{n+1}),s^{-1}d)$ is abelian, it sends $\libc^{\ge 2}(T^n)$ to $0$ and  thus, this diagram induces an isomorphism
$$
\varinjlim_n(s^{-n}\libc(T^n),s^{-n}d)\stackrel{\cong}{\longrightarrow}
\varinjlim_n(s^{-n}T^n,s^{-n}d_1).
$$
In particular, the induced map in stable homology $H^{\rm st}(q)$ is an isomorphism as claimed.

Finally, note  that this argument holds independently of any perturbation of the differential of $M$ by any Maurer-Cartan element of $L$. 
\end{proof}

\begin{example}\label{suspensionl}{\em (The suspension $L$-spectrum)}
Let $M\in\specl$. Recall from Proposition \ref{generalcdgl} that the {\em suspension $L$-spectrum} $\susel M$ is defined as the prolongation of the unstable retractive suspension, i.e.,  $(\susel M)^n=\susl M^n$ with structure maps induced from those  of $M$. By Theorem \ref{freespec} we may choose a model of $M$ which, by the preceding result can be taken to be linear, say $(L\amalg\libc(W),d_1)$. Applying  Theorem \ref{suspensionex} levelwise, we see that $\susel M$ is  stably equivalent to the free $L$-spectrum $(L\amalg\libc(sW),d_1)$, whose $n$th level is $(L\amalg\libc(sW^n),d_1)$ with the differential given by  $d_1sw=-sd_1w$ for $w\in W$, and with the appropriate structure maps.
\end{example}

\subsubsection{The Quillen equivalence}

The content of the previous sections, along with some  results from \S\ref{apendice}, culminates in the following result. Fix a simplicial set  $B$ and denote $L=\lasu_B$. 

As the model and realization functors (\ref{pair}),
$$
\xymatrix{ \catss& \catcdgl, \ar@<0.75ex>[l]^(.50){\langle\,\cdot\,\rangle}
\ar@<0.75ex>[l];[]^(.50){\lasu}\\}
$$
 form a right-transferred pair, they induce, by Proposition \ref{transre}, a right-transferred Quillen adjunction
\begin{equation}\label{funtorbajo}
\xymatrix{ \retb & \cdgll. \ar@<0.75ex>[l]^(.45){{\langle\,\cdot\,\rangle}_{\sslash B}}
\ar@<0.75ex>[l];[]^(.50){\lasu_{\sslash B}}\\}
\end{equation}
\begin{theorem}\label{primero}This pair 
induces a right-transferred Quillen adjunction
\begin{equation}\label{casimejor}
\xymatrix{ \specb& \specl \ar@<0.75ex>[l]^(.50){\bm{\langle\,\cdot\,\rangle}}
\ar@<0.75ex>[l];[]^(.50){\lasub}\\}
\end{equation}
between the corresponding retractive spectra with respect to their stable model structures. Moreover, $\bm{\langle \,\cdot\,\rangle}$ is the prolongation of ${\langle\,\cdot\,\rangle}_{\sslash B}$ and $\,\lasub\,\susb^{\infty-k}= {\susl}^{\infty-k}\,\lasub$ for all $k$.
\end{theorem}

To avoid excessive notation, and following common practice in stable homotopy theory, we use the same notation, but in boldface, for the functors induced between stable categories as for their unstable counterparts.

\begin{proof}

Recall that $\susb$ and  $\lupb$ induce the genuine suspension and loop functors, respectively, in  $\Ho\retb$. Likewise, $\susl$ and $\lupl$ induce the suspension and loop functors  in  $\Ho\cdgll$. Since every object of $\retb$ is cofibrant, every object of  $\cdgll$ is fibrant, and (\ref{funtorbajo})  is a Quillen pair, we obtain conjugate natural tranformations
$$
\lasu_{\sslash B}\,\susb\stackrel{\sim}{\longrightarrow} \susl \,\lasu_{\sslash B}\quad\text{and}\quad {\langle\,\cdot\,\rangle}_{\sslash B}\,\lupl\stackrel{\sim}{\longrightarrow}\lupb\, {\langle\,\cdot\,\rangle}_{\sslash B}
$$
which take values in weak equivalences. We may therefore apply Theorem \ref{prolon}(i) to deduce the existence of the required Quillen adjunction.

Finally, the general framework developed in the previous sections allows us  to invoke Proposition \ref{trans} and Corollary \ref{transcor},  ensuring  that this is indeed a right-transferred pair.
\end{proof}

As a first consequence we obtain:

\begin{proposition}\label{espnoesp1} For any $k\ge 0$, the following diagrams of functors commute up to isomorphisms:
$$
\xymatrix{\retb\ar[d]_{\Sigma^{\infty-k}_B}\ar[r]^{\lasu}& \cdgll \ar[d]^{\Sigma^{\infty-k}_L} & \quad & \retb&\cdgll \ar[l]_{\langle\,\cdot\,\rangle\sslash B} 
\\ \spec_B\ar[r]_{\lasub}&\specl, & \quad & \spec_B\ar[u]^(.40){\eva_k}&\specl\ar[l]^{\bm{\langle\,\cdot\,\rangle}}\ar[u]_(.40){\eva_k}. } 
$$
\end{proposition}
\begin{proof} The right square consists of the right adjoint functors corresponding to the functors in the left square. Since the right diagram  commutes up to isomorphism, the left one must also commute up to isomorphism by adjunction.
\end{proof}

On the other hand, the homological characterization of stable equivalences in $\specl$  is a consequence of the following:

\begin{theorem}\label{homotohomolo}
For any $M\in\specl$ and any $k\in\bz$,
$$
\pie_k\bm\langle M\bm\rangle\cong \{H^{\rm st}_{k-1}(M,a)\}
$$
as $a$ ranges over the non-trivial Maurer-Cartan elements of $L$.
In particular, this set is an invariant of $\Ho\specl$.
\end{theorem}

\begin{proof}
  The right-hand square in Proposition \ref{espnoesp1} implies that, for every
$M\in\specl$ and every $n\ge0$,
$$
\bm\langle M\bm\rangle_n=\langle M^n\rangle_{\sslash B}
$$
which, by definition, is obtained as the  pullback
\begin{equation}\label{pullback}
\xymatrix{
\bm\langle M\bm\rangle_n\ar[d]\ar[r]&B\ar[d]^(.45)\eta\\
\langle M^n\rangle\ar[r]_{\langle p^n\rangle}&\langle L\rangle.
}
\end{equation}
Note that any $L$-spectrum is fibrant in the projective model structure, and therefore its realization spectrum is fibrant as well. Therefore, the $n$th level $\bm\langle M\bm\rangle^b_n$ of the fiber spectrum of $\bm\langle M\bm\rangle$ at $b\in B$ is given by the fiber of $\bm\langle M\bm\rangle_n\to B$ at $b$
$$
\xymatrix{
\bm\langle M\bm\rangle^b_n\ar[d]\ar[r]&{*}\ar[d]^(.45)b\\
\bm\langle M\bm\rangle_n\ar[r]&B,
}
$$
which is then the fiber of $\langle M^n\rangle\to\langle L\rangle$ at $\eta(b)$,
$$
\xymatrix{
\bm\langle M\bm\rangle^b_n\ar[d]\ar[r]&{*}\ar[d]^{\eta(b)}\\
\langle M^n\rangle\ar[r]_{\langle p^n\rangle}&\langle L\rangle.
}
$$
Recall that the map $\eta\colon B\to\langle L \rangle=\langle \lasu_B \rangle$ is the unit of the adjunction $\lasu\dashv\langle\,\cdot\,\rangle$ which induces a bijection between the $0$-simplices of $B$ and the   $0$-simplices of $\langle L \rangle$ corresponding to non-trivial MC elements of $L=\lasu_B$. Thus,  $\eta(b)=b$. 

Therefore, taking into account that $\langle L\rangle \cong \langle L,d_b\rangle$, see (\ref{nodif}), we can reformulate the above diagram as 
$$
\xymatrix{
\bm\langle M\bm\rangle^b_n\ar[d]\ar[r]&\langle 0\rangle\ar[d]_{\langle 0\rangle}\\
\langle M^n,d_b\rangle\ar[r]_{\langle p^n\rangle}&\langle L,d_b\rangle.
}
$$
Observe that perturbing the differential is essential: otherwise the map $\eta(b)\colon \langle 0\rangle\to \langle L\rangle$ does not arise as the realization of a cdgl morphism as $0$ is sent to $b$. 

Thus, by the splitting in (\ref{splittingmc}) and since the realization functor preserves limits, we conclude that
\begin{equation}\label{realizafibra}
\bm\langle M\bm\rangle^b_n\simeq \langle K^n,d_b\rangle.
\end{equation}
Therefore, by (\ref{muchas}) and (\ref{homoto}), and for any $m\ge 1$,
$$
\pi_m\bm\langle M\bm\rangle^b_n\cong H_{m-1}(K^n,d_b)^{(0)}.
$$
  Hence, see Definition \ref{stableho}(ii), the $k$-th fiberwise stable homotopy group of the  realization spectrum of $M$ is  
$$
\pie_k\bm\langle M\bm\rangle = \{\varinjlim_n\pi_{k+n}\bm\langle M\bm\rangle^b_n\}_{b\in B}=\{\varinjlim_n H_{k+n-1}(K^n,d_b)^{(0)}\}_{b\in B}.
$$
The result follows from Definition \ref{sorpresa2} and Remark \ref{remarkestable}.
\end{proof}
 As a result, in light of Theorem \ref{primero}, we obtain:
 
 \begin{corollary}\label{sorpresa}
  A map $f\colon M\to N$ in $\specl$ is a stable equivalence if and only if it induces an isomorphism in stable homology,
   $$
   H^{\rm st}(f)\colon H^{\rm st}(M,a)\stackrel{\cong}{\longrightarrow}H^{\rm st}(N,a).
    $$
  at any non trivial  $a\in\mc(L)$.
    \hfill$\square$
  \end{corollary}

The main result of this section reads: 

 \begin{theorem}\label{primeroq}
The Quillen pair {\em (\ref{casimejor})} induces a  right-transferred Quillen equivalence
\begin{equation}\label{mejor}
\xymatrix{ \spec_B^\bq& \specl. \ar@<0.75ex>[l]^(.50){\bm{\langle\,\cdot\,\rangle}}
\ar@<0.75ex>[l];[]^(.50){\lasub}\\}
\end{equation}
\end{theorem}
\begin{proof}
By Proposition \ref{transferbous}, to show that the pair (\ref{casimejor}) induces a right-transferred Quillen pair as stated, it suffices to verify that for any $M\in\specl$, its realization spectrum $\bm\langle M\bm\rangle$ is $\bq$-local. By Definition \ref{homoq}, this is equivalent to saying that, for each $k\in\bz$, the fiberwise stable homotopy groups $\pi_k^{\rm st}\bm\langle M\bm\rangle$ form a collection of rational vector spaces.  This follows immediately from Theorem \ref{homotohomolo}.

Next, we show that the induced right-transferred pair (\ref{mejor}) is a Quillen equivalence. To this end, since the realization $\bm{\langle\,\cdot\,\rangle}\colon \specl\to\spec_B^\bq$ reflects weak equivalences, it suffices to prove, see for instance \cite[Cor.~1.3.16(c)]{ho0}, that the unit of the adjunction (\ref{mejor}) is a weak equivalence for any rational spectrum. In other words, for any $X\in \specb$, the unit map  
$$
\eta\colon X\longrightarrow \bm\langle \lasub_X\bm\rangle$$
is a rational equivalence. 

To this end, recall that, by definition, the $n$th level stage of $\bm\langle\lasu_X\bm\rangle$ fits into the pullback
$$
\xymatrix{
\bm\langle \lasub_X\bm\rangle_n\ar[d]\ar[r]&B\ar[d]\\
\langle (\lasub_X)^n\rangle\ar[r]_(.55){\langle p^n\rangle}&\langle L\rangle.
}
$$
On the one side, see (\ref{bousfikan}), the right vertical arrow is weakly equivalent to $B\to B^\acento_\bq\amalg \{*\}$. On the other side, by Remark \ref{uyuy}, there is a weak equivalence $(\lasub_X)^n\stackrel{\sim}{\to} \lasu_{X_n}$ so that the bottom map is weakly equivalent to the map 
$$
(q_n)^\acento_\bq\amalg\{*\}\colon (X_n)^\acento_\bq\amalg\{*\}\longrightarrow B^\acento_\bq\amalg \{*\}
$$
induced by the retraction $q_n\colon X_n\to B$.

Assume now that each $X_n$ is fibrant in the projective structure, i.e.,  $q_n$ is a fibration. Then, the existence of a section  of $q_n$ ensures that, for any $b\in B$, the action of $\pi_1(B,b)$ on $X^b_n$ is trivial so we may invoke   \cite[Lemma 5.1]{bouskan} to conclude that $(q_n)^\acento_\bq$ is a fibration whose fiber is homotopy equivalent to $(X_n^b)^\acento_\bq$. As a result, and for each $b\in B$, we have a homotopy commutative diagram 
$$
\xymatrix{
\bm\langle \lasub_X\bm\rangle_n^b\ar[r]\ar[d]_\simeq^{\varphi_n} &
\bm\langle \lasub_X\bm\rangle_n\ar[d]\ar[r]& B\ar[d]\\
(X^b_n)^\acento_\bq\ar[r]&(X_n)^\acento_\bq\ar[r]&B^\acento_\bq
}
$$
in which both horizontal lines are fibration sequences.

On the other hand, observe that for each $b\in B$, the retraction $q_n\colon X_n\to B$ together with the $n$th level map $\eta_n\colon X_n\to \bm\langle\lasub_X\bm\rangle_n$ of $\eta$ determines a commutative diagram
$$
\xymatrix{
X^b_n\ar[r]\ar[d]_{\eta_n^b}&
X_n\ar[d]^{\eta_n}\ar[r]^{q_n}& B\ar@{=}[d]\\
\bm\langle \lasub_X\bm\rangle_n^b\ar[r]&\bm\langle \lasub_X\bm\rangle_n\ar[r]&B,
}
$$
where $\eta_n^b$ denotes the $n$th level of the map $\eta^b\colon X^b\to\bm\langle \lasub_X\bm\rangle^b$ induced by $\eta$ between the respective fiber spectra. 

The composition of both diagrams above yield a  homotopy commutative diagram of fibrations
$$
\xymatrix{
X^b_n\ar[r]\ar[d]_{\varphi_n\eta_n^b}&
X_n\ar[d]_{(\,\cdot\,)^\acento_\bq}\ar[r]^{q_n}& B\ar[d]^{(\,\cdot\,)^\acento_\bq}\\
(X^b_n)^\acento_\bq\ar[r]& (X_n)^\acento_\bq\ar[r]_(.53){(q_n)^\acento_\bq}&B^\acento_\bq.
}
$$
The existence of a section of $q_n$ also implies that the inclusion of the homotopy fiber of $X_n^b\to X_n$, and hence that of $(X^b)^\acento_\bq\to(X_n)^\acento_\bq$, is null-homotopic so there is no obstruction to constructing homotopies
 between all maps $X^b_n\to  (X_n^b)^\acento_\bq$ fitting in the diagram above. In particular, $\varphi_n\eta_n^b$ must be homotopic to the $\bq$-completion.

Now,  without loss of generality, we may assume that $X$ is fibrant in the stable model category,  so that each fiber spectrum $X^b=b^* X$ is also fibrant for every $b\in B$. It follows that, for each $n$, both $X_n$ and $X^b_n$ have the homotopy type  of  loop spaces and therefore all their path components are nilpotent spaces. In particular their $\bq$-completions agree, up to homotopy, with classical rationalization. Hence, up to homotopy, the composite $\varphi_n\eta^b_n\colon X^b_n\to(X_n^b)^\acento_\bq$ is the rationalization map. Since $\varphi_n$ is a homotopy equivalence,  it follows that  $\eta^b_n$ is a rational equivalence for any $b\in B$ and any $n\ge 0$.  As a consequence, $\eta\colon X\to\bm\langle\lasub_X\bm\rangle$ is a rational equivalence in $\spec_B$.
\end{proof}

The previous result motivates the following notion.

\begin{definition}\label{modelspec}
Let $X\in\spec_B$. We define the {\em model}, {\em linear model} and {\em indecomposable model} of $X$ to be the model, linear model and indecomposable model of $\lasub_X$, respectively.
In particular, the  equivalence induced by Theorem \ref{primeroq}
$$
\Ho\spec_B^\bq\cong \Ho\specl
$$
assigns, up to isomorphisms, to each $B$-spectrum $X$  any chosen  model of $X$, including linear or indecomposable ones (see Proposition \ref{estahomolin}).  \hfill$\square$
\end{definition}

 Recall from (\ref{pair2}) in \S\ref{homotocdgl} that the model and realization functor (\ref{pair}) can be slightly altered to yield a Quillen pair
\begin{equation}\label{otropair}
\xymatrix{ \catss^*& \catcdgl \ar@<0.75ex>[l]^(.45){\langle\,\cdot\,\rangle}
\ar@<0.75ex>[l];[]^(.53){\lasu^*}\\}
\end{equation}
where now $\lasu^*_B=\lasu_B/(b)$, with $(b)$ denoting the ideal of $\lasu_B$ generated by the based $0$-simplex $b\in B$. Then, the procedure described in this section applies mutatis mutandis to this setting to establish the analogue of Theorems \ref{primero} and \ref{primeroq} starting from (\ref{otropair}):

\begin{theorem}\label{primeroqq}  There is a Quillen pair
$$
\xymatrix{ \specb& \specl, \ar@<0.75ex>[l]^(.50){\bm{\langle\,\cdot\,\rangle}}
\ar@<0.75ex>[l];[]^(.50){\lasub^*}\\}
$$
which induces a Quillen equivalence
\begin{equation*}
\begin{aligned}
\xymatrix{ \spec_B^\bq& \specl. \ar@<0.75ex>[l]^(.50){\bm{\langle\,\cdot\,\rangle}}
\ar@<0.75ex>[l];[]^(.50){\lasub^*}}
\end{aligned}
\end{equation*}
\hfill$\square$
\end{theorem}

We now address the naturality of the preceding results.
Let $f\colon B\to B'$ be a map of simplicial sets, write $\varphi\colon L\to  L'$ for the induced map $\lasu_f\colon\lasu_B\to\lasu_{B'}$, and consider the functors
 $$
 \widetilde f_!\colon \specb\to\spec_{B'}\quad \text{and}\quad \widetilde\varphi_!\colon \specl\to\spec_{L'}
 $$
provided by Propositions \ref{changespec} and \ref{changespecl} respectively. Then:

\begin{proposition}\label{natural1} The following commutes: 
$$
\xymatrix{
\Ho\specb\ar[d]_{\widetilde f_!}\ar[r]^{\lasub}&\Ho\specl\ar[d]^{\widetilde\varphi_!}\\
\Ho\spec_{B'}\ar[r]_(.48){\lasub}&\Ho\spec_{L'}.
}
$$
\end{proposition}

\begin{proof} The proof follows from the definitions: in the unstable setting, since the model functor preserves colimits, one  checks that the diagram
$$
\xymatrix{
\retb\ar[d]_{f_!}\ar[r]^{\lasu}&\cdgll\ar[d]^{\varphi_!}\\
\catss_{\sslash B'}\ar[r]_(.48){\lasu}&\cdgl_{\sslash L'}
}
$$
commutes. Moreover,  $\widetilde f_!$ and $\widetilde\varphi_!$ are the prolongation of $f_1$ and $\varphi_!$ respectively (see Propositions \ref{changespec} and \ref{changespecl}), while in the homotopy categories, $\lasub$ is likewise the prolongation of $\lasu$ (see Remark \ref{uyuy}). The claim then follows by direct inspection.
\end{proof}
We conclude with some essential observations.

\begin{rem}\label{casoconexo} {\em (Refinements over a connected base)$\,$}\label{conexidad}

\smallskip

Assume now that $B$ is a connected simplicial set which we may take to be reduced. Let $b$ denote its unique $0$-simplex, and as before set  $L=\lasu^*_B$. 

\smallskip

(i) In this case, $L$ is a connected cdgl and thus,  by Proposition \ref{quillenchange}, we have a Quillen equivalences,
$$
\xymatrix{ \spec_{\libc} & \specl \ar@<0.75ex>[l]^(.50){\widetilde{\varphi^*}}
\ar@<0.75ex>[l];[]^(.50){\widetilde{\varphi}_!}\\}
$$
for any    Lie model $\varphi\colon\libc\stackrel{\simeq}{\to} L$ of $B$. Therefore, for all practical purposes, we may choose any Lie model of $B$ as the base cdgl to develop our theory.

\smallskip

(ii) Moreover, the pointed model functor $\lasu^*$ sends the $0$-simplex $b$ to the trivial Maurer-Cartan element of $L$. In this particular context note that Theorem \ref{homotohomolo} now states that for any $M\in\specl$
$$
\pi^{\rm st}\bm\langle M\bm\rangle\cong s\, H^{\rm st} (M).
$$
Equivalently, for any $X\in\spec_B$,
$$
\pi^{\rm st}(X)\otimes\bq\cong s\, H^{\rm st} (\lasub^*_X).
$$
Accordingly,
Corollary \ref{sorpresa} now asserts that a map $f\colon M\to N$ in $\specl$ is a stable equivalence if and only if 
$$
H^{\rm st}(f)\colon H^{\rm st}(M)\stackrel{\cong}{\longrightarrow}  H^{\rm st}(N)
$$ is an isomorphism. 

\smallskip

(iii) Another important feature arises when $L$ is connected: for any 
$$
M=\{M^n\}_{n\ge 0}\cong\{L\oplus K^n\}_{n\ge 0}\in\specl 
$$
consider its connected cover
$$
M^{(0)}=\{(M^n)^{(0)}\}_{n\ge 0}\cong \{L\oplus (K^n)^{(0)}\}_{n\ge 0}
$$
which is a well defined $\cdgll$ spectrum. Then, the inclusion $M^{(0)}\hookrightarrow M$ is a stable equivalence as it induces an isomorphism in stable homology.  In other words, every spectrum in $\specl$ is stably equivalent to a levelwise connected spectrum.

\smallskip

(iv) Furthermore, let $M=(L\amalg \libc(W),d)$ be a free $L$-spectrum and consider its indecomposable reduction $M_{\rm ind}=(L\amalg W,d_1)$, see (ii) of Definition \ref{indelineal}. By taking the connected cover at each level  we find a levelwise connected $L$-spectrum  denoted by $(L\amalg W^{(0)},d_1)$. The inclusion
$$
(L\amalg W^{(0)},d_1)\stackrel{\sim}{\longrightarrow} M_{\rm ind}
$$
is a stable equivalence as it  induces an isomorphism in stable homology. On the other hand, By Proposition \ref{estahomolin} the projection 
$$
M\stackrel{\sim}{\longrightarrow}M_{\rm ind}=(L\amalg W,d_1)
$$
is also a stable equivalence.  Thus, $M$
is stably equivalent to $(L\,\amalg\, W^{(0)},d_1)
 $,
 which in turn is stably equivalent to the free spectrum 
 $
(L\amalg \libc(W^{(0)}),d_1)
 $.
In particular the equivalence
$$
\Ho\spec_B^\bq\cong \Ho\specl
$$
 identifies any $B$-spectrum $X$ with the levelwise connected spectra $(L\,\amalg\, W^{(0)},d_1)$ and  $(L\amalg\,\libc(W^{(0)}),d_1)$ where $(L\amalg \libc(W),d)$ is a model of   $X$. 
\end{rem}

%

\section{From spectra of retractive Lie algebras to spectra of Lie modules}
\subsection{Complete modules over a dga}\label{assmod}

We begin by  collecting the essential properties of the category of complete differential graded modules over a differential graded algebra, dga hereafter. Throughout, any  dga $A$ considered will be associative, unital and augmented, with $ \bar A$ denoting its augmentation ideal. Unless otherwise stated, a (left) $A$-module always refers to a differential graded module over $A$.

\begin{definition} \label{filtermodalg}
A {\em filtration} of a dga $A$ is a decreasing sequence of differential ideals $\{G^n\}_{n\ge 0}$ such that
$$G^0=A,\quad G^1=\bar A\quad\text{and}\quad G^pG^q\subset G^{p+q}\quad\text{for all $p,q\ge 0$}.
$$
Given a dga $A$ filtered by $\{G^n\}_{n\ge 0}$, a {\em filtration} of a differential graded $A$-module $R$ is a decreasing sequence of submodules $\{S^n\}_{n\ge 0}$ where
$$S^0=R\quad\text{and}\quad G^pS^q\subset S^{p+q}\quad\text{for all $p,q\ge 0$}.
$$
 A {\em complete differential graded module} (cdgm henceforth) is a filtered $A$-module $R$ for which the natural map 
$$
R\stackrel{\cong}{\longrightarrow}\varprojlim_n R/S^n
$$
is an isomorphism. Morphisms of cdgm's are maps of differential graded modules that preserve the given filtrations. We denote the corresponding category by $\cdgma$. 
\end{definition}

\begin{rem}\label{adic}
If $S=\{S^n\}_{n\ge 0}$ and $S'=\{S'^n\}_{n\ge 0}$ are filtrations of $R$ with $S^n\subset {S'}^n$ for all $n$, and  $R$ is complete respect to $S$, then it is also complete respect to $S'$.  In particular, and for any filtration of $A$, any cdgm $R$ is complete with respect to the {\em adic} filtration  defined as $\{\bar A^nR\}_{n\ge0}$, with the convention $\bar A^0R=R$.
\end{rem}

\begin{definition}\label{completion} Given a filtered $A$-module $R$, its {\em completion} is defined as
$$
\widehat R=\varprojlim_n R/S^n.
$$
This is always a complete $A$-module with respect to the filtration $\{{\widehat S}^n\}_{n\ge 0}$, 
\begin{equation}\label{filtrationcomp}
\widehat S^n=\ker(\widehat R\to R/S^n),
\end{equation} 
since, for each $n\ge 0$, there is a natural  isomorphism
\begin{equation}\label{filtracioncociente}
\widehat R/\widehat S^n\cong R/S^n. 
\end{equation}
Observe also that a map $R\to S$ between filtered modules preserves the filtration if and only if the induced map $\widehat R\to\widehat S$ does.
\end{definition}

\begin{proposition}\label{cocomp}
 The category $\cdgma$ is bicomplete.
\end{proposition}
\begin{proof} Binary equalizers and coequalizers in $\cdgma$, as well as small products,  are the usual ones, which are complete with respect to the obvious filtrations. Finally, the coproduct of  a family $\{R_i\}_{i\in I}$ of cdgm's, each of which is filtered  by $\{S_i^n\}_{n\ge 0}$, is given by
$$
\amalg_{i\in I}R_i=\varprojlim_n  (\oplus_{i\in I} R_i)/(\oplus_{i\in I} S_i^n).
$$
\end{proof}

The following is a straightforward adaptation of its non-complete counterpart, see for instance \cite[\S6]{fehaltho}:
\begin{definition}\label{semifreem} Given a dga $A$ filtered by $\{G^n\}_{n\ge 0}$, an $A$-module is {\em semifree} if it is of the form $(A\otimesc W,d)$ where:
$$
A\otimesc W=\varprojlim_n A/G^n\otimes W,
$$
and $W=\cup_{k\ge 0} W(k)$ is the union of an increasing sequence of graded vector spaces such that  
\begin{equation}\label{diferencial}
dW(0)=0\quad\text{and}\quad dW(k)\subset A{\otimesc} \bigl(\oplus_{i<k}W(i)\bigr), \quad k\ge1 .
\end{equation}
A {\em (complete) semifree resolution} of a  cdgm $R$ is a quasi-isomorphism with semifree domain of the form
$$
(A\otimesc W,d)\stackrel{\simeq}{\longrightarrow} R.
$$
Observe that any  semifree $A$-module is complete with respect to the filtration $\{\mathcal{G}^n\}_{n\ge 0}$ where
$$
\mathcal{G}^n=\ker\bigl(A\otimesc W\to (A\otimes W)/(G^n\otimes W)\cong A/G^n\otimes W\bigr).
$$
Indeed, by (\ref{diferencial}), these are differential modules and 
$$
(A\otimesc W)/\mathcal{G}^n\cong (A\otimes W)/(G^n\otimes W)\cong A/G^n\otimes W.
$$
\end{definition}

\begin{rem}\label{cgm} Let $A$ be a  graded algebra with trivial differential filtered by $\{G^n\}_{n\ge 0}$, and consider the category $\cgm_A$ of  complete $A$-modules also endowed with trivial differential. Given $R\in\cgm_A$, any linear map $W\to R$ from a graded vector space produces a $\cgm_A$ morphism $A\otimesc W\to R$. This shows that the functor $A\otimesc -$ from the category of graded vector spaces  to $\cgm_A$ is left adjoint to the forgetful functor. In particular,
$$
A\otimesc (\oplus_{i\in I}W_i)\cong \amalg_{i\in I} A\otimesc W_i
$$
for any family $\{W_i\}_{i\in I}$ of graded vector spaces.
\end{rem}

Adapting the proof of \cite[Prop.~6.6(i)]{fehaltho} to our context provides:

\begin{proposition}\label{semifreer}
Every complete $A$-module has a semifree resolution.
\end{proposition}

\begin{proof}
Let $R$ be a complete $A$-module filtered by $\{S^n\}_{n\ge0}$ and let $W(0)$ be a copy of the vector space of cocycles of $R$. Set $dW(0)=0$ and observe that the map $\varphi_0\colon A\otimesc W(0)\to R$, induced by the inclusion $W(0)\hookrightarrow R$, is a cdgm morphism since $\varphi_0(\mathcal{G}^n)\subset G^nR\subset S^n$, and it respects the differentials. 
Assume 
$$
\varphi_{k-1}\colon (A\otimesc W(k-1),d)\longrightarrow R
$$
 has been constructed  and let $U$ be the suspension of a copy of $\ker H(\varphi_{k-1})$. Define $d$ on $U$ in the natural way, set $W(k)=W(k-1)\oplus U$, and extend $\varphi_{k-1}$ to a cdgm morphism $\varphi_k\colon (A\otimesc W(k),d)\to R$ accordingly, making use of Remark \ref{cgm}. Let $W=\cup_{k\ge 0} W(k)$. Then the induced cdgm morphism 
  $$
\varphi=\varinjlim_k\varphi_k\colon (A\otimesc W,d)\stackrel{\simeq}{\longrightarrow} R
$$
is a semifree resolution of $R$.

\end{proof}

\begin{rem}\label{surjective}
As in the non-complete context, note that semifree resolutions can be chosen to be surjective if desired. Indeed,  any semifree resolution  $ (A\otimesc W,d)\stackrel{\simeq}{\to} R$ can be extended to a surjective one
$$
(A\otimesc (W\oplus R\oplus dR),d)\cong (A\otimesc W,d)\amalg A\otimesc (R\oplus dR)\stackrel{\simeq}{\longrightarrow} R
$$
in the obvious way.
\end{rem}

 We
equip $\cdgma$ with the {\em projective model structure} in which weak equivalences and fibrations are quasi-isomorphisms and surjections. For the following result, we refer to \cite[\S3.2]{barmayriehl}, which  presents a  version of the well-known uncompleted counterpart, and whose particular proof adapts to the complete setting.

\begin{proposition}\label{cdgmalocally} For any filtered dga $A$, the category $\cdgma$ is proper and combinatorial.\hfill$\square$
\end{proposition} 

This structure is cofibrantly generated by the sets 
$$
\{A\otimesc\langle w_n\rangle\to \{A\otimesc\langle u_{n+1} w_n\rangle\}_{n\ge 0}\quad \text{and}\quad \{0\to A\otimesc\langle u_{n+1} w_n\rangle\}_{n\ge 0}
$$
 of cofibrations and trivial cofibrations, respectively, where $w_n$ is a cycle of degree $n$ and $du_{n+1}=w_n$. Any semifree resolution of a  complete $A$-module is, in particular, a cofibrant replacement. Moreover, by \ref{semifreer} and Remark \ref{surjective},  the class of semifree $A$-modules forms a set of compact objects that generates $\cdgma$ by colimits.

\subsection{Complete Lie modules and their connection with retractive cdgl's}\label{liemod}

\subsubsection{Complete Lie modules}

Recall that, given a differential graded Lie algebra $L$, a (left) $L$-module, or simply {\em module} when no confusion arises,  is a differential graded vector space $R$ equipped with a differential graded linear map 
$
L\otimes R\to  R$, $ x\otimes a\mapsto xa$,
which encodes the Lie bracket on $L$, that is,  
$[x,y]a=x(ya)-(-1)^{|x||y|}y(xa)$ for any $x,y\in L$ and $a\in R$. The category  of $L$-modules is equivalent to the category of (left) differential graded modules over the universal enveloping algebra $UL$. Specifically, given an $L$-module $R$, the corresponding dgl representation $\rho\colon L\to \End(R)$  extends uniquely to a differential graded algebra morphism $\overline\rho\colon UL\to \End(R)$   for which the following commutes,
$$
\xymatrix{L\ar[rd]_(.40)\rho\ar@{^{(}->}[r]&UL \ar[d]^{\overline\rho}\\
& \End(R).
}
$$

We  briefly sketch how the above extends to the completed setting, much of which  is part of the folklore.  

The essential notions in the Lie setting are defined mutatis mutandis from their counterparts in the associative case:

\begin{definition}\label{filtermodlie}
Let $L$ be a dgl filtered by $\{F^n\}_{n\ge 1}$. A {\em filtration} of an $L$-module $R$ is a decreasing sequence of $L$-modules $\{S^n\}_{n\ge 0}$ such that
$$S^0=R\quad\text{and}\quad F^pS^q\subset S^{p+q}\quad\text{for  $p\ge 1$ and $q\ge 0$}.
$$
 A  filtered $L$-module $R$ is {\em complete} if the map 
$$
R\stackrel{\cong}{\longrightarrow}\varprojlim_n R/S^n
$$
is an isomorphism. The category $\cdgml$ of {\em complete differential graded $L$-modules} is defined accordingly.

As in the associative case, and for any filtration of $L$, any complete $L$-module $R$ is always complete with respect to the {\em adic} filtration $\{L^nR\}_{n\ge0}$, with the convention $L^0R=R$.

The {\em completion} of an $L$ module $R$ filtered by $\{S^n\}_{n\ge 0}$ is defined again as $\widehat R=\varprojlim_n R/S^n$ and this is always a complete $L$-module with respect to the induced filtration, see (\ref{filtrationcomp}) and (\ref{filtracioncociente}).

\end{definition}

\begin{example}\label{kmod} (1) Any filtered dgl $L$ is naturally a filtered $L$-module with respect to the filtration
$\{F^n\}_{n\ge0}$, where $F^0=L$ and, for $n\ge1$, $F^n$ denotes the $n$th term of the given filtration of $L$. Furthermore, if $L$ is complete as a dgl, then it is also complete as an $L$-module.

(2) Let $L\stackrel{s}{\to}M\stackrel{p}{\to}L$ be a retractive cdgl. Then,   $M$ becomes  a complete $L$-module via the action $x y=[s(x),y]$ for $x\in L$, $y\in K$, with respect to the filtration $\{F^nM\}_{n\ge 0}$ where $\{F^n\}_{n\ge 1}$ is the filtration on $L$ and $F^0M=M$. As a consequence, $K=\ker p$ is also a complete $L$-module.
\end{example}

A particularly important instance is  the following:

\begin{definition}\label{ulcomp}
Any filtration $\{F^n\}_{n\ge 1}$ of a given dgl $L$ induces in $UL $ the filtration $\{F^nUL \}_{n\ge 0}$ where
$$
F^0UL=UL ,\quad F^nUL =\Span\{ x_1\dots x_p,\,x_i\in F^{n_i},\,\sum n_i=n\},\quad n\ge 1.
$$
Define the {\em complete universal enveloping algebra} as the completion of $UL $ with respect to this filtration,
$$
\compul=\varprojlim_nUL /F^nUL .
$$
 This completed dga satisfies the following universal property: let $A$ be a dga filtered by $\{G^n\}_{n\ge 0}$ and let $\varphi\colon L\to A$ be a filtered dgl morphism, where the Lie bracket on $A$ is given by the commutator. Then,  the unique extension $UL \to A$ of $\varphi$ preserves the corresponding filtrations and thus, it induces  a unique morphism of complete dga's $\compul\to\widehat A$ such that the following commutes:
$$
\xymatrix{\widehat L\ar[rd]_(.45){\widehat\varphi}\ar[r]&\compul\ar[d]\\
& \widehat A.
}
$$
\end{definition}

\begin{example}\label{freecomp}
Let $L=\lib(V)$ be a free Lie algebra. Then, $UL \cong T(V)$ while $F^nUL =T^{\ge n}(V)$ for $n\ge 1$. Hence,
$$
\compul\cong \widehat{T}(V)=\varprojlim_n T(V)/T^{\ge n}(V)={\textstyle \prod}_n\, T^n(V).
$$
\end{example}

\begin{rem}\label{semifreeuu} (i)
Let $L$ be a filtered dgl and let $W$ be a graded vector space. In light  of  (\ref{filtracioncociente}), we have
$$
\compul\otimesc W\cong UL \otimesc W.
$$
 In particular, any semifree complete $\compul$-module is of this form, equipped with an appropriate  differential.
 
 (ii) On the other hand, for any pair  of filtered dgl's $L$ and $L'$, the classical isomorphism $ U(L\oplus L')\cong UL\otimes UL'$ extends to the completions
 $$
 \widehat{U}(L\oplus L')\cong\compul\otimesc\widehat{UL'}.
 $$
\end{rem}

We now sketch a proof of the following:

\begin{theorem}\label{equivalencia} For any filtered dgl $L$  the categories $\cdgm_{\widehat L}$ and $\cdgm_{\compul}\,$ are equivalent. In particular if $L$ is a cdgl, $\cdgm_{L}$ and $\cdgm_{\compul}\,$ are also equivalent.
\end{theorem}

\begin{proof} Let $R\in\cdgm_{\widehat L}$ filtered by $\{S^n\}_{n\ge 0}$ and consider the dga $\Endog(R)=\Endo(R)\oplus \bq$ of endomorphisms of $R$ equipped with a formal augmentation. Note that $\Endog(R)$ is filtered by $\{G^m\}_{m\ge 0}$
where
$$
G^0=\Endog(R), \quad G^n=\{f\in\Endo(R),\, f(S^q)\subset S^{q+n}, \, q\ge0\},\quad n\ge 1.
$$
Observe also that the map $L\to \Endo(R)$ induced by the module structure on $R$ preserves the respective filtrations so that, as in Definition \ref{ulcomp}, we obtain a commutative diagram
$$
\xymatrix{\widehat L\ar[rd]\ar[r]&\compul\ar[d]\\
& \widehat\Endog(R).
}
$$
However, if $R$ is complete one  checks that $\Endog(R)$ is also complete so that $\widehat\Endog(R)\cong\Endog(R)$ and the vertical map equips $R$ with a structure of complete $\compul$-module.

Conversely, the composition  of a filtered action $\compul\to \widehat\Endog(R)=\Endog(R)$ with the canonical map $\widehat L\to\compul$ provides a complete $\widehat L$-module structure on $R$. 
\end{proof}

As a consequence of Proposition \ref{cdgmalocally}, we obtain:

\begin{proposition}\label{cdgmlocally} For any cdgl $L$, the category $\cdgmul$, equipped with the projective module structure, is proper and combinatorial.\hfill$\square$
\end{proposition}

Now assume  $L$ is a connected cdgl, 
consider the full subcategory $\cdgmul^0$ of $\cdgmul$ consisting of connected $\compul$-modules and  equip it  with the {\em bounded projective model structure}, in which weak equivalences are quasi-isomorphisms and fibrations are surjections in positive degrees. This is also a proper combinatorial model structure. Moreover, the functors
\begin{equation}\label{conemod}
\xymatrix{ \cdgmul^0 & \cdgmul, \ar@<0.75ex>[l]^(.43){(\,\cdot\,)^{(0)}} \ar@<0.75ex>@{^(->}[l];[] \\ }
\end{equation}
where $(\,\cdot\,)^{(0)}$ denotes the connected cover,
form a Quillen pair, with the inclusion  reflecting cofibrations and weak equivalences.

\subsubsection{Bridging $\compul$-modules and retractive cdgl's}

The connection between retractive  Lie algebras and complete Lie modules is provided by the following:

\begin{definition}\label{calker}
Let $L$ be a cdgl filtered by $\{F^n\}_{n\ge 1}$. 
We  define a pair of adjoint functors
  \begin{equation}\label{mejillon}
\xymatrix{ \cdgmul& \cdgll \ar@<0.75ex>[l]^(.43){\calker}
\ar@<0.75ex>[l];[]^(.55){\libc_L}\\}
\end{equation}
As follows: given a retractive cdgl $M$, let $ \calker(M)=K$
     where, as usual, $K$ denotes the kernel of  the retraction of $M$. By Example \ref{kmod} and Theorem \ref{equivalencia} this yields a well-defined functor. Because $\cdgmul$ is complete and locally presentable (see Proposition \ref{cdgmlocally}), and  $\calker$ preserves limits, it admits a left adjoint which we now describe explicitly.
\end{definition}

 Let $(UL \otimesc W,d)$ be a semifree $\compul$-module and consider the injective  map
$$
UL \otimes W\stackrel{j}{\hookrightarrow}L\amalg \lib(W)\quad
$$
defined by
$$
j(w)=w,\quad j\bigl((x_1\dots x_n)\otimes w\bigr)=\bigl[x_1,[x_2,[\dots,[x_{n},w]\bigr]\dots\bigr],
$$
for $x_i\in L$, $i=1,\dots,n$ and $w\in W$. Observe that in the cdgl coproduct $L\amalg\libc(W)$,  the  filtration $\{G^n\}_{n\ge 0}$ is the one induced by the filtration $\{F^n\}_{n\ge 0}$ in $L$ and the bracket-length filtration  in $\lib(W)$ which, with a slight abuse of notation,  we denote by $\{W^n\}_{n\ge 1}$. Explicitly, see \cite[\S3.1]{bufemutan0},
$$
G^n=\sum_{p_1+q_1+\cdots+p_r+q_r=n} \bigl[F^{p_1},[W^{q_1},[\dots,[F^{p_r},W^{q_r}]\bigr]\dots\bigr].
$$
  Note that, for $n\ge 1$, 
  $$
  j (F^nUL \otimes W)\subset [F^n,W]\subset G^n
  $$
   so that $j$ induces an injection 
   $$
  \widehat j\colon UL \otimesc W\hookrightarrow L\amalg\libc(W)
  $$
   between the respective completions,
   where the coproduct is now taken in the category of complete  Lie algebras. In fact, by \cite[\S3.1]{bufemutan0}
  $$
  \varprojlim_n \bigl(L\amalg\lib(W)\bigr)/G^n\cong L\amalg\libc(W).
  $$
  Finally, endow $L\amalg\libc(W)$ with a differential by declaring $L$ to be a sub-cdgl  and setting
$
dw=\widehat j(dw)$ for all $w\in W$. Under these conditions
$$
\widehat j\colon (UL \otimesc W,d)\hookrightarrow (L\amalg \libc(W),d)
$$
becomes a morphism of differential graded objects.

Now, let $R\in \cdgm_{\compul}$ and choose a surjective semifree resolution of $R$ which, by (i) of Remark \ref{semifreeuu}, must be  of the form
$$
(UL \otimesc W,d)\stackrel{\simeq}{\longrightarrow} R,
$$
We may assume this resolution is surjective so that, as $\compul$-modules,
\begin{equation}\label{todomodulo}
R\cong (UL \otimesc W,d)/I.
\end{equation}
Moreover, observe that $\widehat j$ induces a map, which we denote by the same symbol to avoid additional notation,
$$
\widehat j\colon R\hookrightarrow (L\amalg\libc(W),d)/J,
$$
where $J$ is the differential ideal of $(L\amalg\libc(W),d)$ generated by $\widehat j(I)$.
\begin{proposition}\label{primeraprop}
$\libc_L(R)\cong(L\amalg\libc( W),d)/J$.
\end{proposition}
\begin{proof} Note that, by construction, $(L\amalg\libc( W),d)/J$ contains $L$ as a subcdgl and the projection $L\amalg\libc( W)\to L$ sends $J$ to 0, making the  quotient a retractive cdgl over $L$. 

On the other hand, given $L\stackrel{s}{\to}M\stackrel{p}{\to}L$ a retractive cdgl, any morphism of complete $UL $-modules $ f\colon (UL \otimesc W,d)/I\to K$   is  induced by a  map  $ \overline f\colon (UL \otimesc W,d)\to K$ in $\cdgmul$ which vanishes in $I$. This, in turn, defines a morphism of retractive cdgl's $(L\amalg\libc( W),d)\to M$ by sending each $w\in W$ to $\overline f(1\otimes w)$. Since this map annihilates $J$, we obtain a morphism $ (L\amalg\libc( W),d)/J\to M$ in $\cdgll$. This construction defines a bijection 
$$
 \Hom_{\cdgll}\bigl( (L\amalg\libc( W),d)/J,M\bigr)\stackrel{\cong}{\longrightarrow}\Hom_{\cdgml}\bigl( (UL \otimesc W,d)/I,\calker(M)\bigr)
$$
which establishes the proposition.
\end{proof}

Furthermore,  as for the splitting (\ref{splitting}), given $R\in\cdgmul$, consider the free cdgl $(\libc(\mathcal{R}),d)$ where $\mathcal{R}$ is an isomorphic copy of $R$ and  with indecomposable differential induced also by $R$. Then we have:
\begin{proposition}\label{segundaprop} As retractive cdgl's,
$$
 \libc_L(R)\cong L\oplus (\libc(\mathcal{R}),d).
 $$
in which the bracket $[L,\mathcal{R}]\subset\mathcal{R}$ is determined by the $\compul$-module structure. Moreover, given a morphism of $\compul$-modules $f\colon R\to S$ the map
$$
\libc_L(f) \colon L\oplus \libc(\mathcal{R})\longrightarrow L\oplus \libc(\mathcal{S})
$$
acts as the identity on $L$ and is induced by $f$ on $\libc(\mathcal{R})$.
\end{proposition}

\begin{proof}
Assume  that $R\cong (UL \otimesc W,d)$ is  semifree. By Proposition \ref{primeraprop}, $\libc_L(R)\cong (L\amalg\libc( W),d)$ which, in view of  (\ref{esencialk}), can be written as
$$
\libc_L(R)=(L\amalg\libc(W),d)\cong L\oplus(\libc(\mathcal{R}),d)
$$
where the differential and Lie bracket in the right-hand side are as stated. 

In the general case write $R\cong (UL \otimesc W,d)/I$ so that, by Proposition \ref{primeraprop}, 
$$
\libc_L(R)=\libc_L(S)/J,\quad\text{with}\quad S=(UL \otimesc W,d).
$$
By the previous case
$$
\libc_L(S)=L\oplus (\libc(\mathcal{S}),d).
$$
To conclude, observe that $J$ is a sub differential vector space of $\mathcal{S}$ such that $\mathcal{S}/J=\mathcal{R}$. 

The second claim follows immediately.
\end{proof}

The following is an important observation.

\begin{rem}\label{ejemplo1} Let $(L\amalg\libc(W),d)\cong L\oplus(\libc(T),d)$ be a free retractive cdgl. Note that $T\cong \widehat j(UL \otimesc W)$ and we may write
$$
(L\amalg\libc(W),d)\cong L\oplus(\libc(UL \otimesc W),d).
$$
In particular, if the differential $d$ is retractive linear  (Definition \ref{libre}) we obtain
$$
(L\amalg\libc(W),d)\cong \libc_L(UL \otimesc W,d).
$$
\end{rem}

\begin{corollary}\label{quism} A morphism of $\compul$-modules $f\colon R\to S$ is a quasi-isomorphism if and only if  $\libc_L(f)$ is a quasi-isomorphism.\hfill$\square$
\end{corollary}

\begin{proposition}\label{loquecuesta} For any $\compul$-module $R$, the inclusion $L\hookrightarrow \libc_L(R)$ induces a bijection on the $\widetilde\mc$ sets.
 \end{proposition}
 
 \begin{proof} We begin with a general fact. Given any retractive cdgl $M\cong L\oplus K$ , 
  every Maurer-Cartan element  in $M$  can be uniquely expressed as 
\begin{equation}\label{nomasmc}
 x+y,\quad x\in\mc(L),\quad y\in \mc(K,d_x).
\end{equation}
Indeed, any element in $\mc(M)$ is necessarily written as $x+y$ with  $x\in \mc(L)$.  Recall that the perturbed differential $d_x$  in $M$ is well defined on $K$  and that
$$
\mc(M,d_x)=\{z-x,\,\,z\in \mc(M)\}.
$$
Taking $z=x+y$, we deduce that $y\in\mc(K,d_x)$. 

On the other hand, by Proposition \ref{segundaprop}, write
 $
 \libc_L(R)=L\oplus (\libc(\mathcal{R}),d)
 $
and note that, for any $x\in L$ and  $a\in\mathcal{R}=\widehat j(R)$, we have $[x,a]\in \mathcal{R}$. Thus, if $x\in \mc(L)$, the perturbed differential $d_x$ in $(\libc(\mathcal{R}),d)$, like $d$,  preserves $\mathcal{R}$. Hence, write $\mathcal{R}=A\oplus B\oplus d_xB$ where $d_xA=0$ and apply \cite[Prop.~8.10]{bufemutan0} to conclude that the projection
$$
(\libc(\mathcal{R}),d_x)\stackrel{\sim}{\longrightarrow} (\libc(A),0)
$$
is a weak equivalence. Consequently,
$$
\widetilde\mc(\libc(\mathcal{R}),d_x)\cong \widetilde\mc(\libc(A),0)=\{0\}.
$$
This, together with (\ref{nomasmc}), completes the proof.
 \end{proof}
\begin{corollary}
The map $\libc_L(f)$ is a cdgl weak equivalence if and only if,  for any $x\in \mc(L)$, the map induced by $f$ at connected covers 
$$
 (\libc(\mathcal{R}),d_x)^{(0)}\stackrel{\simeq}{\longrightarrow} (\libc(\mathcal{S}),d_x)^{(0)}
$$
is a quasi-isomorphism. \hfill$\square$
\end{corollary}

The following shows that, in its current form, the pair (\ref{mejillon}) cannot serve our purposes as a bridge between  $\cdgmul$ and $\cdgll$.

\begin{example}\label{falla} In general, there is no model structure on $\cdgmul$ for which  the adjunction $\libc_L\dashv\mathcal{K}$ forms a Quillen pair with  $\libc_L$ reflecting cofibrations and weak equivalences.

For instance choosing $L=0$, this adjunction reduces to the free and forgetful functors
$$
\xymatrix{ \catvect& \cdgl \ar@<0.75ex>[l]^(.46){\mathcal{U}}
\ar@<0.75ex>[l];[]^(.52){\libc}\\}
$$
where $\catvect$ denotes the category of differential graded vector spaces. Suppose we attempt to transfer the model structure from $\cdgl$ by declaring a map $f$ in $\catvect$ to be a cofibration or a weak equivalence whenever $\libc(f)$ is. Then,  not every object in $\catvect$ would admit a cofibrant replacement.

To see this, note  that any cofibrant  $(W,d)\in\catvect$ must have homology concentrated in non-negative degrees. Otherwise, write $W\cong Z\oplus C\oplus dC$ with $dZ=0$ and let $z\in Z$ have negative degree. Then the square
$$
\xymatrix@R=15pt@C=15pt{
0\ar[d]\ar@{=}[r]&0\ar[d]\\
(\libc(W),d)\ar[r]&(\libc(z),0)
}
$$
has no lift, even though the right vertical map is a  weak equivalence in $\cdgl$.

Now consider  $V=\Span(a,b)$ with trivial differential, $|a|=-1$, $|b|=1$, and assume $f\colon(W,d)\to (V,0)$ is a cofibrant replacement.  Again, write $W\cong Z\oplus C\oplus dC$ with $dZ=0$, $Z=Z_{\ge 0}$. Since both $\libc(f)$ and the inclusion $(\libc(Z),0)\hookrightarrow (\libc(W),d)$ are weak equivalences, the composition $\lib(Z)\to \lib(a,b)$ must be an isomorphism in non-negative degrees. However, the bracket $[a,b]$ does not lie in its image.
\end{example}

From this point onward, and for the remainder of the section,  we assume that $L$ is connected. Accordingly, we restrict the adjoint pair (\ref{mejillon}) on the left to the connected setting,  obtaining
\begin{equation}\label{mejillon0}
\xymatrix{ \cdgmul^0& \cdgll. \ar@<0.75ex>[l]^(.43){\calker^0}
\ar@<0.75ex>[l];[]^(.55){\libc_L}\\}
\end{equation}
 Here, $\calker^0$ assigns to each $M\in\cdgll$ the connected $\compul$-module $\calker(M^{(0)})$.
It is straightforward to check that  this remains an adjoint pair. 

\begin{rem}\label{doblefun}
Note  that the adjunction (\ref{mejillon0}) also arises as the composition of two adjoint pairs
$$
\xymatrix{ \cdgmul^0& \cdgll^0 \ar@<0.75ex>[l]^(.43){\calker}
\ar@<0.75ex>[l];[]^(.55){\libc_L}&\cdgll \ar@<0.75ex>[l]^(.43){(\,\cdot\,)^{(0)}}
\ar@<0.75ex>@{^(->}[l];[]\\}
$$
where $\cdgll^0$ denotes the category of levelwise connected retractive cdgl's and $(\,\cdot\,)^{(0)}$ is the connected cover functor.  Equipping both  $\cdgll^0$ and $\cdgmul^0$  with the (bounded) projective model structure, we obtain Quillen adjunctions for both pairs. This is immediate for the left-hand pair. For the right one note that the functors $(\,\cdot\,)^{(0)}$ and $\calker$  preserve fibrations and weak equivalences.
\end{rem}

\begin{theorem}\label{modelmod} 
For any connected cdgl $L$ the adjunction $\libc_L\dashv\calker^0$ in (\ref{mejillon0}) is a Quillen pair in which $\libc_L$ reflects cofibrations, and weak equivalences. 
\end{theorem}
\begin{proof}   Let $
\varphi_R\colon (UL\otimesc W,d)\stackrel{\simeq}{\to} R$ be a connected semifree resolution of a given connected $\compul$-module $R$. By Corollary \ref{quism}, $\libc_L(f)\colon (L\amalg\libc(W),d)\stackrel{\simeq}{\to} \libc_L(R)$ is a quasi-isomorphism, hence a weak equivalence in $\cdgll$. On the other hand, since $W$ is non-negatively graded, $L\hookrightarrow (L\amalg \libc(W),d)$ is a cofibration in $\cdgll^0$. This shows that $\varphi_R$ is a cofibrant replacement of $R$.

Denote $C_R=(UL\otimesc W,d)$ and let $\psi\colon R\to S$ be any  morphism of $\compul$-modules. Since $\varphi_S$ can be chosen to be surjective, the classical  {\em lifting lemma} ensures the existence of a morphism $C_\psi\colon C_R\to C_S$ making the following diagram commute: 
$$
\xymatrix{C_R\ar[d]_{\varphi_R}\ar[r]^{C_\psi}&C_S\ar[d]^{\varphi_S}\\
R\ar[r]_{\psi}&
S.\\}
$$
In short,  $\cdgmul^0$ admits  quasi-functorial  cofibrant replacements.

Next, recall the standard cylinder of a non necessarily connected $\compul$-module $R$,
$$
\xymatrix{
R \ar@<0.5ex>[r]^(.37){i_0}\ar@<-0.5ex>[r]_(.37){i_1}&\cyl(R)\ar[r]^(.60)p&R \\}
$$
where  
 $\cyl(R)=(R\oplus R'\oplus sR,d)$, with $R'$ a copy of $R$, the differential on $R$ and $R'$ matching that of $R$, and $dsr=-r-sdr$ for  $r\in R$. The maps  $i_o$ and $i_1$ are the inclusions of $R$ and $R'$, respectively,  while  $p(r)=p(r')=r$ and $p(sr)=0$ for all $r\in R$ and  $r'\in R'$. All these maps  respect the given filtrations so they belong to $\cdgmul$, and moreover, they  are all quasi-isomorphisms.

In particular, If $R$ is connected, 
the above construction gives 
 a (to be) very good cylinder of $C_R$ in $\cdgmul^0$.

All of the above allows us to apply Theorem \ref{transfer} and conclude that there exists a left-transferred model structure on $\cdgmul^0$ for which the adjunction $\libc_L\dashv\calker^0$ is Quillen and whose cofibrations and weak equivalences are created by $\libc_L$. We now verify that this model structure coincides with the projective one.

On the one hand, By Corollary \ref{quism}, the weak equivalences in the transferred model structure on $\cdgmul^0$ are precisely the quasi-isomorphisms.

On the other hand, let $f\colon R\to S$ be a fibration in the transferred structure on $\cdgmul^0$, write $\mathfrak{S}=S_{\ge 1}$ and consider the surjective morphism of $\compul$-modules
$$
R\amalg \bigl(UL \otimesc (\mathfrak{S}\oplus d\mathfrak{S})\bigr)\longrightarrow S
$$
which restricts to $f$ on $R$ and is naturally defined on $\bigl(UL \otimesc (\mathfrak{S}\oplus d\mathfrak{S})\bigr)$. The inclusion $R\hookrightarrow R\amalg \bigl(UL \otimesc (\mathfrak{S}\oplus d\mathfrak{S})\bigr)$ is a trivial cofibration in the transferred structure since its image under $\libc_L$ is, see \cite[Prop.~8.10(i)]{bufemutan0}. Hence, the following commutative square
$$
\xymatrix{R\ar[d]\ar@{=}[r]&R\ar[d]^{f}\\
R\amalg \bigl(UL \otimesc (\mathfrak{S}\oplus d\mathfrak{S})\bigr)\ar[r]&
S\\}
$$
admits a lift, which shows that $f$ is surjective in positive degrees.

Conversely, let $f\colon R\to S$ be a map of connected $\compul$-modules that is surjective in positive degrees, and suppose it fits into a commutative square
$$
\xymatrix{P\ar@{^(->}^(.46){\sim}_j[d]\ar[r]^g&R\ar[d]^(.43){f}\\
Q\ar[r]_h&
S\\}
$$
where $j$ is a trivial cofibration in the transferred model structure. Applying the functor $\libc_L$, which is $L\oplus \libc(\,\cdot\,)$ by Proposition \ref{segundaprop}, we obtain a commutative square in $\cdgll^0$ of the form 
$$
\xymatrix{L\oplus \libc(\mathcal{P})\ar@{^(->}^(.46){\sim}_{L\oplus\libc(j)}[d]\ar[r]^{L\oplus\libc(g)}&L\oplus \libc(\mathcal{R})\ar[d]^(.43){L\oplus \libc(f)}\\
L\oplus \libc(\mathcal{Q})\ar[r]_{L\oplus\libc(h)}\ar@{-->}[ur]^\varphi&
L\oplus \libc(\mathcal{S})\\}
$$
where the left vertical map is a trivial cofibration and the right vertical arrow is surjective in positive degrees. Since a trivial cofibration between connected cdgl's  is also a trivial cofibration in the projective structure in $\cdgll^0$, the dotted lift $\varphi$ exists.

Now, for any indecomposable element $a\in \mathcal{Q}$, we must have
$$
\varphi(a)=b+\Gamma,\quad b\in \mathcal{R},\quad \Gamma\,\,\text{cycle in}\,\,\libc^{\ge 2}(\mathcal{R}).
$$
We define the map
$$
\psi\colon Q\longrightarrow R,\quad \psi(a)=b,
$$
which commutes with differentials. Moreover, for $x\in L$ and $a\in Q$,  
$$
\varphi[x,a]=[x,\varphi(a)]=[x,b]+[x,\Gamma].
$$
Since  $[x,a]\in \mathcal{Q}$, $[x,b]\in \mathcal{R}$, and $[x,\Gamma]\in\libc^{\ge 2}(\mathcal{R})$, we find that, 
$\psi(xa)=xb=x\psi(a)$. That is, $\psi$ is a map of $\compul$-modules which fits in the following commutative diagram,
$$
\xymatrix{P\ar@{^(->}^(.46){\sim}_j[d]\ar[r]^g&R\ar[d]^(.43){f}\\
Q\ar[r]_h\ar@{-->}[ur]_\psi&
S.\\}
$$
This shows that $f$ is a fibration in the transferred structure.
\end{proof}

\subsection{Parametrized cdgl spectra are spectra of connected \texorpdfstring{$\compul$}{l}-modules}

\subsubsection{Spectra of $\compul$-modules}
The suspension and desuspension functors $s$ and $s^{-1}$ are mutually inverse equivalences of $\cdgmul$. These are the functors used for developing spectra in this category, which is always assumed to be equipped with the projective model structure. 

\begin{definition}  A {\em spectrum} in $\cdgmul$  is an object  of $\spec(\cdgmul)$ which we denote by $\specmodul$. That is, it is a family of $\compul$-modules $R=\{R^n\}_{n\ge0}$  
together with {\em structure maps}, $\sigma\colon sR^n\to R^{n+1}$ for $n\ge 0$,  or equivalently, their adjoints $\eta\colon R^n\to s^{-1}R^{n+1}$.
A morphism $f\colon R\to S$ of $\specmodul$ is a family $\{f^n\colon R^n\to S^n\}_{n\ge0}$ of maps of $\compul$-modules compatible with the structure maps.  
\end{definition}

\begin{rem}\label{moduequi} Since $s$ and $s^{-1}$ are equivalences, the projective and stable model structures on $\specmodul$ coincide. Furthermore, their  prolongations to spectra define mutually inverse equivalences on $\specmodul$. In addition, as in equation (\ref{equiv}), there is a Quillen equivalence
\begin{equation}\label{equivmodul}
 \xymatrix{ \cdgmul & \specmodul.\ar@<0.75ex>[l]^(.43){\eva_0}
\ar@<0.75ex>[l];[]^(.55){s^{\infty}}\\}
\end{equation}
\end{rem}

 \smallskip
 
Assume now $L$ is connected. Then, upon restriction to  connected $\compul$-modules,  the truncated suspension and desuspension
$$
\xymatrix{ \cdgmul^0 & \cdgmul^0,\ar@<0.75ex>[l]^(.48){\tilde{s}^{-1}}
\ar@<0.75ex>[l];[]^(.49){\tilde{s}}\\}
$$
$$
(\tilde{s}R)_n=\begin{cases}0,& n=0,\\ R_{n-1},& n>0,\end{cases}\qquad (\tilde{s}^{-1}R)_n=\begin{cases} \ker( d\colon R_1\to R_0),& n=0,\\ R_{n+1},& n>0,\end{cases}
$$
 are   no longer equivalences. Nevertheless, they are adjoint functors, $\tilde{s}$ provides a functorial section of $\tilde{s}^{-1}$, and the composition $\tilde{s}\tilde{s}^{-1}(M)$ gives the simply connected cover of $M$.
Moreover, since $\tilde s^{-1}R=(s^{-1}R)^{(0)}$, it follows that $\tilde{s}^{-1}$ preserves weak equivalences and fibrations  so that the adjunction $\tilde{s}\dashv \tilde{s}^{-1}$ becomes a Quillen pair of endofunctors on $\cdgmul^0$. This allows us to define the category $\specmodu$ of spectra over $\cdgmul^0$ to which the general background from \S\ref{specmod} applies, yielding:

\begin{proposition}\label{generalcdgm} 
There are adjoint endofunctors 
$$
\xymatrix{ \specmodu & \specmodu\ar@<0.75ex>[l]^(.45){\mathbf{\tilde s^{-1}}}
\ar@<0.75ex>[l];[]^(.50){\mathbf{\tilde s}}\\}
$$
obtained by prolongation of  $\tilde s$ and $\tilde s^{-1}$. Furthermore, for any $k\ge 0$, there are adjoint functors
$$
\xymatrix{ \cdgmul^0 & \specmodu\ar@<0.75ex>[l]^(.40){\eva_k}
\ar@<0.75ex>[l];[]^(.53){\tilde s^{\infty-k}}\\}
$$
\hfill$\square$
\end{proposition}

\begin{rem}\label{remarkmod} Moreover, by Theorem \ref{prolon}(i), the prolongation of the adjunction  (\ref{conemod}) to spectra induces a Quillen pair
$$
\xymatrix{ \specmodu & \specmodul \ar@<0.75ex>[l]^(.43){\bm{(\,\cdot\,)^{(0)}}} \ar@<0.75ex>[l];[] \\ }
$$
which will later be shown to be a Quillen equivalence.
\end{rem}

\subsubsection{The Quillen equivalence}

For the remainder of this section, we assume that $L$ is connected.

\begin{theorem}\label{mainmod} The adjunction $\libc_L\dashv\calker^0$ induces a left-transferred Quillen pair with respect to the stable structures,
$$
\xymatrix{ \specmodu & \specl, \ar@<0.75ex>[l]^(.43){\calker^0}
\ar@<0.75ex>[l];[]^(.55){\libcb_L}\\}
$$
in which the right adjoint is the prolongation of $\calker^0$ and  $\libc_L\,\tilde s^{\infty-k}=\Sigma_L^{\infty-k}\,\libc_L$ for all $k$.
\end{theorem} 
As in previous sections, we will henceforth use the same notation for functors in the unstable categories and their induced counterparts in the corresponding spectral categories, since no confusion arises and to avoid excessive notation.
\begin{lemma}\label{calkersur} The functor $\calker^0\colon\specl\to\specmodu$ is full and essentially surjective.
\end{lemma}

\begin{proof}
Let $P\in\specmodu$. For each $n\ge 0$ consider the retractive dgl $L\oplus \mathcal{P}^n$ where $\mathcal{P}^n$ is an abelian Lie algebra isomorphic to $P^n$ whose Lie bracket with $L$ is determined by the module structure. From the explicit description of the functor $\lupl$ given in Proposition \ref{lupquasi}, the structure map $\eta_n\colon P^n\to s^{-1}P^{n+1}$ induces a morphism in $\cdgll$ given by:
\begin{equation}\label{calkermap}
\begin{aligned}
\eta_n&\colon L\oplus\mathcal{P}^n\longrightarrow \lupl(L\oplus \mathcal{P}^{n+1})\cong L\oplus \bigl(\mathcal{P}^{n+1}\otimesc (C\oplus dC\oplus \bq dt)\bigr),\\
&\qquad\qquad \quad x\in\mathcal{P}^n\mapsto y\,dt \quad\text{with}\quad \eta(x)=s^{-1}y.\\
\end{aligned} 
\end{equation}
Therefore the sequence $L\oplus\mathcal{P}=\{L\oplus\mathcal{P}^n\}_{n\ge 0}$, together with the structure maps $\{\eta_n\}_{n\ge 0}$, defines a spectrum in $\specl$ such that 
$$
P=\calker^0(L\oplus\mathcal{P}).
$$
Now, given a morphism $h\colon P\to Q$ in $\specmodu$, the map $g\colon L\oplus\mathcal{P}\to L\oplus\mathcal{Q}$, which is $h^n$ on each $\mathcal{P}^n$, defines a morphism in $\specl$ for which $\calker^0(g)=h$.
\end{proof}

\begin{proof}[Proof of Theorem \ref{mainmod}] By Proposition \ref{lupquasi}, and for any   $M=L\oplus K\in\cdgll$, the $\compul$-module $\calker^0(\lupl M)$ is quasi-isomorphic to $(s^{-1}K)^{(0)}=\tilde s^{-1}K^{(0)}$ which is in turn quasi-isomorphic, and therefore weakly equivalent, to $\tilde s^{-1}\calker^0(M)$. Applying Theorem \ref{prolon}(1) 
we obtain a Quillen adjunction as stated, satisfying the required properties. 

We now verify that this adjunction is left-transferred by showing that both, the stable and the transferred structures in $\specmodu$ have the same cofibrations and fibrant objects.

Any stable cofibration of $\specmodu$ is  a cofibration in the transferred structure since ${\libc_L}$ is a left Quillen functor. Conversely, suppose that $f\colon R\to S$ is a morphism  in $\specmodu$ such that ${\libc_L}(f)$ is a cofibration of $\specl$, and recall  the stable and projective structures on  $\specmodu$ have the same cofibrations and thus, the same trivial fibrations. 

Consider a commutative diagram in $\specmodu$
$$
\xymatrix{R\ar[d]_(.46)f\ar[r]&P\ar[d]^h\\
S\ar[r]&
Q\\}
$$
 in which $h$ is a trivial fibration.  By Lemma \ref{calkersur} we can write $h=\calker^0(g)$ for some $g\colon L\oplus\mathcal{P}\to L\oplus\mathcal{Q}$. Therefore, the original square admits a lift if and only if its adjoint
$$
\xymatrix{\libc_L(R)\ar[d]_(.46){\libc_L(f)}\ar[r]&L\oplus\mathcal{P}
\ar[d]^{g}\\
\libc_L(S)\ar[r]&
L\oplus\mathcal{Q}\\}
$$
admits a lift, which it does since $g$ is a trivial fibration in $\specl$. 

Now, let $R\in\specmodu$  be a fibrant object in the transferred structure and let $f$ be an acyclic cofibration of the stable structure. Then $\libc_L(f)$ is an acyclic cofibration and thus $f$ is also an acyclic cofibration in the transferred structure. Hence $R\to 0$ has the right lifting property with respect to $f$ and $R$ is fibrant in the stable structure.

Conversely, assume that $R=\{R_n\}_{n\ge 0}$  fibrant in the stable category. As in Lemma \ref{calkersur}, write $R=\calker^0(L\oplus\mathcal{R})$ where the adjoint structure maps of $L\oplus \mathcal{R}$ 
$$
\eta_n\colon L\oplus\mathcal{R}^n\longrightarrow \lupl(L\oplus \mathcal{R}^{n+1})
$$
are as given in (\ref{calkermap}).  Composing each of these with the equivalence in  Corollary \ref{lupquasicoro}, we obtain a morphism of retractive cdgl's 
 $$
 \mu_n\colon L\oplus\mathcal{R}^n\stackrel{\simeq}{\longrightarrow} L\oplus \tilde s^{-1}\mathcal{R}^{n+1}
  $$
  which is a quasi-isomorphism since  the adjoint structure map $R^n\stackrel{\simeq}{\to} \tilde s^{-1}R^{n+1}$ is itself a quasi-isomorphism. Now, since each $\mathcal{R}^n$ is abelian, an application of the {\em Goldman--Millson Theorem}, see \ref{homotocdgl}, implies that $\mu_n$ is a weak equivalence. Hence, $\eta_n$ is also a weak equivalence and thus, $L\oplus \mathcal{R}\in\specl$ is a fibrant spectrum.
  
  Next, consider a commutative diagram in $\specmodu$ 
$$
\xymatrix{P\ar[d]_(.46)f\ar[r]&R\ar[d]\\
Q\ar[r]&
0\\}
$$
where $f$ is a transferred acyclic cofibration. As before, this square admits a lift if and only if the  square 
$$
\xymatrix{\libc_L(P)\ar[d]_(.46){\libc_L(f)}\ar[r]&L\oplus\mathcal{R}
\ar[d]^{g}\\
\libc_L(Q)\ar[r]&
L\\}
$$
does, which holds since $L\oplus\mathcal{R}$  is fibrant.
\end{proof}

 The first consequence of Theorem \ref{mainmod} is the analogue of Proposition \ref{espnoesp1} for which the same proof applies:
\begin{proposition}\label{espnoesp12} For any $k\ge 0$, the following diagrams of functors commute up to isomorphisms:
$$
\xymatrix{\cdgmul^0\ar[d]_{\tilde s^{\infty-k}}\ar[r]^(.55){\libc_L}& \cdgll \ar[d]^{\Sigma^{\infty-k}_L} & \quad & \cdgmul^0&\cdgll \ar[l]_(.43){\calker^0} 
\\ \specmodu\ar[r]_{\libc_L}&\specl, & \quad & \specmodu\ar[u]^(.40){\eva_k}&\specl\ar[l]^(.46){\calker^0}\ar[u]_(.40){\eva_k}. } 
$$
\hfill$\square$
\end{proposition}

The following tells us that the stable homology of $\compul$-modules indeed characterizes stable equivalences of spectra in this category. 

\begin{definition}\label{homoestamod} Given $R\in\cdgmul^0$, consider the sequence
\begin{equation}\label{secuencia2}
R^0\stackrel{\nu_0}{\longrightarrow} \tilde s^{-1}R^1\to\dots\to \tilde s^{-n}R^n\stackrel{\tilde s^{-n}\nu_n}{\longrightarrow} \tilde s^{-n-1}R^{n+1}\to\dots
\end{equation}
where $\nu_n\colon R^n\to \tilde s^{-1}R^{n+1}$ is the adjoint structure map.
 The {\em stable homology of $R$} is defined as the homology of the colimit of this sequence:
 $$
 H^{\rm st}(R)=H(\varinjlim_n \tilde s^{-n} R^n).
 $$
 Any map $f\colon R\to S$ in $\specmodu$  induces  a morphism between the corresponding sequences as above, and hence a map
$$
H^{\rm st}(f)\colon H^{\rm st}(R)\longrightarrow H^{\rm st}(S).
$$
Since $
H_k(\tilde s^{-n}R^n)\cong H_{n+k}(R^n)
$ for $k\in\bz$ and homology commutes with directed colimits,
 we have:
$$
H^{\rm st}(R)=\oplus_{k\in\bz}H_k^{\rm st}(R)
$$ where 
$$
H^{\rm st}_k(R)=\varinjlim_{n} H_{n+k}(R^n).
$$
\end{definition}

\begin{proposition}\label{homomodes}  For any $M\in\specl$ and any  $R\in\specmodu$ there are natural isomorphisms,
$$
H^{\rm st}(M)\cong H^{\rm st}\bigl(\calker^0(M)\bigr),\qquad 
H^{\rm st}(R)\cong H^{\rm st}\bigl(\libc_L(R)\bigr).
$$ 
\end{proposition}

\begin{proof}  As any spectrum $M\in \specl$ is stably equivalent to a levelwise connected one (see (iii) of Remark \ref{casoconexo}) we may write 
$$
M=\{M^n\}_{n\ge 0} \quad\text{with}\quad M^n=L\oplus K^n.
$$ 
with $K^n$ connected for every $n$. 
Since $\calker^0$ is the prolongation of $\calker^0$, we have
$$
\calker^0(M)=\{K^n\}_{n\ge0},
$$
and the adjoint structure map $\mu_n\colon K^n\to \tilde s^{-1}K^{n+1}$ is simply the restriction of the map $\nu_n\colon K^n\to s^{-1}K^{n+1}$ from Definition \ref{sorpresa2} to the given codomain. This is well defined as every $K^n$ is connected. Hence, although  the corresponding sequences (\ref{secuencia})  and (\ref{secuencia2}) slightly differ they have the same direct limit and 
 the first isomorphism follows.

For the second, and in light  of Remark \ref{uyuy}, we may assume that $\libc_L$ is the prolongation functor of $\libc_L$. Thus, levelwise  and with the notation in Proposition \ref{segundaprop},
$$
\libc_L(R)^n=L\oplus \libc(\mathcal{R}^n).
$$
As $R$ is connected, this retractive cdgl is cofibrant and, by Theorem \ref{suspensionex}, the structure map $
\susl \bigl(\libc_L(R)\bigr)^n\to \bigl(\libc_L(R)\bigr)^{n+1}
$ is weakly equivalent to a morphism
$$
L\oplus\libc( s\mathcal{R}^n)\longrightarrow L\oplus\libc(\mathcal{R}^{n+1}).
$$
Note that, since the differential on $\libc(\mathcal{R})$ is already linear, the induced differential on the suspension is the suspended differential.
The adjoint map, see Definition \ref{falta}, is then weakly equivalent to   the morphism of retractive cdgl's
$$
L\oplus\libc(\mathcal{R}^n)\longrightarrow L\oplus s^{-1}\libc(\mathcal{R}^{n+1})
$$
whose restriction
$
 \libc(\mathcal{R}^n)\to s^{-1}\libc(\mathcal{R}^{n+1})
$
is induced by the adjoint structure map $R^n\to \tilde s^{-1}R^{n+1}$ of $R$. 
The same argument followed in the proof of Proposition \ref{estahomolin} shows that the stable homology of $\libc_L(R)$ is the homology of the colimit of the sequence
$$
R^0\to \tilde s^{-1}R^1\to\dots\to \tilde s^{-n}R^n\to \tilde s^{-n-1}R^{n+1}\to\dots
$$
which is precisely $H^{\rm st}(R)$.
\end{proof}

In view of Theorem \ref{mainmod}, since $\libc_L$ reflects weak equivalences, the preceding result  implies the following: 

\begin{corollary}\label{elquefalta}
A map $f\colon R\to S$ of $\specmodu$ is a stable equivalence if and only if it induces an isomorphism in stable homology.\hfill$\square$
\end{corollary}

Another  consequence is:

\begin{corollary}\label{elquefalta2}
The functor $\calker^0$ reflects weak equivalences.\hfill$\square$
\end{corollary}

We conclude with the main result of this section.

\begin{theorem}\label{segundo} The adjunction
$$
\xymatrix{ \specmodu & \specl, \ar@<0.75ex>[l]^(.43){\calker^0}
\ar@<0.75ex>[l];[]^(.55){\libc_L}\\}
$$
is a Quillen equivalence. 
\end{theorem}

\begin{proof}
By Theorem \ref{mainmod} $\libc_L$ reflects weak equivalences. Therefore, in light of  \cite[Cor.~1.3.16(b)]{ho0}, it suffices to show that for each spectrum $M\in\specl$ the counit of the adjunction 
$$
\varepsilon_M\colon\libc_L\,\calker^0(M)\longrightarrow M
$$
is a stable equivalence, that is, it induces an isomorphism in stable homology. Again, there is no loss of generality assuming  $M$  connected. Write 
$$
M=\{M^n\}_{n\ge 0} \quad\text{with}\quad M^n=L\oplus K^n
$$
so that
$$
\calker^0(M)=\{K^n\}_{n\ge0}.
$$
By Remark \ref{uyuy}, we may assume that $\libc_L$ is the prolongation of the unstable $\libc_L$. Then, using the notation in Proposition \ref{segundaprop},
$$
\libc_L\calker^0(M)=\{L\oplus \libc(\mathcal{K}^n)\}_{n\ge 0},
$$
and the counit of the adjunction  is, levelwise, the family of morphisms in $\cdgll$,
$$
\varepsilon_M^n\colon L\oplus \libc(\mathcal{K}^n)\longrightarrow L\oplus\libc( K^n),\quad n\ge 0,
$$
which extends the isomorphism $\mathcal{K}^n\stackrel{\simeq}{\to} K^n$. To conclude, apply Proposition
\ref{homomodes} which shows that
$$
H^{\rm st}\bigl(\libc_L\calker^0(M)\bigr)\cong H^{\rm st}\bigl(\calker^0(M)\bigr)\cong H^{\rm st}(M).
$$
Under this isomorphism, $H^{\rm st}(\varepsilon_M)$ is the identity.

\end{proof} 
We now provide a computable description of the equivalence established above. 
 
\begin{corollary}\label{corocomputable2} The equivalence induced by the preceding result
$$
 \Ho\specl\,{\cong}\,\Ho\specmodu
$$
assigns to each $L$-spectrum $M$ with levelwise connected model $(L\amalg\libc(W),d)$ a spectrum of $\compul$-modules of the form
$$
(UL\otimesc W,d_1)=\{(UL\otimesc W^n, d_1)\}_{n\ge 0}
$$
where the differentials and structure maps are naturally inherited from the  indecomposable or linear reduction of $(L\amalg\libc(W),d)$
\end{corollary}

\begin{proof} By (iv) of Remark \ref{casoconexo}, any $L$-spectrum $M$ admits indeed a levelwise connected model. As usual, write
$$
(L\amalg\libc(W),d)=\{(L\amalg\libc(W^n),d)\}_{n\ge 0}= \{L\oplus (\libc(T^n),d)\}_{n\ge 0}
$$
and let
$$
(L\amalg W,d_1)=\{(L\amalg W^n,d_1)\}_{n\ge 0}=\{L\oplus (T^n,d_1)\}_{n\ge 0}
$$
 be its indecomposable reduction(see (ii) Definition \ref{indelineal}).
 By Corollary \ref{elquefalta2} the functor  $\calker^0$ reflects weak equivalences so  the equivalence in the statement associates to  $M$, up to isomorphism, the spectrum of $\compul$-modules  $\calker^0(L\amalg W,d_1)$. Since $\calker^0$ is a prolongation,
 $$
 \calker^0(L\amalg W,d_1)=\{\calker^0\bigl(L\oplus (T^n,d_1)\bigr)\}_{n\ge 0}=\{(T^n,d_1 )\}_{n\ge 0}.
 $$
 On the other hand, the adjoint $n$th structure map
 $$
 (T^n,d_1)\longrightarrow \tilde s^{-1}(T^{n+1},d_1)
  $$
  is, by definition, the composition of 
  $$
  \calker^0(\bar\eta_n)\colon 
\calker^0\bigl(L\oplus (T^n,d_1)\bigr)\longrightarrow\calker^0\lupl\bigl(L\oplus(T^{n+1},d_1)\bigr),
$$
where $\bar\eta_n$ is the corresponding structure map of the indecomposable reduction (see (ii) of Definition \ref{indelineal}), with the weak equivalence
$$
\calker^0\lupl\bigl(L\oplus(T^{n+1},d_1)\bigr)\stackrel{\sim}{\longrightarrow} \tilde s^{-1}\calker^0\bigl(L\oplus(T^{n+1},d_1)\bigr)
$$
induced by the weak equivalence of Corollary \ref{lupquasicoro}.

To finish write
$$
(T^n,d_1)\cong (UL\otimesc W^n,d_1),\quad n\ge 0.
$$
Alternatively, one could apply Remark \ref{ejemplo1} levelwise  on  the linear reduction $(L\amalg\libc(W),d_1)$  to obtain
$$
M\sim (L\amalg\libc(W),d_1)\sim \libc_L(UL\otimesc W,d_1).
$$
\end{proof}
\section{From spectra of connected \texorpdfstring{$\compul$}{l}-modules to \texorpdfstring{$\compul$}{l}-modules}

 The correspondence between modules over a ring and symmetric spectra of connected modules is classical and well understood \cite[Prop.~4.7]{ship}. Nonetheless, for completeness, we include a brief proof in the complete, non-symmetric case. In what follows, $\compul$ may be replaced by any  filtered differential graded algebra and everything also holds in the non-complete setting. 


\begin{definition}\label{pairclas} Define the    functors
$$
\xymatrix{ \specmodu & \cdgmul, \ar@<0.75ex>[l]^(.54){\mathcal{C}}
\ar@<0.75ex>[l];[]^(.44){\mathcal{D}}\\}
$$
as follows: for $R\in\specmodu$, set
$$
\mathcal{D}R=\varinjlim_n s^{-n+1}R^n
$$
where the colimit is taken over the maps 
$$
s^{-n+1}R^n\stackrel{s^{-n+1}(\eta)}{\longrightarrow}s^{-n+1}(\tilde s^{-1}R^{n+1})\hookrightarrow s^{-n+2}R^{n+1}
$$ with $\eta$ denoting the adjoint structure map.

On the other hand,
 $\mathcal{C}$ may be viewed as the connected cover of $s^{\infty-1}$: for $R\in\cdgmul$ define $(\mathcal{C}R)^n=(s^{n-1}R)^{(0)}$ with structure maps  given by the injections $\tilde s(s^{n-1}R)^{(0)}\hookrightarrow (s^{n}R)^{(0)}$. 
 
 It is straightforward to check that  $\mathcal{D}$ is left adjoint to $\mathcal{C}$.
\end{definition}

 \begin{rem}
 The informed reader may have noticed an extra suspension and desuspension, respectively, in the usual definitions of $\mathcal{D}$ and $\mathcal{C}$. These adjustments are necessary to account for the degree shift between the stable homology of any $L$-spectrum and the rational stable homotopy of the $B$-spectrum it represents, see Theorem \ref{homotohomolo} and its connected version in (ii) of Remark \ref{casoconexo}.
 \end{rem}

\begin{proposition}\label{homolohomolo} For any $R\in\specmodu$ there is a natural isomorphism
$$
H^{\rm st}(R)=s^{-1} H(\mathcal{D}R).
$$
\end{proposition}
\begin{proof} Observe that
$$
H_k(\tilde s^{-n}R^n)\cong H_{k+n}(R^n)=H_k(s^{-n}R^n)
$$
for any $k\in\bz$ and $n\ge 0$.
Hence,
since homology commutes with directed colimits,
$$
H^{\rm st}(R)=\varinjlim_n H(\tilde s^{-n}R^n)\cong  \varinjlim_n H(s^{-n}R^n)\cong s^{-1} H(\mathcal{D}R).
$$
\end{proof}

Under the projective model structure in $\cdgmul$ we have:

\begin{proposition}\label{casifin} The adjunction $\mathcal{D}\dashv\mathcal{C}$ is a Quillen equivalence.
\end{proposition}

\begin{proof}
Let $f\colon R\to S$ be a map in $\specmodu$. If $f$  is  a cofibration, then it is in particular a levelwise cofibration and thus, each map $s^{-n}R^n\to s^{-n}S^n$ is also a cofibration in $\cdgmul$. Since  cofibrations of any combinatorial model category are closed under filtered colimits the map $\mathcal{D}f$ is also a cofibration. 

On the other hand, by Proposition \ref{homolohomolo}, $H^{\rm st}(f)$ is an isomorphism if and only if $H(\mathcal{D}f)$ is. By Corollary \ref{elquefalta}, this is equivalent to $f$ being a weak equivalence if and only if $\mathcal{D}f$ is. 

Hence, $\mathcal{D}\dashv\mathcal{C}$ is a Quillen pair in which $\mathcal{D}$ reflects weak equivalences.
Finally, it is also straightforward to verify that for any $\compul$-module $R$, $\mathcal{D}\mathcal{C}(R)=\varinjlim_n R^{(-n)}\cong R$ so the counit of the adjunction is an isomorphism. By \cite[Cor.~1.3.16(b)]{ho0}, the proposition follows.
\end{proof}

\begin{rem}  As in  \cite[Prop.~4.9]{ship}, we may alternatively decompose the Quillen equivalence from Proposition \ref{casifin} into two separate Quillen equivalences, corresponding to those in  Remark \ref{remarkmod} and  (\ref{equivmodul}): 
$$
\xymatrix{ 
\specmodu & 
\specmodul \ar@<0.75ex>[l]^(.43){(\,\cdot\,)^{(0)}}
           \ar@<0.75ex>[l];[]\ar@<-0.80ex>[r]_(.50){\eva_0} & 
\cdgmul.
       \ar@<-0.60ex>[l]_{s^{\infty}} \\
}
$$
  \end{rem}
  
  \section{The conclusion and practical aspects}
The main results of the preceding sections can be summarized as follows. Let  $B$ be  a reduced simplicial set with unique $0$-simplex $b$, set $L=\lasu^*_B=\lasu_B/(b)$, and consider the functor
  $$
 \Psi\colon  \Ho\spec_B \longrightarrow\Ho\cdgmul
  $$
induced by the composition of the following sequence of Quillen adjunctions, 
\begin{equation}\label{composmono}
\xymatrix{ 
\spec_B 
  \ar@<0.75ex>[r]^(.50){\lasub^*} 
& 
\specl 
  \ar@<0.75ex>[l]^(.45){\bm{\langle\,\cdot\,\rangle}}
    \ar@<-0.80ex>[r]_(.45){\calker^0} 
& 
\specmodu 
  \ar@<-0.60ex>[l]_(.50){\libc_L} 
  \ar@<0.50ex>[r]^(.45){\mathcal{D}} 
& 
\cdgmul, 
  \ar@<0.90ex>[l]^(.55){\mathcal{C}} 
}
\end{equation}
Since all the functors in this sequence have been shown to preserve all weak equivalences, regardless of fibrancy or cofibrancy, we may write directly,
$$ \Psi(X)=\mathcal{D}\calker^0(\lasub^*_X),\quad X\in \spec_B.
$$
\begin{theorem}\label{mainmain} The  induced functor
$$
 \Psi_\bq\colon  \Ho\spec_B^\bq\stackrel{\cong}{\longrightarrow}\Ho\cdgmul.
  $$
is an equivalence of categories.
\end{theorem}

\begin{proof} By Theorem \ref{segundo} and Proposition \ref{casifin} the last two pairs in (\ref{composmono}) are Quillen equivalences. Theorem   \ref{primeroqq} shows  that the first pair also induces a Quillen equivalence on rational spectra.
\end{proof}

  Among the various formulations of the functor $\Psi$ we highlight the following, the first having a distinctly geometric character whereas the second is particularly well suited to computations.

\begin{theorem}\label{teorema1} Up to weak equivalence, the functor $\Psi$ assigns to each spectrum $X\in \spec_B$, with fiber spectrum $X^b\in\spec$, the $\compul$-module $\varinjlim_n s^{-n+1}(K^n)^{(0)}$ where $K^n$ is a suitable Lie model of $X_n^b$ whose module structure is given by a model $L\to L\oplus K^n\to L$ of $B\to X_n\to B$.
\end{theorem}

\begin{proof} 
Write $\lasub^*_X=M$ where each $M^n=L\oplus K^n$  is, by Remark \ref{uyuy},  weakly equivalent to $\lasu^*_{X_n}$. In particular, $\langle M^n\rangle$ is weakly equivalent to the $\bq$-completion $(X_n)^\acento_\bq$.

On the one hand, 
$$
\Psi(X)=\mathcal{D}\calker^0(M)=\varinjlim_n s^{-n+1}(K^n)^{(0)}.
$$

On the other hand, in view of (\ref{pullback}) and (\ref{realizafibra}), $\bm\langle M\bm\rangle_n$ fits into a (homotopy) pullback of the form
$$
\xymatrix{
\bm\langle M\bm\rangle_n\ar[d]\ar[r]&B\ar[d]\\
X^\acento_\bq\ar[r]_{ p^\acento_\bq}& B^\acento_\bq.
}
$$
whose fiber over $b$ is 
$(X_n^b)^\acento_\bq$, which is weakly equivalent to $ \langle K^n\rangle
$.
\end{proof}

\begin{theorem}\label{teocomputable} 
Let $(L\amalg\libc(W),d)$ be a model of $X\in\specb$.   Then, $\calker^0\lasub_X^*$ is weakly equivalent to a $\compul$-module spectrum of the form
$$
(UL\otimesc W, d_1)=\{(UL\otimesc W^n, d_1)\}_{n\ge 0}
$$
where  the differentials and structure maps are naturally inherited from the  indecomposable (or linear) reduction of $(L\amalg\libc(W),d)$. In particular,
$$
\Psi(X)=UL\otimesc (\varinjlim_n s^{-n+1} W^n).
$$
\end{theorem}

\begin{proof}
The first assertion follows directly from Corollary \ref{corocomputable2}. Hence,
 $$
 \Psi(X)=\mathcal{D}( UL\otimesc W, d_1),
 $$
  which, by definition, has the desired form.
\end{proof}
This explicit and computation-oriented formulation  of $\Psi$ allows to determine effectively all the entries in the bijective correspondence it establishes between homotopy invariants of $\compul$-modules and rational spectra. We now present some illustrative examples.

We first prove that retractive suspension and loops  of $B$-spectra correspond, as expected, to suspension and desuspension  of modules, respectively.

\begin{proposition}\label{suslup}For any $X\in\spec_B$,
 $$
 \Psi(\suseb X)=s\Psi(X)\quad\text{and}\quad\Psi(\lupeb X)=s^{-1}\Psi(X).
 $$
\end{proposition}

\begin{proof}
Let $(L\amalg\libc(W),d_1)$ be a linear  model of $X$ which, by (iv) of Remark \ref{casoconexo}, may be assumed  levelwise connected. Example  \ref{suspensionl} shows that $(L\amalg \libc(sW),d_1)$ is a linear model of  $\susel\lasub_X^*$ which is stably equivalent to $\lasub^*_{\suseb X}$. In other words, $(L\amalg \libc(sW),d_1)$ is a model of $\susb X$ and a direct application of Theorem \ref{teocomputable} yields the first assertion. For the second, recall that $\lupeb$ and $s^{-1}$ are the respective inverses of $\suseb$ and $s$ in $\Ho\spec_B$ and $\Ho\cdgmul$ respectively.
\end{proof}

Next, we show that rational stable maps in $\Ho\spec_B$ correspond to the bifunctor $\ext$ in $\Ho\cdgmul$.

\begin{proposition}\label{extfun} For any $X,Y\in\spec_B$ we have a natural isomorphism of graded vector spaces
$$
\{X,Y\}^B\otimes\bq\cong \ext_{\compul}\bigl(\Psi(X),\Psi(Y)\bigr).
$$
\end{proposition}

\begin{proof} Observe that for any $R,S\in\cdgmul$ and any integer $k$,  $ \ext_{\compul}(R,S)_k=\Hom_{\Ho\cdgmul}(R,s^kS)$. The proposition follows from the previous result and Theorem \ref{mainmain}.
\end{proof}

We continue with an important example. 
\begin{proposition}\label{esfera}
$\Psi(\bs_B)\cong\compul$
\end{proposition}
\begin{proof}
Recall from Section \ref{prelimi} that $\bs_B=B\times\bs$. For convenience we replace $\bs_B$ by  its cofibrant replacement  $B\vee \bs$ where each structure map is the weak equivalence $\susb (B\vee S^n)\to B\vee\Sigma S^n=B\vee S^{n+1}$. Then, $\lasub^*_{B\vee \bs}$ is stably equivalent to the free $L$-spectrum $(L\amalg\libc(x),d)=\{(L\amalg\libc(x^{n-1}),d)\}_{n\ge 0}$  where $x^{n-1}$ is a cycle of degree $n-1$. The $n$th structure map is (weakly equivalent to) the isomorphism $(L\amalg\libc(sx^{n-1}),d)\to (L\amalg\libc(x^{n}),d)$ sending $sx^n$ to $x^{n+1}$. The indecomposable reduction of this spectrum is then 
$$
(L\amalg x,d)=\{(L\amalg x^{n-1},d)\}_{n\ge 0}=\{L\oplus (UL\otimes x^{n-1})\}_{n\ge0}.
$$
Using Remark \ref{ya vere}, a direct inspection shows that the $n$th adjoint structure map, composed with the weak equivalence of Corollary \ref{lupquasicoro}, is the morphism
$$
L\oplus (UL\otimesc x^{n-1})\longrightarrow L\oplus s^{-1}(UL\otimesc x^{n}),\quad \Phi\otimes x^{n-1}\mapsto s^{-1}(\Phi\otimes x^n).
$$
Thus, applying Theorem \ref{teocomputable} we obtain
$$
\Psi(\bs_B)=\mathcal{D}(UL\otimesc x)= UL\otimesc (\varinjlim_n s^{-n+1} x^{n-1})=UL\otimesc sx^{-1}\cong \compul.
$$
\end{proof}

On the other hand, parametrized homotopy groups correspond to homology:

\begin{proposition} \label{homotospec} For any $X\in\spec_B$,
$$
\pi^{\rm st}(X)\otimes\bq\cong H\bigl(\Psi(X)\bigr).
$$
\end{proposition}

\begin{proof} The statement follows either by applying Propositions \ref{extfun} and \ref{esfera} together with formula (\ref{stamap}), or directly, using (ii) of Remark \ref{casoconexo} along with Propositions \ref{homomodes} and \ref{homolohomolo}:
$$
\pi^{\rm st}(X)\otimes\bq\cong s\,H^{\rm st}(\lasub^*_X)\cong s\,H^{\rm st}(\calker^0\lasub^*_X)\cong H(\mathcal{D}\calker^0\lasub^*_X)=H\bigl(\Psi(X)\bigr).
$$
\end{proof}
  We also describe how  the change of base adjunction of parametrized spectra translate  into restriction and extension of scalars in the module setting. 
  
  Let $f\colon B\to B'$ be a map of reduced simplicial sets and  write $\varphi\colon L\to  L'$ for the induced map $\lasu_f^*\colon\lasu_B^*\to\lasu_{B'}^*$. Consider the functor
 $$
 \widetilde f_!\colon \specb\longrightarrow\spec_{B'}
 $$
provided by Proposition \ref{changespec} and denote by 
$$
\widehat{\varphi}_!\colon \Ho\cdgmul\longrightarrow \Ho\cdgmulpri
$$
the {\em derived extension of scalars} or  {\em derived base change} induced by the morphism of complete dga's
$$
\widehat{U\varphi}\colon  \compul\longrightarrow \widehat{UL'}.
$$
In fact,  $\widehat{UL'}$ is a $(\widehat{UL'},\compul)$-bimodule via $\widehat\varphi$ and
$$
\widehat{\varphi}_!(R)=\widehat{UL'}\otimesc_{\compul}\,S
$$
where $S\stackrel{\simeq}{\to}R$ is any cofibrant replacement of $R$. In other words,
$
\widehat{\varphi}_!(R)=\widehat{UL'}\otimesc_{\compul}^{\mathbf L}\,R$.

  \begin{proposition}\label{cambiobase}
  The following diagram commutes:
  $$
\xymatrix{
\Ho\specb\ar[d]_{\widetilde f_!}\ar[r]^(.40){\Psi}&\Ho\cdgmul\ar[d]^{\widehat\varphi_!}\\
\Ho\spec_{B'}\ar[r]_(.40){\Psi}&\Ho\cdgmulpri.
}
$$
\end{proposition}

\begin{proof}
First we check that, for any  morphism $\varphi\colon L\to L'$ of connected cdgl's, the following diagram commutes:
$$
\xymatrix{
\Ho\specl\ar[d]_{\varphi_!}\ar[r]^(.40){\mathcal{D}\calker^0}&\Ho\cdgmul\ar[d]^{\widehat\varphi_!}\\
\Ho\spec_{L'}\ar[r]_(.40){\mathcal{D}\calker^0}&\Ho\cdgmulpri.
}
$$
Indeed let $(L\amalg \libc(W),d)$ be a model of any given $L$ spectrum. By Corollary  \ref{corocomputable2},
$$
\mathcal{D}\calker^0(L\amalg \libc(W),d)\cong\mathcal{D}(UL\otimesc W, d_1)=UL\otimesc(\varinjlim_n s^{-n+1}W^n)
$$
which is already a cofibrant $\compul$-module. Hence,
$$
\begin{aligned}
&\widehat\varphi_!\mathcal{D}\calker^0(L\,\amalg\, \libc(W),d)\cong UL'\otimesc (\varinjlim_n s^{-n+1}W^n)\\&\cong \mathcal{D}\calker^0(L'\,\amalg \,\libc(W),d)\cong \mathcal{D}\calker^0\varphi_!(L\,\amalg\, \libc(W),d)
\end{aligned}
$$
The claim follows by combining this with Proposition \ref{natural1}.
\end{proof}

 On the other hand, consider the \emph{restriction of scalars} functor
$$
\widehat\varphi^{\,*}\colon \Ho\cdgmulpri\longrightarrow \Ho\cdgmul
$$
  induced by $\widehat{U\varphi}$. Since  $\widehat\varphi_!\dashv {\widehat\varphi}^{\,*}$ and $f_!\dashv f^*$ are pairs of adjoint functors, upon passing  to the rational stable categories of spectra, the preceding result reads:

  \begin{corollary}\label{coronatural} The following diagrams commute:
 $$
\xymatrix{\Ho\spec_B^\bq\ar[d]_{f_!}\ar[r]^(.43){\Psi}_(.43)\cong& \Ho\cdgmul \ar[d]^{\widehat\varphi_!} & \quad & \Ho\specb^\bq \ar[r]^(.43){\Psi}_(.43)\cong&\Ho\cdgmul
\\ \Ho \spec_{B'}^\bq\ar[r]_(.43){\Psi}^(.43)\cong&\Ho\cdgmulpri & \quad & \Ho \spec_{B'}^\bq\ar[r]_(.43){\Psi}^(.43)\cong\ar[u]^{f^*}&\Ho\cdgmulpri \ar[u]_(.45){\widehat\varphi{\,*}}. } 
$$
\hfill$\square$
\end{corollary}

We finish with the following observation:

\begin{rem}
In the case $L=0$, the classical equivalence between the stable category of  $H\bq$-modules  and the derived category of rational chain complexes (see for instance \cite[Thm.~1.1]{ship} or \cite[Thm.~5.1.6]{schship}) takes a particularly computable form.  In this context, the sequence (\ref{composmono}) becomes
\begin{equation}\label{compos2}
\xymatrix{ 
\spec^{\bq} 
  \ar@<0.75ex>[r]^(.40){\lasub^*} 
& 
\spec(\cdgl) 
  \ar@<0.75ex>[l]^(.58){\bm{\langle\,\cdot\,\rangle}}
    \ar@<-0.80ex>[r]_(.52){\mathcal{U}^0} 
& 
\spec(\catvect^0) 
  \ar@<-0.60ex>[l]_(.50){\libc} 
  \ar@<0.50ex>[r]^(.58){\mathcal{D}} 
& 
\catvect, 
  \ar@<0.90ex>[l]^(.42){\mathcal{C}} 
}
\end{equation}
where the pair $\libc\dashv \mathcal{U}^0$ is induced by the complete free functor and the connected cover of the forgetful functor.
by Theorem \ref{teocomputable}, the resulting equivalence
$$
\Ho(\spec^\bq)\cong \Ho(\catvect),
$$
assigns to each spectrum $X\in\spec$ the chain complex $\varinjlim_ns^{-n+1}W^n$ where $(\libc(W^n),d)$ is a Lie model of $X_n$.

Conversely, since any $(V,d)\in \catvect$ is weakly equivalent to its homology  $(H,0)$, the associated rational spectrum is then the rationalization of the $\cdgl$ spectrum $ \libc\mathcal{C}(H,0)$, which we denote by $(\libc(s^{\infty-1} H)^{(0)},0)$. Note that this spectrum is stably equivalent to the non-truncated one $(\libc(s^{\infty-1} H),0)$. Finally, one checks that its realization, which corresponds to $(V,d)$ under the equivalence (\ref{compos2}),   is a generalized Eilenberg-MacLane spectrum $X$ satisfying
$$
\pi_n^{\rm st}(X)\otimes\bq=H_n(V,d)
$$
 for all $n\in\bz$. In other words,
$$
X=\vee_i\,\Sigma^{k_i}H\bq,
$$
where the index $i$ runs over a basis of $H$ whose elements have degree $k_i$.
\end{rem}

  \section{The strong monoidal character}
  
 As in our general framework let $B$ be a reduced simplicial set with unique $0$-simplex $b$, and denote $L=\lasu^*_B=\lasu_B/(b)$. We show here that the passage from $B$-spectra to $\compul$-modules induces a strong monoidal functor between the respective derived categories. In particular, the stable homotopy category of rational $B$-spectra is equivalent, as strong monoidal categories, to the derived category of $\compul$-modules.
 
 \medskip
 
    Recall that given any dgl $L$, and as for any differential graded Hopf algebra, every left $UL$-module $R$ induces a  right module structure via the antipode map $a\colon UL\to UL$. This is the unique anti-algebra morphism satisfying $a(x)=-x$ for $x\in L$. By anti-algebra map we mean that $a(\phi\psi)=(-1)^{|\phi||\psi|}a(\psi)a(\phi)$ for any  homogeneous elements $\phi,\psi\in UL$. The right action is then defined by $r\cdot\phi=a(\phi)\cdot r$, for $\phi\in UL$ and $r\in R$. This shows that the category of right and left $UL$-modules are equivalent although these two actions do not define a bimodule structure. Therefore, the balanced tensor product $\otimes_{UL}$ does not yield a bifunctor on $\cdgm_{UL}$ with values in the same category and thus, it does not provide a  monoidal structure.  Nevertheless, $\cdgm_{UL}$ is indeed symmetric monoidal under the tensor product over $\bq$ with the action on $R\otimes S$ defined via the diagonal $UL\to UL\otimes UL$. All of the above extends to the complete setting so that $\cdgmul$ is symmetric monoidal with respect to the complete tensor product $\otimesc$.

We prove:

\begin{theorem}\label{monoidal2} The functor $\Psi$ is strong monoidal. In particular $\Psi_\bq$ is a strong monoidal equivalence.
\end{theorem}

We begin by modeling the external and internal smash product in the unstable setting. Let $L$ and $L'$ be cdgl's.   In view of Remarks \ref{semifreeuu}
and \ref{ejemplo1} together with the standard decomposition (\ref{esencialk}), it follows that any free retractive cdgl over $L\oplus L'$ is of the form
\begin{equation}\label{ueue}
\bigl((L\oplus L')\,\amalg\,\libc(Z),d\bigr)\cong (L\oplus L')\oplus \bigl(\libc\bigl((UL\otimes UL')\otimesc Z\bigr),d\bigr).
\end{equation}

Now let $A$ and $B$ be reduced simplicial sets, and let $X\in\reta$ and $Y\in \retb$ be connected  retractive simplicial sets over  $A$ and $B$, respectively. Denote by $L=\lasu_A^*$ and $L'=\lasu_B^*$.

By Corollary \ref{cofirecon}, choose connected retractive models $(L\amalg\libc(V),d)$ and $(L'\amalg\libc(W),d)$ of  $\lasu^*_X$ and $\lasu^*_Y$ respectively. Again by Remark \ref{ejemplo1}, express the retractive linear part of both models, see Definition \ref{libre}, as
$$
(L\amalg\libc(V),d_1)=L\oplus\libc(UL\otimesc V,d_1),\quad (L'\amalg\libc(W),d_1)=L'\oplus \libc(UL'\otimesc W,d_1),
$$
and, for each $v\in V $ and $w\in W$, write
$$
d_1v=\sum_i\alpha_i\otimes v_i,\quad d_1w=\sum_j\beta_j\otimes w_j.
$$

\begin{proposition}\label{external} A retractive model of $X\,\bar\wedge\, Y$ over $A\times B$
 is given by
$$
\bigl((L\oplus L')\amalg\libc\bigl(s(V\otimes W)\bigr),d\bigr)\cong(L\oplus L')\oplus \bigl(\libc\bigl((UL\otimes UL')\otimesc  s(V\otimes W)\bigr),d\bigr)
$$
Furthermore, the linear part of $d$ in the right-hand side is
$$
d_1s(v\otimes w)=-\sum_i\alpha_i\otimes s(v_i\otimes w)-\sum_j (-1)^{|v|+|v||\beta_j|}\,\beta_j\otimes s(v\otimes w_j).
$$

\end{proposition}

\begin{rem}\label{producto}
In what follows we  use the general result from \cite[Cor.~12.12]{bufemutan0}: let $S,T\in\catss^*$ be connected pointed simplicial sets with Lie models 
$$
(\libc(\mathcal{S}),d)\stackrel{\simeq}{\to}\lasu^*_S\quad\text{and}\quad (\libc(\mathcal{T}),d)\stackrel{\simeq}{\to}\lasu^*_T
$$
 respectively. Then, $S\times T$ admits a Lie model of the form
$$
(\libc(\mathcal{S}\oplus\mathcal{T}\oplus s(\mathcal{S}\otimes\mathcal{T}),d)\stackrel{\simeq}{\longrightarrow}\lasu^*_{S\times T}.
$$
whose differential extends those in $(\libc(\mathcal{S}),d)$ and $(\libc(\mathcal{T}),d)$ and such that the map
\begin{equation}\label{producto3}
(\libc(\mathcal{S}\oplus\mathcal{T}\oplus s(\mathcal{S}\otimes\mathcal{T}),d)\stackrel{\simeq}{\longrightarrow}(\libc(\mathcal{S}),d)\oplus (\libc(\mathcal{T}),d),
\end{equation}
which restricts to the identity on $\mathcal{S}$ and  $\mathcal{T}$, and vanishes on $s(\mathcal{S}\otimes\mathcal{T})$, is also a quasi-isomorphism.

Furthermore, define a new grading on $\libc\bigl(\mathcal{S}\oplus\mathcal{T}\oplus s(\mathcal{S}\otimes\mathcal{T})\bigr)$ by setting $|\mathcal{S}|=|\mathcal{T}|=1$ and $|s(\mathcal{S}\otimes\mathcal{T})|=2$, extended bracket-wise. Then, for each $x\in\mathcal{S}$ and $y\in\mathcal{T}$,
\begin{equation}\label{producto1}
ds(x\otimes y)=-s(d_1x\otimes y)-(-1)^{|x|}s(x\otimes d_1 y)-[x,y]+\Gamma,
\end{equation}
where $\Gamma$ has new degree at least 3.

Also, observe that the kernel of this map can be written as 
$$
\libc\bigl(UE\otimesc ([\mathcal{S},\mathcal{T}]\oplus s(\mathcal{S}\otimes\mathcal{T})\bigr)
$$
where $E=\libc(\mathcal{S})\oplus \libc(\mathcal{T})$. In other words,  
\begin{equation}\label{producto2}(\libc(\mathcal{S}\oplus\mathcal{T}\oplus s(\mathcal{S}\otimes\mathcal{T}),d)\cong  (\libc(\mathcal{S}),d)\oplus (\libc(\mathcal{T}),d) \oplus\bigl( \libc\bigl(UE\otimesc ([\mathcal{S},\mathcal{T}]\oplus s(\mathcal{S}\otimes\mathcal{T})\bigr),d\bigr)
\end{equation}

\end{rem}

\begin{proof}[Proof of Proposition \ref{external}]


By the above remark, if we 
write $L=(\libc(U),d)$ and $L'=(\libc(U'),d)$, we obtain Lie models of $A\times B$, $X\times B$, $A\times Y$ and $X\times Y$ of the form
$$
(\libc(\mathcal{U}),d)=\bigl(\libc\bigl(U\oplus U'\oplus s(U\otimes U')\bigr),d\bigr),
$$
$$
 \bigl(\libc\bigl(\mathcal{U}\oplus V\oplus s(V\otimes U')\bigr),d\bigr), \quad  \bigl(\libc\bigl(\mathcal{U}\oplus W\oplus s(U\otimes W)\bigr),d\bigr),
$$
and
$$
P=\bigl(\libc\bigl(\mathcal{U}\oplus V\oplus W\oplus s(V\otimes U')\oplus s(U\otimes W)\oplus s(V\otimes W)\bigr),d\bigr)
$$
respectively.

Since the functor $\lasu^*$ preserves colimits, a  Lie model of the pushout $(X\times B)\cup_{A\times B}(A\times Y)$ has the form
$$
\bigl(\libc\bigl(\mathcal{U}\oplus V\oplus W\oplus s(V\otimes U')\oplus s(U\otimes W)\bigr),d\bigr).
$$
Similarly, a retractive Lie model of the external smash $X\bar\wedge Y$, which is the pushout in (\ref{pushout1}) of \S\ref{prelimi}, is then
$$
 \bigl(\libc\bigl(\mathcal{U}\oplus s(V\otimes W)\bigr),d\bigr).
 $$
 Now, since $\catcdgl$ is proper, the pushout of  
 $$
 L\oplus L'\stackrel{\simeq}{\longleftarrow} (\libc(\mathcal{U}),d)\longrightarrow   \bigl(\libc\bigl(\mathcal{U}\oplus s(V\otimes W)\bigr),d\bigr)
 $$
induces a  quasi-isomorphism
\begin{equation}\label{quasi}
 \bigl(\libc\bigl(\mathcal{U}\oplus s(V\otimes W)\bigr),d\bigr)\stackrel{\simeq}{\longrightarrow}\bigl((L\oplus L')\amalg\libc\bigl(s(V\otimes W)\bigr),d\bigr)
\end{equation}
  where, in view of (\ref{ueue}), the codomain is isomorphic to
  $$
 (L\oplus L')\oplus \bigl(\libc\bigl((UL\otimes UL')\otimesc  s(V\otimes W)\bigr),d\bigr).
  $$
This proves the first statement.

For the second assertion, we write the model of $X\times Y$ as in Remark \ref{producto} where  $\mathcal{S}$ and $\mathcal{T}$ are now $U\oplus V$ and $U'\oplus W$. As in (\ref{producto2}), we have
$$
P\cong (\libc(\mathcal{S}),d)\oplus (\libc(\mathcal{T}),d) \oplus\bigl( \libc\bigl(UE\otimesc ([\mathcal{S},\mathcal{T}]\oplus s(\mathcal{S}\otimes\mathcal{T})\bigr),d\bigr)
$$
A direct computation  shows that, for the element $[v,w]\in[\mathcal{S},\mathcal{T}]$ the difference
$$
d_1[v,w]-  \sum_i\alpha_i\otimes [v_i, w]-\sum_j (-1)^{|v|+|v||\beta_j|}\,\beta_j\otimes [v, w_j]
$$
belongs to the submodule generated by $[U,W]$ and $[V,U']$. hence, in the codomain of the morphism
$$
P\longrightarrow  (\libc(\mathcal{U}\oplus s(V\otimes W),d)\stackrel{\simeq}{\longrightarrow}\bigl((L\oplus L')\amalg\libc\bigl(s(V\otimes W)\bigr),d\bigr),
$$ 
which is both, the composition of the model  of the map $X\times Y\to X\bar\wedge Y$  with  the quasi-isomorphism in (\ref{quasi}), and  a projection, it follows from (\ref{producto1}) that
$$
d_1s(v\otimes w)=-\sum_i\alpha_i\otimes s(v_i\otimes w)-\sum_j (-1)^{|v|+|v||\beta_j|}\,\beta_j\otimes s(v\otimes w_j).
$$
\end{proof}

 Assume now $X,Y\in\retb$ so that, with the notation in the previous result, $L=L'=\lasu_B^*$.
 
 \begin{proposition}\label{internal} 
 A retractive model of  $X\wedge_B Y\in\retb$ is given by
 $$
 L\oplus  \bigl(\libc\bigl((UL\otimes UL)\otimesc  s(V\otimes W)\bigr),d\bigr)
 $$
 where $L$ acts on the second summand via the diagonal adjoint action. 
 \end{proposition}

 \begin{proof}
 By Proposition \ref{external} the map $X\bar\wedge Y\to B\times B$ is model by the projection 
$$
(L\oplus L)\oplus \bigl(\libc\bigl((UL\otimes UL)\otimesc  s(V\otimes W)\bigr),d\bigr)\longrightarrow L\oplus L.
$$
In general, modelling does not commute with pullback but it does with inclusions of simplicial sets, see \cite[Prop.~7.8]{bufemutan0}. By (\ref{pullback1}) of \S\ref{prelimi}, a model of $X\wedge_BY$ is therefore obtained as the pullback of the above morphism along the diagonal $L\to L\oplus L$. This  is exactly the retractive stated in the proposition.
\end{proof}

 \begin{proof}[Proof of Theorem \ref{monoidal2}]
 
 By Proposition \ref{esfera}, $\Psi$ sends the identity object of $\Ho\spec_B$ to the identitiy object of $\Ho\cdgmul$, up to isomorphism.
 
 Next,
we show that for any $X,Y\in \spec_B$,
 $$
 \Psi(X\wedge_B Y)\cong\Psi(X)\otimesc \Psi(Y).
 $$

 Let 
 $$
 (L\amalg \libc(V),d)\quad\text{and}\quad (L\amalg \libc(W),d)
  $$
  be  models of the spectra $X$ and $Y$ respectively, which we may assume to be levelwise connected by (iii) of Remark \ref{casoconexo}. By Corollary \ref{corocomputable2}, and in the homotopy categories, $\calker^0$ applied to these $L$-spectra become isomorphic to 
  $$
  (UL\otimesc V,d_1) \quad\text{and}\quad (UL\otimesc W,d_1),
  $$
  where differentials and structure maps are induced by the indecomposable reductions of the given models.

   Therefore, 
   $$
   \Psi(X)\cong 
 \mathcal{D}(UL\otimesc V,d_1)=\varinjlim_p\,\, (UL\otimesc s^{-p+1}V^p,s^{-p+1}d_1)\cong UL\otimesc \bigl(\varinjlim_p\,s^{-p+1}V^p\bigr)
 $$
 and
 $$
 \Psi(Y)\cong  \mathcal{D}(UL\otimesc W,d_1)=\varinjlim_q\,\, (UL\otimesc s^{-q+1}W^q,s^{-q+1}d_1)\cong UL\otimesc \bigl(\varinjlim_q\,s^{-q+1}W^q).
 $$
A straightforward computation, using the explicit description of a directed colimit, yields
\begin{equation}\label{dormulafinal}
\begin{aligned}
\Psi(X)\otimesc \Psi(Y)&\cong \Bigl((UL\otimes UL)\otimesc\bigl(\oplus_{p,q}( s^{-p+1}V^p\otimes s^{-q+1}W^q)\bigr)\Bigl)/I\\
& 
\cong\Bigl((UL\otimes UL)\otimesc\bigl(\oplus_{n\ge 0}\,s^{-n+2}(\oplus_{p+q=n} V^p\otimes W^q)\bigr)\Bigl)/I.\\
\end{aligned}
\end{equation}
Here, the differential is induced by successive desuspension of the differential  in $(UL\otimes UL)\otimesc (V^p\otimes W^q)$ and  $I$ is explicitly described as follows:

The adjoint structure maps,
$$
\begin{aligned}
\eta_V^p\colon UL\otimesc V^p\to UL\otimesc s^{-1}V^{p+1},\\
\eta_W^q\colon UL\otimesc W^q\to UL\otimesc s^{-1}W^{q+1},
\end{aligned}
$$
induce a morphism
$$
\eta_{p,q}\colon (UL\otimes UL)\otimesc (V^p\otimes W^q)\to  (UL\otimes UL)\otimesc \bigl((s^{-1}V^{p+1}\otimes W^q)\oplus (V^p\otimes s^{-1}W^{q+1})\bigr).
$$
Then,  $I$ is generated by all successive desuspensions of $\{a-\eta_{p,q}(a)\}$ with $a\in  (UL\otimes UL)\otimesc (V^p\otimes W^q)$. 

In other words, consider the object of $\specmodu$,
$$
\bigl((UL\otimes UL)\otimesc (V\otimes W),d_1\bigr)=\{\bigl((UL\otimes UL)\otimesc (\oplus_{p+q=n} V^p\otimes W^q),d_1\bigr)\}_{n\ge0}
$$
whose $n$th adjoint structure map is $\eta_n=\oplus_{p+q=n}\eta_{p,q}$. Then, in view of (\ref{dormulafinal}) we see that
$$
\Psi(X)\otimesc \Psi(Y)\cong \mathcal{D}\Bigl(\mathbf{s}\bigl((UL\otimes UL)\otimesc (V\otimes W),d_1\bigr)\Bigr).
$$

On the other hand,  using  the explicit description  of the smash product of $B$-spectra from Remark \ref{monoidal} and  Proposition \ref{internal}, we deduce that $
\calker^0\lasub(X\wedge_BY)$ is weakly equivalent to the spectrum $R=\{R^n\}_{n\ge0}\in\specmodu$ where
$$
R^n= \bigl((UL\otimes UL)\otimesc \bigl(\oplus_{p+q=n}\, s(V^p\otimes W^q)\bigr),d_1\bigr).
$$
By Proposition \ref{external}, the differential  is given by
$$
d_1s(v\otimes w)=-\sum_i(\alpha_i\otimes 1)\otimes s(v_i\otimes w)-\sum_j (-1)^{|u|+|u||\beta_j|}\,(1\otimes\beta_j)\otimes s(v\otimes w_j).
$$
Furthermore, a direct computation shows that the $n$th structure map is precisely the suspension of $\eta_n$, modulo the general isomorphism $s\bigl(UL\otimes UL)\otimesc Z\cong (UL\otimes UL)\otimesc sZ$. Thus
$$
R\cong \mathbf{s}\bigl((UL\otimes UL)\otimesc (V\otimes W),d_1\bigr)
$$
and the claim follows.

\smallskip
  
  Finally, a careful yet routine computation  shows that in $\Ho\cdgmul$, the  isomorphisms 
$$
\Psi(X\wedge_B Y)\cong\Psi(X)\otimesc \Psi(Y)\quad\text{and}\quad \Psi(\bs_B)\cong\compul,
$$
satisfy the associative and unit coherence conditions.

 \end{proof}  
Combining the previous theorem and Proposition \ref{homotospec} we get:
\begin{corollary}\label{fffinal}
Given $X,Y\in\spec_B$,
$$
\pie(X\wedge_B Y)\otimes\bq\cong H\bigl(\Psi(X)\bigr)\otimesc H\bigl(\Psi(Y)\bigr).
$$
\hfill$\square$
\end{corollary}

\section{Appendix: spectra in (retractive) model categories}\label{apendice}

In this section, we collect the fundamental results on stabilization in model categories endowed with a Quillen pair of endofunctors, for which Hovey's foundational paper \cite{ho1} (and in some cases, its precursor work by Schwede \cite{schwe}) is the classical reference. These results are elegantly compiled and particularly readable in \S A.3 of Braunack-Mayer's doctoral thesis \cite{brau0}. We also include a brief subsection that gathers general results on (co)slice categories with model structures, along with another subsection that highlights key consequences of transferring a model structure through an adjunction.

In what follows, any statement without a proof or reference can be verified directly.

\subsection{Spectra on model categories}\label{specmod}

Since \cite{ho1}, the natural framework for developing spectra from the categorical point of view is within  a left proper combinatorial model category $\quic$ endowed with a Quillen pair of endofunctors\footnote{Caution: these are not the suspension and loop functors derived from the model structure, as they generally form a well-defined adjoint pair only in the homotopy category. Nevertheless, this notation may help unfamiliar readers follow the arguments and grasp the overall picture.}   $\Sigma\dashv\Omega$. We will write 
$(\quic,\sus\dashv\lup)$
 whenever we need to emphasize the particular Quillen pair. 
 
 \begin{definition}\label{especdtro} A {\em spectrum} in $\quic$ is a sequence $x=\{x_n\}_{n\ge 0}$ of objects of $\quic$ endowed with a sequence of {\em structure maps} $\sigma\colon\sus x_n\to x_{n+1}$, for $n\ge 0$, or equivalently, the corresponding adjoint maps $\sigma^\adj\colon x_n\to \lup x_{n+1}$. A morphism of spectra $f\colon x\to y$ is a collection of maps $\{f_n\}_{n\ge 0}$, $f_n\colon x_n\to y_n$, compatible with the structure maps in the obvious sense. We denote by 
 $\spec(\quic)$ 
 the category of spectra in $\quic$ assuming the corresponding pair of endofunctors is fixed.  This category is bicomplete \cite[Lemma 1.3]{ho1}.
 \end{definition}

 The pair $\sus\dashv\lup$ prolongs  to adjoint endofunctors $\suse\dashv\lupe$ in $\spec(\quic)$ by
$$
 (\suse x)_n=\sus x_n, \qquad (\lupe x)_n=\lup x_{n}
 $$
 with structure maps
 $$
 \sus(\suse x)_n=\sus^2 x_n\stackrel{\sus\sigma}{\longrightarrow}\sus x_{n+1}=(\suse x)_{n+1},\qquad (\lupe x)_n=\lup x_n\stackrel{\lup\sigma^\adj}{\longrightarrow}\lup^2x_{n+1}=\lup(\lupe x)_{n+1}.
 $$

 For each $k\ge0$, the evaluation functor
\begin{equation}\label{eva}
 \eva_k\colon \spec(\quic)\longrightarrow \quic,\quad x\mapsto x_k,
\end{equation}
 has a left adjoint
\begin{equation}\label{susk}
 \sus^{\infty-k}\colon \quic\longrightarrow\spec(\quic),\quad (\sus^{\infty-k}a)_n=\begin{cases} \sus^{n-k}a&\text{if $n\ge k$},\\ \,\,\,0&\text{otherwise},\end{cases}
\end{equation}
being $0$ the initial object of $\quic$.

 \begin{definition}\label{proestable}The {\em projective model structure} in $\quic$ is given by declaring a map $f\in\spec(\quic)$ to be a projective fibration or a projective weak equivalence if $f$ is a levelwise fibration or a levelwise weak equivalence. That is, if for each $n\ge 0$, the map $f_n$ is a  fibration or a weak equivalence, respectively. 
 \end{definition}
 
 Projective cofibrations and trivial cofibrations are characterized as follows, see \cite[Prop.~1.14]{ho1}: a map $f\colon x\to y$ of $\spec(\quic)$ is a projective cofibration if and only if $f_0$ and the induced maps $g_n\colon x_{n}\amalg_{\sus x_{n-1}}\sus y_{n-1}\to y_n$, $n\ge 1$, are cofibrations.  Analogously, $f$ is a projective trivial cofibration if and only if $f_0$ and $g_n$, $n\ge1 $ are trivial cofibrations in $\quic$.
 
 If $\quic$ is cofibrantly generated  by the set $ \mathcal I$ and $\mathcal J$ of cofibrations and trivial cofibrations respectively, then the projective model structure in $\spec(\quic)$ is also cofibrantly generated, see \cite[Thm.~1.13]{ho1},  by the  sets
\begin{equation}\label{generaes}
 {\mathcal I}_{\sus}=\cup_{k\ge 0}\sus^{\infty-k}({\mathcal I})\quad\text{and}\quad {\mathcal J}_{\sus}=\cup_{k\ge 0}\sus^{\infty-k}({\mathcal J}).
\end{equation}
 With respect to the projective model structure the adjoint pairs of functors $\suse\dashv\lupe$ and $\sus^{\infty-k}\dashv \eva_k$, $k\ge0$, are Quillen \cite[Prop.~1.15]{ho1}.

The {\em stable model structure} on $\spec(\quic)$ is constructed ad hoc to ensure that the adjunction $\suse\dashv\lupe$ is a Quillen equivalence.

To that end, we briefly recall the basics of Bousfield localization, see for instance \cite[Chap.~3]{hirsch0}:

\begin{definition}\label{boufieldlo} Let $\quic$ be a model category and  $S$ a set of morphisms of $\quic$. The {\em left Bousfield localization} of $\quic$ with respect to $S$ is a new model structure on $\quic$, denoted  $\quic_S$, together with a left Quillen functor $\quic\to\quic_S$ that is universal among left Quillen functors $F\colon \quic\to\quid$ for which $F(s)$ is a weak equivalence for every $s\in S$.

The  left Bousfield localization  with respect to any set of maps $S$ exists  provided that $\quic$ is left proper, and either cellular or combinatorial, see \cite[Thm.~4.1.1]{hirsch0} and \cite[Thm.~4.7]{bar} respectively.  The localized category $\quic_S$ retains the underlying category and cofibrations of $\quic$, while expanding the weak equivalences to include $S$-local equivalences. Fibrant objects in $\quic_S$ are precisely the $S$-local fibrant objects of $\quic$, and a map between such  objects  is a weak equivalence in $\quic_S$ if and only if it is one in $\quic$. Left properness and the combinatorial character is preserved by localization.
\end{definition}

Starting with a left proper combinatorial model category $\quic$ equipped with a Quillen pair of endofunctors $\sus\dashv\lup$, the category $\spec(\quic)$ is itself left proper and combinatorial, see for instance \cite[Lemma A.3.2]{brau0}. We may therefore  localize this category  with respect to the following set $S$:  for any object $a\in \quic$ and for any $n\ge 0$ write $\sus a =\eva_{n+1}\sus^{\infty-n} a$ so that  $\id_{\sus a}\colon\sus a\to \eva_{n+1}\sus^{\infty-n} a$ has an adjoint denoted by $\zeta_n^a\colon 
 \sus^{\infty-(n+1)} \sus a\to \sus^{\infty-n} a$. As in \cite[Def 3.3]{ho1} consider
\begin{equation}\label{equibous}
S=\{\zeta_n^{Qc}\}_{n\ge0,c}
\end{equation}
where $c$ runs through any domain or codomain of a set of generating cofibrations of $\quic$ and $Q$ denotes a functorial cofibrant replacement. 

\begin{definition} \label{eses} The {\em stable model structure} in $\spec(\quic)$ is its Bousfield localization with respect to $S$.  The fibrant objects in the stable structure are called {\em  $\lup$-spectra}. These are objects $x\in\spec(\quic)$ such that each $x_n$ is fibrant and, for all $n\ge 0$, the adjoint $x_n\to\lup x_{n+1}$ of the structure map is a weak equivalence  (see \cite[Thm.~3.4]{ho1}).
\end{definition}

As intended, with respect to the stable structure, the 
 pair $\suse\dashv\lupe$ is a Quillen equivalence in $\spec(\quic)$ \cite[Thm.~3.9]{ho1}. 
Moreover, if the endofunctors $\Sigma \dashv\lup$  already form a Quillen equivalence on $\quic$, then 
the adjunction
\begin{equation}\label{equiv}
 \xymatrix{ \quic & \spec(\quic)\ar@<0.75ex>[l]^(.54){\eva_0}
\ar@<0.75ex>[l];[]^(.45){\Sigma^{\infty}}\\}
\end{equation}
is also a Quillen equivalence when $\spec(\quic)$ is equipped with the stable structure.

\medskip

In what follows, given two pairs of adjoint functors $P\dashv Q$ and $P'\dashv Q'$ and a natural transformation $\eta\colon P\to P'$ (respec. $\eta\colon Q'\to Q$) we denote by $\eta=\tau^{\scriptscriptstyle\vee}\colon Q'\to Q$ (respec. $\tau=\eta^{\scriptscriptstyle\vee}\colon P\to P'$) its adjoint or dual natural transformation.

\begin{rem}\label{tauro} It is well known that $(\tau^{\scriptscriptstyle\vee})^{\scriptscriptstyle\vee}=\tau$ and that $\tau$ consists of isomorphisms if and only if $\eta$ does. As a consequence, recall from \cite[Cor.~1.4.1(b)]{ho0}, that if  adjunctions in question are Quillen pairs, then $\tau_a$ is a weak equivalence for any cofibrant object $a$ if and only if $\eta_b$ is a weak equivalence for any fibrant object $b$. Indeed, if $\tau$ has this property, it induces an isomorphism in the homotopy category. Consequently, by adjunction, $\eta$ induces an isomorphism in the homotopy category and is therefore a weak equivalence on fibrant objects.
 \end{rem}

 Let $(\quic,\sus\dashv\lup)$ and $(\quid,\sus'\dashv\lup')$ be two left proper, combinatorial model categories endowed with the corresponding pairs of Quillen endofunctors. Any given adjunction
$$ 
\xymatrix{ \quic & \quid \ar@<0.75ex>[l]^(.50){G}
\ar@<0.75ex>[l];[]^(.50){F}\\}
$$
naturally induces the adjoints pairs $F\sus\dashv \lup G$ and $\sus'F\dashv G\lup'$. Suppose further that  $F\dashv G$ is a Quillen pair  and let $\tau\colon F\sus\to \sus' F$ be a natural transformation such that $\tau_a$ is a weak equivalence for every cofibrant object $a\in \quic$. Equivalently, by Remark \ref{tauro} above, its dual transformation $\eta\colon G\Omega'\to \Omega G$ is a weak equivalence on fibrant objects. Under these assumptions,  Lemma 5.3, Proposition 5.5 and Theorem 5.7 of \cite{ho1} are summarized as follows:

    \begin{theorem}\label{prolon} {\em }  (i) There is a Quillen pair with respect to the stable model structures,
$$
\xymatrix{ \spec(\quic) & \spec(\quid) \ar@<0.75ex>[l]^(.50){\widetilde G}
\ar@<0.75ex>[l];[]^(.50){\widetilde F}\\}
$$
in which $\widetilde G$ is the prolongation  of $G$ and $\widetilde F\sus^{\infty-k}= {\sus^\prime}^{\infty-k} F$ for all $k$. Moreover, if $\tau_a$ is an isomorphism for any $a\in\quic$ (equivalently, $\eta_b$ is an isomorphism for any $b\in \quid$), then $\widetilde F$ is also the prolongation of $F$.

\medskip

 (ii) Additionally,  if $F\dashv G$ is a Quillen equivalence, and either the domains of the generating cofibrations are cofibrant, or else, $\tau_a$ is a weak equivalence for all $a\in\quic$ (equivalently, $\eta_b$ is a weak equivalence for all $b\in\quid$), then  $\widetilde F\dashv\widetilde G$ is also a Quillen equivalence.\hfill$\square$
\end{theorem}

Recall that the {\em prolongation} of $G$ is given by  $(\widetilde Gx)_n=Gx_n$ with (adjoint) structure maps 
$$
(\widetilde Gx)_n=Gx_n\stackrel{G\sigma^\adj}{\longrightarrow} G\lup'x_{n+1}\stackrel{\rho_{x_{n+1}}}{\longrightarrow}\lup G x_{n+1}=\lup(\widetilde G x)_{n+1},
$$
in which  $\rho=\tau^{\scriptscriptstyle\vee}\colon G\lup'\to \lup G$.

\begin{rem}\label{uyuy} Under the hypothesis of Theorem \ref{prolon}(i), we assume the weaker condition that $\tau_a$ is a weak equivalence for any $a\in \quic$, or alternatively, that the domains of the generating cofibrations are cofibrant. In this case $\widetilde F$ may not be the actual prolongation of $F$ but there exists a natural weak equivalence
$$
Fx_n\stackrel{\sim}{\longrightarrow}(\widetilde Fx)_n
$$
for all $n\ge 0$ and for any cofibrant spectra $x\in\spec(\quic)$. Furthermore, if $\widetilde\sigma$ denotes the $n$th structure map  of $\widetilde Fx$ and 
$$
\gamma\colon \Sigma'Fx_n\stackrel{\tau_{x_n}}{\longrightarrow}F\Sigma x_n\stackrel{F\sigma}{\longrightarrow}Fx_{n+1}
$$
is the structure map of the prolongation of $F$ applied to $x$, then the following diagram commutes in the homotopy category:
$$
\xymatrix{\Sigma'Fx_n\ar[d]_(.46)\gamma\ar[r]^\sim&\Sigma'(\widetilde Fx)_n\ar[d]^{\widetilde\sigma}\\
Fx_{n+1}\ar[r]_\sim&
(\widetilde Fx)_{n+1}.\\}
$$
In other words, the induced functor
$$
\widetilde F\colon \Ho\spec(\quic)\to \Ho\spec(\quid)
$$
is, up to isomorphism, naturally equivalent to  the functor induced by the prolongation of $F$. This is made explicit in the proof of \cite[Thm.~5.7]{ho1}.
\end{rem}

\subsection{Retractive  model categories and their spectra}\label{apret}

Let $\quic$ be a category and let $c\in\quic$. 

\begin{definition}\label{catret} The {\em retractive category over $c$} denoted by $\quic_{\sslash c}$ is the  category under $\id_c$ of the category $\quic$ over $c$. That is, an object of $\quicr$ consists of a map (retraction or projection) $x\to c$ of $\quic$ endowed with a section $c\to x$. A morphism $x\to y$ of $\quicr$ is encoded by a commutative diagram
$$
\xymatrix@R=10pt@C=10pt{
&c\ar[dl]\ar[dr]&\\
x\ar[dr]\ar[rr]&&y\ar[dl]\\
&c&
}
$$
where the diagonal arrows define $x$ and $y$ as objects of $\quicr$. 
\end{definition}

The zero object of this category is $c\stackrel{\id}{\to}c\stackrel{\id}{\to}c$ and if $\quic$ is locally presentable then so is $\quicr$.

Whenever $\quic$ has coproducts, the forgetful functor $\quicr\to\quic$ has a left adjoint which sends $x\in\quic$ to the coproduct $x\amalg c$ with the obvious section and retraction. We denote it by $x_{+c}$. 

Let 
\begin{equation}\label{ap1}
\xymatrix{ \quic & \quid \ar@<0.75ex>[l]^(.50){G}
\ar@<0.75ex>[l];[]^(.50){F}\\}
\end{equation} 
be a pair of adjoint functors, in which $\quic
$ has pullbacks, and let $c\in\quic$. Then, there is an induced pair of adjoint functors
\begin{equation}\label{ap2}
\xymatrix{ \quicr & \quid_{\sslash F(c)} \ar@<0.75ex>[l]^(.50){G_{\sslash c}}
\ar@<0.75ex>[l];[]^(.50){F_{\sslash c}}\\}
\end{equation} 
 where $F_{\sslash c}$ results from applying $F$ to the section and retraction of any object in $\quicr$ and $G_{\sslash c}$ applies first $G$ to the section and retraction of any object in $\quid_{\sslash F(c)}$,  and then pulls back the resulting retraction along the unit $c\to G\bigl(F(c)\bigr)$.

Dually, if $\quic$ has pushouts, any object $d\in\quid$ induces a pair of adjoint functors,
\begin{equation}\label{ap3}
\xymatrix{ \quic_{\sslash G(d)}& \quidr \ar@<0.75ex>[l]^(.50){G_{\sslash d}}
\ar@<0.75ex>[l];[]^(.50){F_{\sslash d}}\\}
\end{equation} 
in which $G_{\sslash d}$ simply applies $G$ to any object in $G_{\sslash d}$ and $F_{\sslash d}$ first applies $F$ to any object in $\quic_{\sslash G(d)}$, and then it pushes out the resulting section along the counit $F\bigl(G(d)\bigr)\to d$.

\begin{definition}\label{change} Let $\quic$ be a category with pullbacks and pushouts and let  $f\colon b\to c$ be a map of $\quic$. The {\em change of base adjunction} is the pair of adjoint functors,
\begin{equation}\label{ap4}
\xymatrix{ \quic_{\sslash b} & \quicr, \ar@<0.75ex>[l]^(.50){f^*}
\ar@<0.75ex>[l];[]^(.50){f_!}\\}
\end{equation} 
defined as follows:  $f_!$ pushes out $f$ along a section of a given object in $\quic_{\sslash b}$. On the other hand, $f^*$ pulls back $f$ along a retraction of a given object in $\quicr$.
 \end{definition}

\smallskip

Let now $\quic$ be a model category and let $c\in \quic$. Then, the retractive category $\quicr$ inherits a model structure in which a map is a fibration, cofibration or weak equivalence if its image through the forgetful functor $\quicr\to \quic$ is a fibration, cofibration or weak equivalence, see \cite[Thm.~7.6.5]{hirsch1} or \cite[Thm.~15.3.6]{maypon}.
It follows, as shown in \cite[Thms.~1.20,~1.24,~2.19,~2.24]{hirsch1} or \cite[Thm.~15.3.6 and Rem.~15.3.7]{maypon}, that if $\quic$ is left proper, right proper or cofibrantly generated,  then  $\quicr$ is as well. In fact, if $f\colon x\to y$ is either a generating cofibration or trivial cofibration of $\quic$ then,
\begin{equation}\label{gencof}
f_{+c}\colon x_{+c}\to y_{+c}
\end{equation}
 is a generating cofibration or trivial cofibration of $\quicr$. An immediate consequence is:
 
 \begin{corollary}\label{combret} 
 If $\quic$ is a proper combinatorial model category so is $\quicr$.\hfill$\square$
 \end{corollary}

If the adjoint functors in (\ref{ap1}) form a Quillen pair, then the ones in  (\ref{ap2}) and (\ref{ap3}) do as well \cite[Prop.~2.2(ii)]{li}. Also, by \cite[Prop.~2.2(i)]{li}, the change of base (\ref{ap4}) is always Quillen. Moreover, if $f$ is a weak equivalence and $\quic$ is  proper, then  (\ref{ap4}) is a Quillen equivalence. This is a particular instance of  \cite[Prop.~3.1(a)(i) and (ii)]{li} or \cite[Prop.~2.5 and Rem.~2.6]{rezk}

\medskip

As for the stabilization of retractive model categories,  let $\quic$ be   a left proper combinatorial model category endowed with a Quillen pair of endofunctors $\sus\dashv\lup$ and let $c\in\quic$. By the preceding discussion we have an induced Quillen pair of endofunctors $\sus_{\sslash c}\dashv\lup_{\sslash c}$ in $\quicr$, which we denote simply by $\sus_c\dashv\lup_c$. In view of  \S\ref{specmod} we introduce the following:  

\begin{definition}\label{specret} The data $(\quicr,\sus_c\dashv\lup_c)$ define a left proper combinatorial projective and stable model structures on  $\spec(\quicr)$ which will be denoted by  $\spec_c(\quic)$, or simply $\spec_c$ when there is no risk of confusion.
\end{definition}

Finally, let $f\colon c\to d$ be a map of $\quic$ for which there is a natural transformation $\tau\colon f_!\,\sus_c\to \sus_d \,f_!$ such that $\tau_a$ is a weak equivalence for any cofibrant object $a\in\quicr$. By Remark \ref{tauro}, this is equivalent to requiring that the dual natural transformation $\eta\colon f^*\Omega_d\to \Omega_c f^*$ is a weak equivalence for any fibrant object $b\in \quic_{\sslash d}$. Applying Theorem \ref{prolon} to the change of base adjunction (\ref{ap4}) we obtain:

\begin{corollary} \label{changeret} (i) There is a Quillen pair with respect to the stable model structures,
$$
\xymatrix{ \spec_c(\quic) & \spec_d(\quic) \ar@<0.75ex>[l]^(.50){\widetilde{f^*}}
\ar@<0.75ex>[l];[]^(.50){\widetilde{f_!}}\\}
$$
in which $\widetilde{f^*}$ is the prolongation of $f^*$ and $\widetilde{f_!}\,\sus_c^{\infty-k}= \sus_d^{\infty-k}\, f_!$ for all $k$. Moreover, if $\tau$ (equivalently $\eta$) consists of isomorphisms, then $\widetilde{f_!}$ is also the prolongation of $f_!$. 

\medskip

(ii) Additionally,  if $f$ is a weak equivalence, and either the domains of the generating cofibrations of $\quic$ are cofibrant, or else $\tau_a$ is a weak equivalence for all $a\in\quicr$ (equivalently, $\eta_b$ is a weak equivalence for all $b\in\quic_{\sslash d}$), then  $\widetilde{f_!}\dashv\widetilde{f^*}$ is a Quillen equivalence.\hfill${\square}$
\end{corollary}

\subsection{Some properties of transferred model structures}

Among various formulations of the well-known procedure of left or right transferring a model structure along an adjunction \cite[Thm.~11.3.2]{hirsch0} we choose the following:

\begin{theorem}\label{transfer}
{\em \cite[\S2.6]{bermoer}, \cite[Thm.~2.2.1]{hess}}
Let
$$
\xymatrix{
\quic &
\quid
\ar@<0.75ex>[l]^(.50){G}
\ar@<0.75ex>[l];[]^(.50){F}
}
$$
be a pair of adjoint functors between locally presentable categories. Then:

\smallskip

\emph{(Right transfer)}
If $\quic$ is a cofibrantly generated model category with generating cofibrations $\mathcal I$ and generating trivial cofibrations $\mathcal J$, and if $\quid$ admits a fibrant replacement functor together with functorial path objects for fibrant objects, then there exists a cofibrantly generated model structure on $\quid$ whose weak equivalences and fibrations are created by $G$. Moreover, this model structure is cofibrantly generated by $F(\mathcal I)$ and $F(\mathcal J)$.

\smallskip

\emph{(Left transfer)}
If $\quid$ is a cofibrantly generated model category, and if $\quic$ admits quasi-functorial cofibrant replacements together with very good cylinder objects for them, then there exists a cofibrantly generated model structure on $\quic$ whose weak equivalences and cofibrations are created by $F$.

\hfill$\square$
\end{theorem}

\begin{definition} Under the hypotheses of Theorem \ref{transfer}, we say that the resulting Quillen adjunction is a {\em right-transferred} or {\em left-transferred}  pair, respectively. 
\end{definition}

In this paper, we need the following consequences of Theorem \ref{transfer}.

\begin{proposition} \label{transferbous} Let \xymatrix{ \quic & \quid \ar@<0.55ex>[l]^(.50){G}
\ar@<0.55ex>[l];[]^(.50){F}\\} be a right-transferred pair and $S$ a set of maps in $\quic$ such that $G(a)$ is $S$-local for any $a\in\quid$. Then, the  pair 
\begin{equation}\label{necesito}\xymatrix{ \quic_S & \quid \ar@<0.55ex>[l]^(.50){G}
\ar@<0.55ex>[l];[]^(.50){F}\\}
\end{equation}
is also right-transferred.
\end{proposition}

\begin{proof} Trivially $F\colon \quic_S\to \quid$ preserves cofibrations. On the other hand, if $f$ is a fibration in $\quid$, $G(f)$ is a fibration in $\quic$ with $S$-local domain and codomain. Hence, by \cite[Thm.~3.3.16(1)]{hirsch0}, it is also a fibration in $\quic_S$ since its domain and codomain are $S$-local, and (\ref{necesito}) is indeed a Quillen pair.

On the other hand, if $G(f)$ is a fibration in $\quic_S$, then it is trivially a fibration in $\quic$ and thus $f$ is a fibration.
 
 Finally, recall that a weak equivalence in $\quic_S$ between local objects is a weak equivalence in $\quic$. Hence, if $G(f)$ is a weak equivalence in $\quic_S$ then it is also a weak equivalence in $\quic$, and consequently, $f$ itself is a weak equivalence.
\end{proof}
\begin{proposition} \label{transre} For any right-transferred pair \xymatrix{ \quic & \quid \ar@<0.55ex>[l]^(.50){G}
\ar@<0.55ex>[l];[]^(.50){F}\\} and  
 any $c\in\quic$, 
$$
\xymatrix{ \quicr & \quid_{\sslash F(c)} \ar@<0.75ex>[l]^(.50){G_{\sslash c}}
\ar@<0.75ex>[l];[]^(.50){F_{\sslash c}}\\}
$$
is also right-transferred.
\end{proposition}

\begin{proof} Simply observe that the cofibrations and trivial cofibrations of the transferred model structure in $\quid_{\sslash F(c)}$ are precisely the cofibrations and trivial cofibrations of the structure induced in $\quid_{\sslash F(c)}$ by that of $\quid$:

 Let $f$ be either a generating cofibration or trivial cofibration of  $\quic$. Then, in view of (\ref{gencof}), $f_{+c}$ is a generating cofibration or trivial cofibration of $\quicr$. On the one hand $F(f_{+c})$ is a generating cofibration or trivial cofibration of $\quid_{\sslash F(c)}$ with the structure transferred by $F_{\sslash c}\dashv G_{\sslash c}$. On the other hand, again by (\ref{gencof}), $F(f)_{+F(c)}$ is a generating cofibration or trivial cofibration of $\quid_{\sslash F(c)}$ with the structure induced by $\quid$. But, since $F$ preserves colimits, $F(f_{+c})=F(f)_{+F(c)}$. 
\end{proof} 

Finally, under the hypothesis and notation of Theorem \ref{prolon}, we have:

\begin{proposition}\label{trans} If \xymatrix{ \quic & \quid \ar@<0.55ex>[l]^(.50){G}
\ar@<0.55ex>[l];[]^(.50){F}\\} is right-transferred, the Quillen pair given in (i) of Theorem \ref{prolon},
$$
\xymatrix{ \spec(\quic) & \spec(\quid), \ar@<0.75ex>[l]^(.50){\widetilde G}
\ar@<0.75ex>[l];[]^(.50){\widetilde F}\\}
$$
is also right-transferred.
\end{proposition}

\begin{proof}
We check again that both model structures have the same cofibrations and trivial cofibrations. Let $\mathcal I$ and $\mathcal J$ be sets of generating cofibrations and trivial cofibrations of $\quic$. In view of (\ref{generaes}) a set of generating cofibrations  in the transferred model structure on $\spec(\quid)$ is 
$$
\widetilde F({\mathcal I}_\sus)=\cup_{k\ge 0}\widetilde F\bigl(\sus^{\infty-k}({\mathcal I})\bigr).
$$
By Theorem \ref{prolon}, given any $f\in \quic$ and $k\ge 0$, $\widetilde F\sus^{\infty-k} f={\sus'}^{\infty-k}Ff$. Hence,
  $$
 \widetilde F({\mathcal I}_\sus)=  \cup_{k\ge 0} {\sus'}^{\infty-k}F({\mathcal I})={\mathcal I}_{\sus'}.
  $$
  The last equality comes from the fact that $F\dashv G$ is a tranferred Quillen pair. Analogously 
  $$
 \widetilde F({\mathcal J}_\sus)={\mathcal J}_{\sus'}.
 $$
 Finally, recall from (\ref{generaes}) that ${\mathcal I}_{\sus'}$ and ${\mathcal J}_{\sus'}$ are generating sets of cofibrations and trivial cofibrations in the stable model structure on $\spec(\quid)$.
\end{proof}

Under the same premises, Propositions \ref{transre} and \ref{trans} readily imply:

\begin{corollary}\label{transcor} If \xymatrix{ \quic & \quid \ar@<0.55ex>[l]^(.50){G}
\ar@<0.55ex>[l];[]^(.50){F}\\} is right-transferred, then so is 
$$
\xymatrix{ \spec_c(\quic) & \spec_{F(c)}(\quid) \ar@<0.75ex>[l]^(.50){\widetilde{ G_{\sslash c}}}
\ar@<0.75ex>[l];[]^(.50){\widetilde{ F_{\sslash c}}}\\}
$$
for every $c\in\quic$.\hfill$\square$
\end{corollary}

\bigskip\bigskip\bigskip\bigskip

\noindent {\sc Institut de Math\'ematiques et Physique, Universit\'e Catholique de Louvain, Chemin du Cyclotron 2,
1348 Louvain-la-Neuve,
         Belgique}.

\noindent\texttt{yves.felix@uclouvain.be}

\bigskip

\noindent{\sc Departamento de \'Algebra, Geometr\'{\i}a y Topolog\'{\i}a, Facultad de Ciencias, Universidad de M\'alaga, Blvr. Louis Pasteur 31, 29010 M\'alaga, Spain.}

\noindent
\texttt{aniceto@uma.es}

\noindent
\texttt{alexsaiz221@gmail.com}


\begin{thebibliography}{99}

\bibitem{barmayriehl} T. Barthel, J.P. May and E. Riehl, \emph{Six model structures for DG-modules over DGAs: model category theory in homological action}, New York J. Math. {\bf 20} (2014), 1077--1159.

\bibitem{bar} C. Barwick, \emph{On left and right model categories and left and right Bousfield localizations}, Homology Homotopy Appl. {\bf 1} (2010), 1--76.

\bibitem{baber} M.~A. Batanin and  C. Berger, \emph{Homotopy theory for algebras over polynomial monads}, Theory Appl. Categ. {\bf 32} (2017), 148--253.
    
\bibitem{bor} F. Borceux, \emph{Handbook of Categorical Algebra 2},  Encyclopedia of Mathematics and its
Applications, Vol. 51, Cambridge University Press, Cambridge, 1994.


\bibitem{bermoer} C. Berger and I. Moerdijk, \emph{ Axiomatic homotopy theory for operads}, Comment. Math. Helv. {\bf 78} (2003), 805--831.
    
\bibitem{bouskan} A.~K. Bousfield and A.M. Kan, \emph{Homotopy limits, completions and
  localizations}, Lecture Notes in Mathematics, Vol. 304, Springer-Verlag,
  Berlin-New York, 1972.
  
  \bibitem{bous} A.~K. Bousfield, \emph{The localization of spectra with respect to homology}, Trans. Amer. Math. Soc. \textbf{148} (1979), 473--552.
  
   \bibitem{brau0}  V. Braunack-Mayer, \emph{Rational parametrised stable homotopy theory}, PhD thesis, Zurich University, 2018.
  
  \bibitem{brau1} V. Braunack-Mayer, \emph{Combinatorial parametrised spectra}, Algebr. Geom. Topol. {\bf 21} (2021), 801--891.
      
  \bibitem{brau2} V. Braunack-Mayer, \emph{Strict algebraic models for rational parametrised spectra I}, Algebr. Geom. Topol. {\bf 21} (2021), 917--1019.

\bibitem{bufemutan1} U. Buijs, Y. F\'elix, A. Murillo, and D. Tanr\'e,  \emph{{Homotopy theory of complete Lie algebras and Lie models of
  simplicial sets}}, J. Topol. {\bf 11} (2018), 799--825.

  \bibitem{bufemutan2} U. Buijs, Y. F\'elix, A. Murillo, and D. Tanr\'e, \emph{Lie models of simplicial sets and representability of the Quillen functor},  Isr. J. Math. {\bf 238} (2020) 313--358.

  \bibitem{bufemutan0} U. Buijs, Y. F\'elix, A. Murillo, and D. Tanr\'e,  \emph{Lie Models in Topology},   Progress in Mathematics \textbf{335},  Birkh\"auser-Springer Nature, 2021.
      
      \bibitem{crabbja} M. Crabb and I. James, \emph{Fibrewise homotopy theory}, Springer Monographs in Mathematics, Springer,   1998.

  \bibitem{fefuenmu0} Y. F\'elix, M. Fuentes and A. Murillo, \emph{All known realizations of complete Lie algebras coincide}, Algebr. Geom. Topol. \textbf{25} (2025),  1155--1167.
      
       \bibitem{fefuenmu1} Y. F\'elix, M. Fuentes and A. Murillo, \emph{A Lie characterization of the Bousfield-Kan $\bq$-completion and $\bq$-good spaces},   	arXiv:2407.02812.
      
   \bibitem{fehaltho}     Y. F\'elix, S. Halperin and J.-C. Thomas,
\emph{Rational Homotopy Theory}, Graduate Texts in Mathematics \textbf{205}, Springer, 2001.

\bibitem{femutan}  Y. F\'elix, A. Murillo, and D. Tanr\'e,  \emph{Fiberwise stable rational homotopy}, J. Topol. {\bf 3} (2010), 743--758.

  
    \bibitem{hess} K. Hess, M. Kedziorek, E. Riehl and B. Shipley,  \emph{A necessary and sufficient condition for induced model structures},  J. Topol. \textbf{10} (2017),  324--369.
    
  \bibitem{ho0} M. Hovey, \emph{Model Categories}, Mathematical Surveys and Monographs \textbf{63}, AMS, 1999.  
  
  \bibitem{ho1} M. Hovey, \emph{Spectra and symmetric spectra in general model categories}, J. Pure Appl.
Algebra {\bf 165} (2001), 63--127.

\bibitem{hoshis} M. Hovey, B. Shipley and J. Smith, \emph{Symmetric spectra}, J. Amer. Math. Soc. {\bf 13} (1999), 149--208.

    \bibitem{hirsch0} P.~S. Hirschhorn, \emph{Model Categories and Their Localizations}, Mathematical Surveys and Monographs {\bf 99}, AMS 2002.
    
 \bibitem{hirsch1} P.~S. Hirschhorn, \emph{Overcategories and undercategories of cofibrantly generated model categories},  J. Homotopy Relat. Struct. {\bf 16} (2021), 753--768.

\bibitem{lawsu} R. Lawrence and D. Sullivan,
\emph{A formula for topology/deformations
  and its significance}, Fund. Math. \textbf{225} (2014), 229--242.
    
 \bibitem{li}   Z. Li, \emph{A note on the model (co-)slice categories}, Chinese Annals of Mathematics, Series B {\bf 37}  (2016), 95--102.
    
    \bibitem{mal} C. Malkiewich, \emph{Parametrized spectra, a low-tech approach},  	arXiv:1906.04773.
        
   \bibitem{manmayschshi} M. Mandell, P. May, S. Schwede and  B. Shipley, \emph{Model categories of diagram spectra}, Proc. London Math. Soc. {\bf 82} (2001), 441--512.     
        
        \bibitem{maypon} J.~P. May and K. Ponto, \emph{More Concise Algebraic Topology}, Chicago Lectures in Mathematics, University of Chicago Press, Chicago, 2012. 
    
    \bibitem{maysi} J.~P. May and J. Sigurdsson, \emph{Parametrized Homotopy Theory}, Mathematical Surveys and Monographs {\bf 132}, American Mathematical Society, Providence, RI, 2006. 

\bibitem{qui} D. Quillen,
\emph{Rational Homotopy Theory}, Ann. of Math. \textbf{90} (1969), 205--295.

\bibitem{ree} C.~L. Reedy, \emph{Homotopy theory of model categories}, unpublished
manuscript (1974). 

\bibitem{rezk} C. Rezk, \emph{Every homotopy theory of simplicial algebras admits a proper model}, Topol. Appl. \textbf{119} (2002), 65--94.
    
  \bibitem{sch} U. Schreiber, \emph{Rational parametrized spectra}, 
MathOverflow, \url{https://mathoverflow.net/questions/261747} (2017).
    
    \bibitem{schwe} S. Schwede, \emph{Spectra in model categories and applications to the algebraic cotangent complex}, J. Pure
Appl. Algebra {\bf120}  (1997), 77--104.
    
    \bibitem{schship} S. Schwede and B. Shipley, \emph{Stable model categories are categories of modules}, Topology {\bf42}  (2003), 103--153.

    \bibitem{ship} B. Shipley, \emph{$H\bz$-Algebra spectra are differential graded algebras}, Amer. J. Math. {\bf 129} (2007), 351--379.



\bibitem{whi} D. White, \emph{Model structures on commutative monoids in general model
categories}, J. Pure
Appl. Algebra {\bf 221}  (2017), 3124--3168.
\end{thebibliography}
\end{document}